%% file: nonalg.tex
\documentclass[10pt,leqno]{amsart}

\newtheorem{lemma}{Lemma}
\newtheorem{theorem}{Theorem}
\newtheorem{corollary}{Corollary}
\newtheorem{proposition}{Proposition}
\newtheorem{definition}{Definition}

\usepackage{graphics}
\usepackage{epsfig}
\usepackage{amssymb,amsfonts,amscd}
\usepackage{enumerate}
\def\C{{\mathbb C}}
\def\K{{\mathbb K}}
\def\N{{\mathbb N}}
\def\R{{\mathbb R}}

\begin{document}


\title[Nonalgebraizable real analytic
tubes in $\C^n$]{Nonalgebraizable real analytic tubes in $\C^n$}

\author{Herv\'e Gaussier and Jo\"el Merker}

\address{(Gaussier and Merker) 
CNRS, Universit\'e de Provence, LATP, UMR 6632, CMI, 
39 rue Joliot-Curie, 13453 Marseille Cedex 13, France}

\email{[gaussier,merker]@cmi.univ-mrs.fr} 

\subjclass{Primary: 32V40. Secondary 32V25, 32H02, 32H40, 32V10}

\date{\number\year-\number\month-\number\day}

\begin{abstract}
We give necessary conditions for certain real analytic tube generic
submanifolds in $\C^n$ to be locally algebraizable. As an application,
we exhibit families of real analytic non locally algebraizable tube
generic submanifolds in $\C^n$. During the proof, we show that the
local CR automorphism group of a minimal, finitely nondegenerate real
algebraic generic submanifold is a real algebraic local Lie group.  We
may state one of the main results as follows. Let $M$ be a real
analytic hypersurface tube in $\C^n$ passing through the origin,
having a defining equation of the form $v=\varphi(y)$, where $(z,w)=
(x+iy,u+iv)\in\C^{n-1}\times \C$.  Assume that $M$ is Levi
nondegenerate at the origin and that the real Lie algebra of local
infinitesimal CR automorphisms of $M$ is of minimal possible dimension
$n$, {\it i.e.}  generated by the real parts of the holomorphic vector
fields $\partial_{z_1}, \dots, \partial_{z_{n-1}}, \partial_w$. Then
$M$ is locally algebraizable only if every second derivative
$\partial^2_{y_ky_l}\varphi$ is an algebraic function of the
collection of first derivatives
$\partial_{y_1}\varphi,...,\partial_{y_m}\varphi.$
\end{abstract}

\maketitle

\begin{center}
\begin{minipage}[t]{10cm}
\baselineskip =0.35cm
{\scriptsize

\centerline{\bf Table of contents~:}

\smallskip

{\bf \S1.~Introduction \dotfill 1.}

{\bf \S2.~Preliminaries \dotfill 5.}

{\bf \S3.~Proof of Theorem~1.1 \dotfill 8.}

{\bf \S4.~Local Lie group structure for the CR automorphism group \dotfill 17.}

{\bf \S5.~Minimality and finite nondegeneracy \dotfill 19.}

{\bf \S6.~Algebraicity of local CR automorphism groups \dotfill 22.}

{\bf \S7.~Description of explicit families of strong tubes in $\C^n$
\dotfill 28.}

{\bf \S8.~Analyticity versus algebraicity \dotfill 34.}

}\end{minipage}
\end{center}

\bigskip

\section*{\S1.~Introduction}

A real analytic submanifold $M$ in $\C^n$ is called {\it algebraic}\,
if it can be represented locally by the vanishing of a collection of
Nash algebraic real analytic functions. We say that $M$ is {\it
locally algebraizable}\, at one of its points $p$ if there exist some
local holomorphic coordinates centered at $p$ in which $M$ is
algebraic. For instance, every totally real, real analytic submanifold
in $\C^n$ of dimension $k\leq n$ is locally biholomorphic to a
$k$-dimensional linear real plane, hence locally algebraizable. Also,
every complex manifold is locally algebraizable.  Although every real
analytic submanifold $M$ is clearly locally equivalent to its tangent
plane by a {\it real}\, analytic (in general not holomorphic)
equivalence, the question whether $M$ is biholomorphically
equivalent to a real algebraic submanifold is subtle.  In this
article, we study the question whether every real analytic
CR submanifold is locally algebraizable.
One of the interests of
algebraizability lies in the reflection principle, 
which is better understood in the algebraic category. 
Indeed, in the fundamental works of Pinchuk [Pi1975], [Pi1978] and of
Webster [We1977], [We1978] and in the recent works of Sharipov-Sukhov
[SS1996], Huang-Ji [HJ1998], Verma [Ve1999], Coupet-Pinchuk-Sukhov
[CPS2000], and Shafikov [Sha2000], [Sha2002], the extendability of
germs of CR mappings with target in a real algebraic hypersurface is
achieved. On the contrary, even if some results previously shown under
an algebraization hypothesis were proved recently under general
assumptions (see the strong result obtained by Diederich-Pinchuk
[DP2003]), most of the results cited above are still open in the case
of a real analytic target hypersurface.

\subsection*{1.1.~Brief history of the question}
By the work of Moser and Webster [MW1983, Thm.~1], it is known that every
real analytic two-dimensional surface $S\subset\C^2$ at an isolated
elliptic (in the sense of Bishop) complex tangency $p\in S$ is
biholomorphic to one of the surfaces $S_{\gamma,\delta,s}:=\{
(z_1,z_2)\in\C^2:\, y_2=0,\, x_2=z_1\bar z_1+
(\gamma+\delta(x_2)^s)(z_1^2+\bar z_1^2)\}$, where $p$ corresponds to
the origin, where $0<\gamma<1/2$ is Bishop's invariant and where
$\delta=\pm 1$ and $s\in\N$ or $\delta=0$. The quantities $\gamma$,
$\delta$, $s$ form a complete system of biholomorphic invariants for
the surface $S$ near $p$. In particular, every elliptic surface
$S\subset \C^2$ is locally algebraizable. 
To the authors' knowledge,
it is unknown whether there exist nonalgebraizable hyperbolic surfaces
in $\C^2$. In fact, 
very few examples of nonalgebraizable submanifolds
are known.  In [Eb1996], the author constructed a nonminimal (and non
Levi-flat) real analytic hypersurface $M$ through the origin in $\C^2$
which is not locally algebraizable ({\it cf.}
[BER2000, p.~330]).  In a recent article [HJY] the authors prove that
the strongly pseudoconvex real analytic hypersurface ${\rm Im}\, w=
e^{\vert z\vert^2}-1$ passing through the origin in $\C^2$ is not
locally algebraizable at any of its points.  Using an associated
projective structure bundle $\mathcal{Y}$ introduced by Chern, they
show that for every rigid algebraic hypersurface in $\C^n$, there
exists an algebraic dependence relation between seven explicit
Cartan-type holomorphic invariant functions on $\mathcal{Y}$. However
a computational approach shows that when $M$ is of the specific form
${\rm Im}\, w= e^{\vert z\vert^2}-1$, no algebraic relation can be
satisfied by these seven invariants.

\subsection*{1.2.~Presentation of the main results}
Our aim is to present a geometrical approach of the problem, valid in
arbitrary dimension and in arbitrary codimension, and to exhibit a
large class of nonalgebraizable real analytic generic submanifolds. We
consider the class $\mathcal{T}_n^d$ of generic real analytic
submanifolds in $\C^n$ passing through the origin, of codimension
$d\geq 1$ and of CR dimension $m=n-d\geq 1$, whose local CR
automorphism group is $n$-dimensional, generated by the real parts of
$n$ holomorphic vector fields having holomorphic coefficients
$X_1,\dots,X_n$ which are linearly independent at the origin and which
commute: $[X_{i_1},X_{i_2}]=0$. We shall call $\mathcal{T}_n^d$ the
class of {\it strong tubes} of codimension $d$. Indeed, since there
exists a straightened system of coordinates $t=(t_1,\dots,t_n)$ over
$\C^n$ in which $X_i=\partial_{t_i}$, we observe that every
submanifold $M\in\mathcal{T}_n^d$ is tubifiable at the origin. By
this, we mean that there exist holomorphic coordinates $t=(z,w)=
(x+iy,u+iv)\in\C^{m}\times\C^d$ vanishing at the origin in which $M$
is represented by $d$ equations of the form $v_j=\varphi_j(y)$. Hence
$M$ is a tube, {\it i.e.} a product of the submanifold
$\{v_j=\varphi_j(y),\ j=1,\dots,d\}\subset \R_{y,v}^n$ by the
$n$-dimensional real space $\R_{x,u}^n$. Since $M\in\mathcal{T}_n^d$,
the only infinitesimal CR automorphisms of $M$ are the real parts of
the vector fields
$\partial_{z_1},\dots,\partial_{z_{m}},\partial_{w_1},\dots,\partial_{w_d}$,
explaining the terminology.  Notice that not every tube belongs to the
class $\mathcal{T}_n^d$. For instance in codimension $d=1$, the
Heisenberg sphere $v=\sum_{k=1}^{n-1}\,y_k^2$ and more generally the
Levi nondegenerate quadrics $v=\sum_{k=1}^{n-1}\,\varepsilon_k\,
y_k^2$, where $\varepsilon_k=\pm 1$, have a CR automorphism group of
dimension $(n+1)^ 2-1>n$ and so do not belong to $\mathcal{T}_n^1$. We
assume that $M\in\mathcal{T}_n^d$ is minimal at the origin, namely the
local CR orbit of $0$ in $M$ contains a neighborhood of $0$ in
$M$. Furthermore, we assume that $M\in\mathcal{T}_n^d$ is finitely
nondegenerate at $0$, namely that there exists an integer $\ell\geq 1$
such that ${\rm Span}\, \{\overline{L}^\beta\, \nabla_t (r_j)(0,0): \,
\beta \in \mathbb N^m, \,\vert \beta\vert \leq \ell, \,
j=1,\dots,d\}=\C^n$, where $r_j(t,\bar t)=0$, $j=1,\dots,d$ are
arbitrary real analytic defining functions for $M$ near $0$ satisfying
$\partial r_1\wedge \cdots\wedge \partial r_d\neq 0$ on $M$, where
$\nabla_t (r_j)(t,\bar t)$ is the holomorphic gradient with respect to
$t$ of $r_j$ and where $\overline{L}^\beta$ denotes
$(\overline{L}_1)^{\beta_1}\cdots (\overline{L}_m)^{\beta_m}$ for an
arbitrary basis $\overline{L}_1,\dots,\overline{L}_m$ of
$(0,1)$-vector fields tangent to $M$ in a neighborhood of $0$.  In
particular Levi nondegenerate hypersurfaces are finitely
nondegenerate.  Finally, assuming only that $\varphi_j(0)=0$,
$j=1,\dots,d$, we shall observe in Lemma~3.2 below that a tube
$v_j=\varphi_j(y)$ of codimension $d$ is finitely nondegenerate at the
origin if and only if there exist multi-indices
$\beta_*^1,\dots,\beta_*^m\in\N^m$ with $\vert\beta_*^k\vert\geq 1$
and integers $1\leq j_*^1,\dots,j_*^m\leq d$ such that the real
mapping
\def\theequation{1.1}\begin{equation}
\psi(y):=\left(
{\partial^{\vert \beta_*^1\vert}\varphi_{j_*^1}(y)
\over \partial y^{\beta_*^1}},\ldots,
{\partial^{\vert \beta_*^m\vert}\varphi_{j_*^m}(y)
\over \partial y^{\beta_*^m}}
\right)=:y'\in\R^m
\end{equation}
is of rank $m$ at the origin in $\R^m$. Our main theorem provides a
necessary condition for the local algebraizability of strong tubes :

\def\thetheorem{1.1}\begin{theorem} 
Let $M$ be a real analytic generic tube of codimension $d$ in $\C^n$
given in coordinates $(z,w)=(x+iy,u+iv)\in\C^m\times\C^d$ by the
equations $v_j=\varphi_j(y)$, where $\varphi_j(0)=0$, $j=1,\dots,d$.
Assume that $M$ is minimal and finitely nondegenerate at the origin,
so the real mapping $\psi(y)=y'$ defined by~\thetag{1.1} is of rank
$m$ at the origin in $\R_y^m$ and let $y=\psi'(y')$ denote the local inverse
in $\psi(y)$.  Assume that $M\in\mathcal{T}_n^d$, namely $M$ is a
strong tube of codimension $d$. If $M$ is locally algebraizable at the
origin, then all the derivative functions $\partial_{y_k'}
\psi_l'(y')$, where $1\leq k,l\leq m$, are real algebraic
functions of $y'$. Equivalently, every second derivative
$\partial^2_{y_ky_l}\varphi_j(y)$ is an algebraic function of the
collection of first derivatives
$\partial_{y_1}\varphi_j,\dots,\partial_{y_m}\varphi_j.$
\end{theorem}

By contraposition, every real analytic strong tube
$M\in\mathcal{T}_n^d$ for which one of the derivative functions
$\partial_{y_k'}\psi_l'$ is not real algebraic is not locally
algebraizable. We will argue in \S8 that this is generically the case
in the sense of Baire.  It is however natural to look for explicit
examples of nonalgebraizable real analytic submanifolds in $\mathbb
C^n$. Since the real parts of the vector fields $\partial_{z_1},
\dots, \partial_{z_{m}}, \partial_{w_1},\dots,\partial_{w_d}$ are
infinitesimal CR automorphisms of every tube $v=\varphi(y)$, we must
provide some sufficient conditions insuring that the dimension of the
Lie algebra of such a tube is exactly $n$. We shall establish in
\S\S7-8 below:

\def\thecorollary{1.2}\begin{corollary}
The tube hypersurface $M_{\chi_1,\dots,\chi_{n-1}}$ in $\C^n$ of equation
$v=\sum_{k=1}^{n-1}[\varepsilon_ky_k^2+y_k^6+y_k^9y_1 \cdots
y_{k-1}+y_k^{n+8}\chi_k(y_1,\dots,y_{n-1})]$, where
$\chi_1,\dots,\chi_{n-1}$ are arbitrary real analytic functions, belongs
to the class $\mathcal{T}_n^1$ of strong tubes. Two such tubes
$M_{\chi_1,\dots,\chi_{n-1}}$ and
$M_{\widehat{\chi}_1,\dots,\widehat{\chi}_{n-1}}$ are biholomorphically
equivalent if and only if $\chi_j=\widehat{\chi}_j$ for every
$j$. Furthermore, for a generic choice in $\chi_1,\dots,\chi_{n-1}$ in
the sense of Baire {\rm (to be precised in \S8)},
$M_{\chi_1,\dots,\chi_{n-1}}$ is not locally algebraizable at the
origin.
\end{corollary}

Here we annihilate some Taylor coefficients in $\varphi$ and keep
some others to be nonzero to insure that $M_\chi$ is a strong tube.
Furthermore, the terms $y_k^9y_1\cdots y_{k-1}$ insure that the
$M_\chi$ are pairwise not biholomorphically equivalent. Using a
classical direct algorithm ({\it cf.}~[Bs1991], [St1991]), or the Lie
theory of symmetries of differential equations, combined with
Theorem~1.1 we may provide some other explicit strong tubes which are
not locally algebraizable (see \S\S7-8 for the proof):

\def\thecorollary{1.3}\begin{corollary} The following five explicit
tubes belong to $\mathcal{T}^1_2$ and are not locally algebraizable at
the origin\,{\rm :} $v=\sin(y^2)$, $v=\tan(y^2)$, $v=e^{e^y-1}-1$,
$v=\sinh (y^2)$ and $v=\tanh(y^2)$.
\end{corollary}

In these five examples, the algebraic independence in
$\partial_y\varphi$ and in $\partial^2_{yy}\varphi$ is clear; however,
checking that each hypersurface is indeed a strong tube requires some
formal computations, {\it see}~\S7.  One may also check by a direct
computation that in a neighborhood of every point $p=(z_p,w_p)$ with
$z_p\neq 0$, the hypersurface $M_{\rm HJY}$ of global equation ${\rm
Im}\, w=e^{\vert z\vert^2}-1$ is a strong tube ({\it see}~\S7.5).
Since it can be represented in a neighborhood of $p$ under the tube
form $v'=e^{\vert z_p\vert^2(e^{y'}-1)}-1$ by means of the local
change of coordinates $z'=2i\, \ln(z/z_p)$, $w'=(w-w_p)\, e^{-\vert
z_p\vert^2}$, applying Theorem~1.1 and inspecting the function
$e^{\vert z_p\vert^2(e^{y'}-1)}-1$, we may check that it is {\it
not}\, algebraizable at such points $p$ with $z_p\neq 0$ ({\it
see}~\S7.5). It follows trivially that the hypersurface $M_{\rm HJY}$
is also not locally algebraizable at all the points $p$ with $z_p=0$,
giving the result of [HJY, Theorem~1.1]. Using the same strategy as
for Theorem~1.1, we obtain more generally the following criterion:

\def\thetheorem{1.4}\begin{theorem}
Let $M_\varphi:v=\varphi(z\bar z)$ be a Levi nondegenerate real
analytic hypersurface in $\C^2$ passing through the origin whose Lie
algebra of local infinitesimal CR automorphisms is generated by
$\partial_w$ and $iz\, \partial_z$. If $M_\varphi$ is locally
algebraizable at the origin, then the first derivative
$\partial_r\varphi$ in $\varphi$ $(r\in \R)$ is algebraic. For
instance, the following seven explicit examples are not locally
algebraizable at the origin\,{\rm :} $v=e^{z\bar z}-1$, $v=\sin(z\bar
z)$, $v=\tan(z\bar z)$, $v=\sinh(z\bar z)$, $v=\tanh(z\bar z)$,
$v=\sin(\sin(z\bar z))$ and $v=e^{e^{z\bar z}-1}-1$.
\end{theorem}

Finally, using the same recipe as for Theorems~1.1 and~1.4, we shall
provide a very simple criterion for the local nonalgebraizability of
some hypersurfaces having a local Lie CR automorphism group of
dimension equal to one exactly.  We consider the class $\mathcal{R}_n$
of Levi nondegenerate real analytic hypersurfaces passing through the
origin in $\C^n$ ($n\geq 2$) such that the Lie algebra of
infinitesimal CR automorphisms of $M$ is generated by exactly one
holomorphic vector field $X_1$ with holomorphic coefficients not all
vanishing at the origin. We call $\mathcal{R}_n$ the class of {\it
strongly rigid}\, hypersurfaces, in order to distinguish them from the
so-called {\it rigid}\, ones whose local CR automorphism group may be
of dimension larger than $1$.  By straightening $X_1$, we may assume
that $X_1=\partial_w$ and that $M$ is given by a real analytic
equation of the form $v=\varphi(z,\bar z)=\varphi(z_1,\dots,z_{n-1},
\bar z_1,\dots,\bar z_{n-1})$.  By making some elementary changes of
coordinates ({\it cf.}~\S3.3), we can furthermore assume without loss
of generality that $\varphi(z,\bar z)=\sum_{k=1}^{n-1}\,
\varepsilon_k\, \vert z_k\vert^2+\chi(z,\bar z)$, where
$\varepsilon_k=\pm 1$ and $\chi(0,\bar z)\equiv
\partial_{z_k}\chi(0,\bar z)\equiv 0$.

\def\thetheorem{1.5}\begin{theorem}
Let $M: v=\varphi(z,\bar z)=\sum_{k=1}^{n-1}\, 
\varepsilon_k\,\vert z_k\vert^2+\chi(z,\bar z)$ be a
strongly rigid hypersurface in $\C^n$ with $\chi(0,\bar z)\equiv 
\partial_{z_k}\chi(0,\bar z)\equiv 0$. If $M$ is locally 
algebraizable at the origin, then all the first derivatives 
$\partial_{z_k}\varphi$ are algebraic functions of 
$(z,\bar z)$.
\end{theorem}

This criterion enables us to exhibit a whole family of non locally
algebraizable hypersurfaces in $\C^n$ :

\def\thecorollary{1.6}\begin{corollary}
The rigid hypersurfaces $M_{\chi_1,\dots,\chi_{n-1}}$ in
$\C^n$ of equation $v= \sum_{k=1}^{n-1}\, [\varepsilon_k\, \vert
z_k\vert^2+\vert z_k\vert^{10}+
\vert z_k\vert^{14}+
\vert z_k\vert^{16}(z_k+\bar z_k)+
\vert z_k\vert^{18} \vert z_1\vert^2 \cdots \vert z_{k-1}\vert^2
+\vert z_k\vert^{2n+16}\, \chi_k(z,\bar z)]$, where the
$\chi_k$ are arbitrary real analytic functions, belong to the class
$\mathcal{R}_n$ of strongly rigid hypersurfaces. Two such tubes
$M_{\chi_1, \dots, \chi_{n-1}}$ and $M_{\widehat{\chi}_1, \dots,
\widehat{\chi}_{n-1}}$ are biholomorphically equivalent if and only if
$\chi_k= \widehat{\chi}_k$ for $k=1,\dots,n-1$.  Furthermore, for a
generic choice of a $(n-1)$-tuple of real analytic functions
$(\chi_1,\dots,\chi_{n-1})$ in the sense of Baire {\rm (to be precised
in \S8)}, $M_{\chi_1,\dots,\chi_{n-1}}$ is not locally algebraizable
at the origin.
\end{corollary}

Finally, by computing generators of the Lie algebra of local 
infinitesimal CR automorphisms of some explicit examples, we
obtain:

\def\thecorollary{1.7}\begin{corollary}
The following seven explicit examples of hypersurfaces in $\C^2$ are
{\rm strongly rigid} and are not locally algebraizable at the
origin\,{\rm :} $v=z\bar z+z^2\bar z^2\sin(z+\bar z)$, $v=z\bar
z+z^2\bar z^2\exp(z+\bar z)$, $v=z\bar z+z^2\bar z^2\cos(z+\bar z)$,
$v=z\bar z+z^2\bar z^2\tan(z+\bar z)$, $v=z\bar z+z^2\bar
z^2\sinh(z+\bar z)$, $v=z\bar z+z^2\bar z^2\cosh(z+\bar z)$ and
$v=z\bar z+z^2\bar z^2\tanh(z+\bar z)$.
\end{corollary}

\subsection*{1.3.~Content of the paper}
To prove Theorem~1.1 we consider an algebraic equivalent $M'$ of
$M$. The main technical part of the proof consists in showing that an
arbitrary real algebraic element $M'$ of $\mathcal{T}_n^d$ can be
straightened in some local complex {\it algebraic}\, coordinates
$t'\in\C^n$ in order that its infinitesimal CR automorphisms are the
real parts of $n$ holomorphic vector fields of the form
$X_i'=c_i'(t_i')\, \partial_{t_i'}$, $i=1,\dots,n$, where the
variables are separated and the functions $c_i'(t_i')$ are
algebraic. For this, we need to show that the automorphism group of a
minimal finitely nondegenerate real algebraic generic submanifold in
$\C^n$ is a local real algebraic Lie group, a notion defined in
\S2.3. A large part of this article (\S\S4, 5, 6) is devoted to
provide an explicit representation formula for the local biholomorphic
self-transformations of a minimal finitely nondegenerate generic
submanifold, {\it see}\, especially Theorem~2.1 and Theorem~4.1.
Finally, using the specific simplified form of the vector fields
$X_i'$ and assuming that there exists a biholomorphic equivalence
$\Phi: M\to M'$ satisfying $\Phi_*(\partial_{t_i})= c_i'(t_i')\,
\partial_{t_i'}$, we show by elementary computations that all the
first order derivatives of the mapping $\psi'(y')$ must be algebraic.
We follow a similar strategy for the proofs of Theorems~1.4 and~1.5.
Finally, in \S\S7-8, we provide the proofs of Corollaries~1.2, 1.3,
1.6 and~1.7.

\subsection*{1.4.~Acknowledgment} We acknowledge interesting discussions
with Michel Petitot Fran\c cois Boulier at the University of Lille 1.

\section*{\S2.~Preliminaries} We recall in this section the basic
properties of the objects we will deal with.

\subsection*{2.1.~Nash algebraic functions and manifolds}
In this subsection, let $\K=\R$ or $\C$. Let $(x_1,\dots,x_n)$ denote
coordinates over $\K^n$. Throughout the article, we shall use the norm
$\vert x \vert:= \max(\vert x_1\vert,\dots,\vert x_n\vert)$ for
$x\in\K^n$. Let $\mathcal{K}$ be an open polydisc centered at the
origin in $\K^n$, namely $\mathcal{K}=\{x\in\K^n: \vert x\vert
<\rho\}$ for some $\rho>0$.  Let $f: \mathcal{K}\to \K$ be a
$\K$-analytic function, defined by a power series converging normally
in $\mathcal{K}$. We say that $f$ is (Nash) {\it $\K$-algebraic}\, if
there exists a nonzero polynomial $P(X_1,\dots,X_n,F)\in\K[X_1,\dots,
X_n,F]$ in $(n+1)$ variables such that the relation
$P(x_1,\dots,x_n,f(x_1,\dots,x_n))=0$ holds for all
$(x_1,\dots,x_n)\in \mathcal{K}$. If $\K=\R$, we say that $f$ is {\it
real algebraic}. If $\K=\C$, we say that $f$ is {\it complex
algebraic}. The category of $\K$-algebraic functions is stable under
elementary algebraic operations, under differentiation and under
composition. Furthermore, implicit solutions of $\K$-algebraic
equations (for which the real analytic implicit function theorem
applies) are again $\K$-algebraic mappings.  The theory of
$\K$-algebraic manifolds is then defined by the usual axioms of
manifolds, for which the authorized changes of chart are
$\K$-algebraic mappings only ({\it cf.} [Za1995]). In this paper, we
shall very often use the stability of algebraicity under
differentiation.

\subsection*{2.2.~Infinitesimal CR automorphisms}
Let $M\subset \C^n$ be a generic submanifold of codimension $d\geq 1$
and CR dimension $m=n-d\geq 1$. Let $p\in M$, let $t=(t_1,\dots,t_n)$
be some holomorphic coordinates vanishing at $p$ and for some
$\rho>0$, let $\Delta_n(\rho):=\{t\in \C^n: \vert t\vert < \rho\}$ be
an open polydisc centered at $p$. We consider the Lie algebra
$\mathfrak{Hol}(\Delta_n(\rho))$ of holomorphic vector fields of the
form $X=\sum_{j=1}^n\, a_j(t)\,\partial/\partial t_j$, where the $a_j$
are holomorphic functions in $\Delta_n(\rho)$. Here,
$\mathfrak{Hol}(\Delta_n(\rho))$ is equipped with the usual Jacobi-Lie
bracket operation. We may consider the complex flow $\exp(\sigma
X)(q)$ of a vector field $X\in\mathfrak{Hol}(\Delta_n(\rho))$. 
It is a holomorphic
map of the variables $(\sigma,q)$ which is well defined in some
connected open neighborhood of $\{0\}\times \Delta_n(\rho)$ in $\C\times 
\Delta_n(\rho)$.

Let $K$ denote the real vector field $K:=X+\overline{X}$, considered
as a {\it real}\, vector field over $\R^{2n}\cong \C^n$. Again, the
real flow of $K$ is defined in some connected open neighborhood of
$\{0\}\times \Delta_n(\rho)^\R$ in $\R\times \Delta_n(\rho)^\R$.  We
remind the following elementary relation between the flow of $K$ and
the flow of $X$.  For a real time parameter $\sigma:=s\in \R$, the
flow $\exp(sX)(q)$ coincides with the real flow of $X+\overline{X}$,
namely $\exp(sX)(q)=\exp(s(X+\overline{X}))(q^\R)$, where for
$q\in\C^n$, we denote $q^\R$ the corresponding real point in
$\R^{2n}$. In the sequel, we shall always identify $\Delta_n(\rho)$
and its real counterpart $\Delta_n(\rho)^\R$.

Let now $\mathfrak{Hol} (M,\Delta_n(\rho))$ denote the real subalgebra 
of the vector fields $X\in \mathfrak{Hol}(\Delta_n(\rho))$
such that $X+\overline{X}$ is tangent to $M\cap \Delta_n(\rho)$. We
also denote by $\mathfrak{Aut}_{CR} (M,\Delta_n(\rho))$ the Lie
algebra of vector fields of the form $X+ \overline{X}$, where $X$
belongs to $\mathfrak{Hol}(M, \Delta_n(\rho))$, so
$\mathfrak{Aut}_{CR} (M,\Delta_n(\rho)) = 2\,{\rm Re}\,
\mathfrak{Hol}(M,\Delta_n(\rho))$. By the above considerations, the
local flow $\exp(sX)(q)$ of $X$ with $s\in \R$ real makes a
one-parameter family of local biholomorphic transformations of $M$. In
the sequel, we shall always identify
$\mathfrak{Hol} (M,\Delta_n(\rho))$ and
$\mathfrak{Aut}_{CR} (M,\Delta_n(\rho))$, namely we shall identify $X$
and $X+\overline{X}$ and say by some abuse of language that $X$ itself
is an infinitesimal CR automorphism.

In the algebraic category, the main drawback of infinitesimal CR
automorphism is that they do not have algebraic flow.  For instance,
the complex dilatation vector field $X=iz\,\partial_z$ has transcendent
flow, even if it is an infinitesimal CR automorphism of every algebraic
hypersurface in $\C^2$ whose equation is of the form $v=
\varphi(z\bar z)$, even if the coefficient of $X$ is
algebraic. Thus instead of infinitesimal CR automorphisms which
generate one-parameter groups of biholomorphic transformations of $M$,
we shall study algebraically dependent one-parameter
families of biholomorphic transformations (not necessarily making a
one parameter group).  To begin with, we need to introduce some
precise definitions about local algebraic Lie transformation groups.

\subsection*{2.3.~Local Lie group actions in the $\K$-algebraic category}
Often in real or in complex analytic geometry, the interest cannot be
focalized on global Lie transformation groups, but only on local
transformations which are close to the identity.  For instance, the
transformation group of a small piece of a real analytic CR manifold
in $\C^n$ which is not contained in a global, large or compact
CR manifold is almost never a true, global transformation
group. Consequently the usual axioms of Lie transformation groups must
be localized. Philosophically speaking, the local point of view is
often the most adequate and the richest one, because a given
analytico-geometric object often possesses much more local invariant
than global invariants, if any. Historically speaking, the local Lie
transformation groups were first studied, before the introduction of
the now classical notion of global Lie group. Especially, in his
first masterpiece work [Lie1880] on the subject, Sophus Lie
essentially dealt with local ``Lie'' groups: he classified all
continuous local transformation groups acting on an open
subset of $\C^2$. This general classification provided afterwards in
the years 1880--1890 many applications to the local study of
differential equations: local normal forms, local solvability, {\it
etc.}

In this paragraph we define precisely local actions of local Lie groups 
and we focus especially on the $\K$-algebraic
category.

Let $c\in\N_*$, let $g=(g_1,\dots,g_c)\in\K^c$ and let two positive
numbers satisfy $0<\delta_2<\delta_1$.  We formulate the desired
definition by means of the two precise polydiscs
$\Delta_c(\delta_2)\subset\Delta_c(\delta_1)\subset \K^c$.  A {\it local
$\K$-algebraic Lie group of dimension $c$}\, consists of the following
data:
\begin{itemize}
\item[{\bf (1)}]
A $\K$-algebraic {\it multiplication mapping}\, $\mu: \Delta_c(\delta_2)\times
\Delta_c(\delta_2)\to \Delta_c(\delta_1)$ which is locally
associative ($\mu(g,\mu(g',g''))=\mu(\mu(g,g'),g'')$), whenever
$\mu(g',g'')\in\Delta_c(\delta_2)$, $\mu(g,g')\in\Delta_c(\delta_2)$
and which satisfies $\mu(0,g)=\mu(g,0)=g$, where the origin $0\in\K^c$
corresponds to the identity element in the group structure.
\item[{\bf (2)}]
A $\K$-algebraic {\it inversion mapping}\, $\iota: \Delta_c(\delta_2)\to
\Delta_c(\delta_1)$ satisfying
$\mu(g,\iota(g))=
\mu(\iota(g),g)=0$ and $\iota(0)=0$ whenever $\iota(g)
\in\Delta_c(\delta_2)$.
\end{itemize} 
Here, the integer $c\in\N_*$ is the {\it dimension}\, of
$G$. We shall say that composition and
inversion are defined locally in a neighborhood of the identity
element. In the $\K$-analytic category, the corresponding definition
is similar.

Now, we can define the notion of local $\K$-algebraic Lie group
action.  Let $n\in\N_*$, let $x=(x_1,\dots,x_n)\in\K^n$ and let two
positive numbers satisfy $0<\rho_2<\rho_1$. Let $G$ be a local
$\K$-algebraic Lie group as defined just above.  We shall
formulate the desired definition by means of the two precise polydiscs
$\Delta_n(\rho_2)\subset\Delta_n(\rho_1)$.  This pair of polydiscs
represents a {\it local $\K$-algebraic manifold} up to changes of
$\K$-algebraic coordinates.  A {\it local $\K$-algebraic Lie group
action on a local $\K$-algebraic manifold}\, consists of a {\it
$\K$-algebraic}\, action mapping $x'=\Phi(x;g)$ defined over
$\Delta_n(\rho_2)\times\Delta_c(\delta_2)$ with values in
$\Delta_n(\rho_1)$ which satisfies:
\begin{itemize}
\item[{\bf (1)}] 
$\Phi(\Phi(x;g);g'))= \Phi(x;\mu(g, g'))$ whenever
$\Phi(x;g)\in\Delta_n(\rho_2)$ and $\mu(g,g')\in\Delta_c(\delta_2)$,
where the local group multiplication $\mu(g, g')$ is $\K$-algebraic as
above;
\item[{\bf (2)}]
$\Phi(x;e)=x$ and $\Phi(\Phi(x;g);\iota(g))=x$ whenever
$\Phi(x;g)\in\Delta_n(\rho_2)$ and $\iota(g)\in\Delta_c(\rho_2)$, where
the inverse group mapping $g\mapsto \iota(g)$ is $\K$-algebraic as above.
\end{itemize}
In this definition, it is allowed to suppose that $x\in \C^n$ and
$g\in \R^c$, which is the case to be considered in the sequel.  By
differentiation, every local $\K$-algebraic action gives rise to
vector fields defined over $\Delta_n(\rho_2)$ which are infinitesimal
generators of the action.  Indeed, let us consider the algebraically
dependent one-parameter families of complex algebraic biholomorphic
transformations $\Phi(x;0,\dots,0,g_i,0,\dots,0)=: \Phi_i(x;g_i)\equiv
(\Phi_{i,1}(x;g_i),\dots,\Phi_{i,n}(x;g_i))\in\K^n$, which we shall also
denote by $\Phi_{i,g_i}(x)$. In general, such a family does
not make a one-parameter group of transformations, but we
can nevertheless introduce the vector fields
$X_i(\Phi_{i,g_i}(x);g_i):= \partial_{g_i}\Phi_i(x;g_i)=
\sum_{l=1}^n\, \partial_{g_i}\Phi_{i,l}(x;g_i)\,\partial/\partial
x_l$.  We notice that the coefficients of these vector fields do in
general depend on the group parameter $g_i\in G$.

In fact, in the algebraic category, there is no hope to modify the
coordinates on the group in order that the infinitesimal generators of
the action are independent of the parameter coordinates $g_j$.
 For instance,
the trivial one-dimensional action (complex dilatation) defined by
$(z,w)\mapsto ((1+g)z,w)=:\Phi(z,w;g)$, where $(z,w)\in\C^2$ and
$g\in\C$ is clearly an algebraic action. Here, the infinitesimal
generator $X(x;g)=(1+g)^{-1} z\partial_z$ depends on the
parameter $g$. The only way to avoid the dependence upon $g$ of the
coefficient of $X$ is to change coordinates on the group by setting
$1+g:=e^\sigma$, $\sigma\in\C$, whence the action is represented by
$(z,w)\mapsto (e^{\sigma}z,w)=:\Phi(z,w;\sigma)$. Indeed, from the
group property $\Phi(\Phi(z,w;\sigma);\sigma')\equiv
\Phi(z,w;\sigma+\sigma')$, it is classical and immediate to deduce
that if we define the parameter independent vector field
$X^0(z,w):=\partial_\sigma\Phi(z,w;\sigma)\vert_{\sigma=0}=
z\,\partial_z$, then it holds that
$\partial_\sigma\Phi(z,w;\sigma)=e^\sigma
z\, \partial_z=X^0(\Phi(z,w;\sigma))$. So the infinitesimal generator of the
action is independent of the parameter $g$. However, the main
trouble here is that the algebraicity of the action is necessarily
lost since the flow of $X^0$ is
not algebraic (the reader may check that each right (or left) invariant
vector field on an algebraic local Lie group defines in general a
nonalgebraic one-parameter subgroup, {\it e.g.}
for ${\rm SO}(2,\R)$, ${\rm SL}(2,\C)$).

Consequently we may allow the infinitesimal
generators of an algebraic local Lie group action $x'=\Phi(x;g)$,
defined by $X_i(x;g_i):= [\partial_{g_i}\Phi_i]
(\Phi_{i,g_i}^{-1}(x);g_i)$ to depend on the group parameter $g_i$,
even if the families $(\Phi_{i,g_i}(x))_{g_i\in\K}$ do not
constitute one-dimensional subgroups of transformations.

\subsection*{2.4.~Algebraicity of complex flow foliations} 
Suppose now that $M$ is a real algebraic generic submanifold in
$\C^n$, for instance a hypersurface which is Levi nondegenerate at a
``center'' point $p\in M$ corresponding to the origin in the
coordinates $t=(t_1,\dots,t_n)$.  Let $X\in\mathfrak{Hol}(M)$ be an
infinitesimal CR automorphism. Even if, for fixed real $s$, the
biholomorphic mapping $t\mapsto \exp(sX)(t)$ is complex algebraic,
{\it i.e.}  the $n$ components of this biholomorphism are complex
algebraic functions by Webster's theorem [We1977], we know by
considering the infinitesimal CR automorphism $X_1:=i(z+1)\partial_z$
of the strong tube ${\rm Im}\, w=\vert z+1\vert^2+\vert z+1\vert^6-2$
in $\C^2$ passing through the origin, that the flow of $X$ is not
necessarily algebraic with respect to all variables $(s,t)$.

Nevertheless, we shall show that the local CR automorphism group of $M$ is a
local algebraic Lie group whose general transformations are of
the form $t'=H(t;e_1,\dots,e_c)$, where $t\in\C^n$ and
$(e_1,\dots,e_c)\in\R^c$ and where $H$ is algebraic with respect to
all its variables. Thus the ``time'' dependent vector fields defined by
$X_i(t;e_i):=[\partial_{e_i}H_i]( H_{i,e_i}^{-1}(t);e_i)$, where
$H_{i,e_i}(t):=H_i(t;e_i):= 
H(t;0,\dots,0,e_i,0,\dots,0)$, have an algebraic flow,
simply given by $(t,e_i)\mapsto H_i(t;e_i)$.  It follows
that each foliation defined by the complex integral curves of the time
dependent complex vector fields $X_i$, $i=1,\dots,c$, is a complex algebraic
foliation, {\it see}~\S3 below.  Now, we can state the main technical
theorem of this paper, whose proof is postponed to \S4, \S5 and \S6.

\def\thetheorem{2.1}\begin{theorem}
Let $M\subset \C^n$ be a {\rm real algebraic}\, connected
geometrically smooth generic submanifold of codimension $d\geq 1$ and
CR dimension $m=n-d\geq 1$.  Let $p\in M$ and assume that $M$ is
finitely nondegenerate and minimal at $p$. Then for every sufficiently
small nonempty open polydisc $\Delta_1$ centered at $p$, the following
three properties hold{\rm :}
\begin{itemize}
\item[{\bf (1)}]
The complex Lie algebra $\mathfrak{Hol}(M,\Delta_1)$ is of finite
dimension $c\in\N$ which depends only on the local geometry of $M$ in
a neighborhood of $p$.
\item[{\bf (2)}]
There exists a nonempty open polydisc $\Delta_2\subset \Delta_1$ also
centered at $p$ and a $\C^n$-valued mapping
$H(t;e)=H(t;e_1,\dots,e_c)$ with $H(t;0)\equiv t$ which is defined in
a neighborhood of the origin in $\C^n\times \R^c$ and which is
algebraic with respect to both its variables $t\in \C^n$ and
$e\in\R^c$ such that for every holomorphic map $h: \Delta_2\to
\Delta_1$ with $h(\Delta_2\cap M)\subset \Delta_1\cap M$ which is
sufficiently close to the identity map, there exists a unique $e\in
\R^c$ such that $h(t)=H(t;e)$.
\item[{\bf (3)}]
The mapping $(t,e)\mapsto H(t;e)$ constitutes
a $\K$-algebraic local Lie transformation group action. More precisely, there
exist a local multiplication mapping $(e,e')\mapsto \mu(e,e')$ and a
local inversion mapping $e\mapsto \iota(e)$ such that $H$, $\mu$ and
$\iota$ satisfy the axioms of local algebraic Lie group action as
defined in \S2.3.
\item[{\bf (4)}]
The $c$ ``time dependent'' holomorphic vector fields
\def\theequation{2.1}\begin{equation}
X_i(t;e_i):=[\partial_{e_i}H_i]( H_{i,e_i}^{-1}(t);e_i), 
\end{equation}
where $H_{i,e_i}(t):=H_i(t;e_i):= H(t;0,\dots,0,e_i,0,
\dots,0)$, have algebraic coefficients and have an
algebraic flow, given by $(t,e_i)\mapsto H_i(t;e_i)$.
\end{itemize}
In the case where $M$ is real analytic, the same theorem holds true
with the word ``algebraic'' everywhere replaced by the word
``analytic''.
\end{theorem}

 A special case of Theorem~2.1 was proved in [BER1999b] where,
apparently, the authors do not deal with the notion of local Lie
groups and consider the isotropy group of the point $p$, namely the
group of holomorphic self-maps of $M$ fixing $p$.  The consideration
of the complete local Lie group of biholomorphic self-maps of a piece
of $M$ in a neighborhood of $p$ (not only the isotropy group of $p$)
is crucial for our purpose, since we shall have to deal with strong
tubes $M\in\mathcal{T}_n^d$ for which the isotropy group of $p\in M$
is trivial. Sections~\S4,~\S5 and \S6 are devoted to the proof of Theorem
4.1, a precise statement of Theorem~2.1. We mention that our method of
proof of Theorem~2.1 gives a non optimal bound for the dimension of
$\mathfrak{Hol}(M,\Delta_1)$.  To our knowledge,
the upper bound $c\leq (n+1)^2-1$ is optimal only in codimension $d=1$
and in the Levi nondegenerate case.

\section*{\S3.~Proof of Theorem~1.1}

We take in this section Theorem 2.1 for granted. As explained in \S1.3
above, we shall conduct the proof of Theorem 1.1 in two essential
steps (\S\S 3.1 and 3.2). The strategy for the proof of Theorems~1.4
and~1.5 is similar and we prove them in \S\S3.3 and 3.4. Let
$M\in\mathcal{T}_n^d$ be a strong tube of codimension $d$ passing
through the origin in $\C^n$ given by the equations
$v_j=\varphi_j(y)$, $j=1,\ldots,d$. Assume that $M$ is
biholomorphically equivalent to a real algebraic generic submanifold
$M'$.

\medskip
\noindent
{\it First step.}
We show that an arbitrary real algebraic element
$M'\in\mathcal{T}_n^d$ can be straightened in some local complex
algebraic coordinates $t'=(t'_1,...,t'_n)\in\C^n$ in order that its
infinitesimal CR automorphisms are the $n$ holomorphic vector fields
of the specific form $X_i'=c_i'(t_i')\, \partial_{t_i'}$,
$i=1,\dots,n$, where the functions $c_i'(t_i')$ are algebraic.

\medskip
\noindent
{\it Second step.}
Assuming that there exists a biholomorphic equivalence $\Phi: M\to M'$
satisfying $\Phi_*(\partial_{t_i})= c_i'(t')\, \partial_{t_i'}$, we
prove by direct computation that all the first order derivatives of the
mapping $\psi'(y')$ must be algebraic.

\subsection*{3.1.~Proof of the first step}
Let $t'=\Phi(t)$ be such an equivalence, with $\Phi(0)=0$ and
$M':=\Phi(M)$ real algebraic.  Let $X_i:=\partial_{t_i}$,
$i=1,\dots,n$, be the $n$ infinitesimal CR automorphisms of $M$ and
set $X_i':=\Phi_*(X_i)$. Of course, we have
$[X_{i_1}',X_{i_2}']=\Phi_*([X_{i_1},X_{i_2}])=0$, so the CR
automorphism group of $M'$ is also $n$-dimensional and commutative.
Let us choose complex algebraic coordinates $t'$ in a neighborhood of
$0\in M'$ such that $X_i'\vert_0=\partial_{t_i'}\vert_0$.  Let us
apply Theorem~2.1 to the real algebraic submanifold $M'$, noting all
the datas with dashes. There exists an algebraic mapping $H'(t';
e)=H'(t';e_1,\dots,e_n)$ such that every local biholomorphic self-map
of $M'$ writes uniquely $t'\mapsto H'(t';e)$, for some $e\in\R^n$.  In
particular, for every $i=1,\dots, n$ and every small $s\in\R$, there
exists $e_s\in\R^n$ depending on $s$ such that $\exp(s X_i')(t')\equiv
H'(t';e_s)$.  From the commutativity of the flows of the $X_i'$,
i.e. from $\exp(s_1X_{i_1}'(\exp(s_2 X_{i_2}'(t'))))\equiv
\exp(s_2X_{i_2}'(\exp(s_1X_{i_1}'(t'))))$, we get
\def\theequation{3.1}\begin{equation}
H'(H'(t';e_2);e_1)\equiv
H'(H'(t';e_1);e_2).
\end{equation} 
This shows that the biholomorphisms $t'\mapsto
H_e'(t'):=H'(t';e)$ commute pairwise. In particular, if we define 
\def\theequation{3.2}\begin{equation}
G_i'(t';e_i):=H'(t';0,\dots,0,e_i,0,\dots,0),
\end{equation}
we have $G_{i_1}'(G_{i_2}'(t';e_2);e_1)\equiv
G_{i_2}'(G_{i_1}'(t';e_1);e_2)$.

Next, after making a linear change of coordinates in the $e$-space, we
can insure that $\partial_{e_i}G_i'(0;e_i)\vert_{e_i=0}=
\partial_{t_i'}\vert_0=X_i'\vert_0$ for $i=1,\dots,n$. Finally,
complexifying the real variable $e_i$ in a complex variable
$\epsilon_i$, we get mappings $G_i'(t_i';\epsilon_i)$ which are
complex algebraic with respect to both variables $t'\in\C^n$ and
$\epsilon_i\in\C$ and which commute pairwise. We can now state and
prove the following crucial proposition (where we have dropped the
dashes) according to which we can straighten commonly the $n$
one-parameter families of biholomorphisms $t'\mapsto
G_i'(t';\epsilon_i)$.

\def\theproposition{3.1}\begin{proposition} Let $t\mapsto
G_i(t;\epsilon_i)$, $i=1,\dots,n$, be $n$ one complex parameter
families of complex algebraic biholomorphic maps from a neighborhood
of $0$ in $\C^{n}$ onto a neighborhood of $0$ in $\C^n$ satisfying
$G_i(t;0)\equiv t$,
$\partial_{\epsilon_i}G_i(0;\epsilon_i)\vert_{\epsilon_i=0}=
\partial_{t_i}\vert_0$ and pairwise commuting: 
$G_{i_1}(G_{i_2}(t;\epsilon_2);\epsilon_1) \equiv
G_{i_2}(G_{i_1}(t;\epsilon_1);\epsilon_2)$. Then there exists a
complex algebraic biholomorphism of the form $t'\mapsto \Phi'(t')=:t$
of $\C^n$ fixing the origin with $d\Phi'(0)={\rm Id}$ such that if we
set $G_i'(t';\epsilon_i):={\Phi'}^{-1}(G_i(\Phi'(t');\epsilon_i))$,
where $t'=\Phi(t)$ denote the inverse of $t=\Phi'(t')$, then we have
\def\theequation{3.3}\begin{equation}
G_i'(t';\epsilon_i)\equiv
(t_1',\dots,t_{i-1}',G_{i,i}'(t_i';\epsilon_i),t_{i+1}',\dots,t_n'),
\end{equation}
where the functions $G_{i,i}'$ are complex algebraic, {\rm depend only on}
$t_i'$ $($and on $\epsilon_i)$ and satisfy $G_{i,i}'(t_i';0)\equiv
t_i'$ and $\partial_{\epsilon_i} G_{i,i}'
(0;\epsilon_i)\vert_{\epsilon_i=0}=1$.
\end{proposition}

\proof
First of all, we define the complex algebraic
biholomorphism
\def\theequation{3.4}\begin{equation}
\Phi_1': \ \ 
(t_1',t_2',\dots,t_n')\longmapsto
G_1(0,t_2',\dots,t_n';t_1')=:t.
\end{equation}
We have $d\Phi_1'(0)={\rm Id}$, because $G_1(t;0)\equiv t$ and
$\partial_{\epsilon_1} G_1(0;\epsilon_1)
\vert_{\epsilon_1=0}=\partial_{t_1}\vert_0$. Furthermore, since
$\partial_{\epsilon_1} G_1(0;\epsilon_1) \vert_{\epsilon_1=0}$ is
transversal to $\{(0,t_2,\dots,t_n)\}$, it also follows that a small
neighborhood of the origin in $\C_t^n$ is algebraically foliated by the
$(n-1)$-parameter family of complex curves
$\mathcal{C}_{t_2',\dots,t_n'}':= \{G_1(0,t_2',\dots,t_n';t_1'): \,
\vert t_1'\vert < \delta\}$ where $\delta>0$ is small and $t_2',\dots,
t_n'$ are fixed. The existence of this foliation
shows that the relation
\def\theequation{3.5}\begin{equation}
t^ *\sim t 
\ \ \ {\rm iff} \ \ \ \text{\rm there
exists} \  \epsilon_1 \ \text{\rm such that} \
t^*=G_1(t; \epsilon_1)
\end{equation}
is a local equivalence relation, whose equivalence classes are the
leaves $\mathcal{C}_{t_2',\dots,t_n'}'$ ({\it see} {\sc Figure~1}).

\bigskip
\begin{center}
\input straightening.pstex_t
\end{center}
\bigskip

Consequently, 
as we clearly have
\def\theequation{3.6}\begin{equation}
(0,t_2',\dots,t_n')\sim G_1(0,t_2',\dots,t_n';t_1')\sim
G_1(G_1(0,t_2',\dots,t_n';t_1');\epsilon_1),
\end{equation}
using the transitivity of the relation $\sim$,
it follows that there exists a complex number
$\varepsilon_{t',\epsilon_1}$ depending on $t'$ and 
on $\epsilon_1$ such that
\def\theequation{3.7}\begin{equation}
G_1(G_1(0,t_2',\dots,t_n';t_1');\epsilon_1)=
G_1(0,t_2',\dots,t_n';\varepsilon_{t',\epsilon_1}).
\end{equation}
By the very definition~\thetag{3.4} of $\Phi_1'$, this is equivalent to 
\def\theequation{3.8}\begin{equation}
\Phi_1(G_1(\Phi_1'(t');\epsilon_1))=
(\varepsilon_{t',\epsilon_1},t_2',\dots,t_n'),
\end{equation}
where $t'=\Phi_1(t)$ denotes the inverse of $t=\Phi_1'(t')$.
Finally, since the left hand side of~\thetag{3.8} is 
clearly a complex algebraic mapping of
$(t';\epsilon_1)$, it follows that there exists 
a complex algebraic function $G_{1,1}'(t';\epsilon_1)$ such that
we can write
\def\theequation{3.9}\begin{equation}
\Phi_1(G_1(\Phi_1'(t');\epsilon_1))\equiv
(G_{1,1}'(t';
\epsilon_1),t_2',\dots,t_n').
\end{equation}
So we have straightened the first family by means of
$\Phi_1'$. 

Next, we drop the dashes and 
we restart with $G_1(t;\epsilon_1)=
(G_{1,1}(t;\epsilon_1),t_2,\dots,t_n)$.
Then, similarly as above, 
by introducing the complex algebraic biholomorphism
\def\theequation{3.10}\begin{equation}
\Phi_2': \ \ 
(t_1',t_2',t_3',\dots,t_n')\longmapsto 
G_2(t_1',0,t_3',\dots,t_n';t_2'),
\end{equation}
which satisfies $d\Phi_2'(0)={\rm Id}$, and by denoting
by $t'=\Phi_2(t)$ the inverse of $t=\Phi_2'(t')$, we get again that
if we set $G_2'(t';\epsilon_2):=\Phi_2(G_2(\Phi_2'(t');
\epsilon_2))$, then 
\def\theequation{3.11}\begin{equation}
G_2'(t';\epsilon_2)\equiv
(t_1',G_{2,2}'(t';\epsilon_2),t_3',\dots,t_n'),
\end{equation}
where the complex algebraic function $G_{2,2}' (t';\epsilon_2)$
satisfies $\partial_{\epsilon_2} G_{2,2}'(0;\epsilon_2)
\vert_{\epsilon_2=0}= 1$ and $G_{2,2}'(t';0) \equiv t_2'$.

We also have to consider the modification of the first family of
biholomorphisms $G_1'(t';\epsilon_1):=\Phi_2(G_1(\Phi_2'(t');
\epsilon_1))$. Using in an essential way the commutativity, we may
compute
\def\theequation{3.12}\begin{equation}
\left\{
\aligned
\Phi_2'(G_1'(t';\epsilon_1))& \ =
G_1(\Phi_2'(t');\epsilon_1)\\
& \ 
=G_1(G_2(t_1',0,t_3',\dots,t_n';t_2');\epsilon_1)\\
& \
=G_2(G_1(t_1',0,t_3',\dots,t_n';\epsilon_1);t_2')\\
& \
=G_2(G_{1,1}(t_1',0,t_3',\dots,t_n';\epsilon_1),0,t_3',\dots,t_n';t_2')\\
& \ 
=\Phi_2'(G_{1,1}(t_1',0,t_3',\dots,t_n';\epsilon_1),t_2',t_3',\dots,t_n').
\endaligned\right.
\end{equation}
It follows that 
\def\theequation{3.13}\begin{equation}
G_1'(t';\epsilon_1)\equiv
(G_{1,1}(t_1',0,t_3',\dots,t_n';\epsilon_1),
t_2',t_3',\dots,t_n')
\end{equation}
whence 
\def\theequation{3.14}\begin{equation}
G_{1,1}'(t';\epsilon_1):=G_{1,1}(t_1',0,t_3',\dots,t_n';\epsilon_1)
\end{equation}
does not depend on $t_2'$. Finally, inserting~\thetag{3.11}
and~\thetag{3.13} in the commutativity relation
$G_1'(G_2'(t';\epsilon_2);\epsilon_1)\equiv
G_2'(G_1'(t';\epsilon_1);\epsilon_2)$, we find
\def\theequation{3.15}\begin{equation}
\left\{
\aligned
G_{1,1}'(t_1',G_{2,2}'(t';\epsilon_2),t_3',\dots,t_n';\epsilon_1)\equiv & \ 
G_{1,1}'(t';\epsilon_1),\\
G_{2,2}'(G_{1,1}'(t';\epsilon_1),t_2',t_3',\dots,t_n';\epsilon_2)\equiv & \
G_{2,2}'(t';\epsilon_2).
\endaligned\right.
\end{equation}
The first relation gives nothing, since we already know that
$G_{1,1}'$ is independent of $t_2'$. By differentiating the second
relation with respect to $\epsilon_1$ at $\epsilon_1=0$, we find that
$G_{2,2}'$ is independent of $t_1'$.

In summary, after the change of coordinates
$\Phi_2'\circ\Phi_1'(t')=t$ which is tangent to the identity map at
$t'=0$, we obtained that
\def\theequation{3.16}\begin{equation}
\left\{
\aligned
G_1'(t';\epsilon_1) & \ 
=(G_{1,1}'(t_1',t_3',\dots,t_n'),t_2',t_3',\dots,t_n'),\\
G_2'(t';\epsilon_2) & \
=(t_1',G_{2,2}'(t_2',t_3',\dots,t_n'),t_3',\dots,t_n').
\endaligned\right.
\end{equation}
Using these arguments, the proof of Proposition~3.1 clearly
follows by induction.
\endproof

Now, we come back to our CR manifold $M'$ having the one-parameter
families of algebraic biholomorphisms $G_i'(t';e_i)$ given
by~\thetag{3.2} and pairwise commuting. Applying Proposition~3.1,
after a change of complex algebraic coordinates of the form
$t'=\Psi''(t'')$, we may assume that the $G_i''(t'';\epsilon_i)$ are
algebraic and can be written in the specific form
\def\theequation{3.17}\begin{equation}
G_i''(t'';\epsilon_i)\equiv (t_1'',\dots,t_{i-1}'',
G_{i,i}''(t_i'';\epsilon_i),t_{i+1}'',\dots,t_n''),
\end{equation} 
with
$\partial_{\epsilon_i}G_{i,i}''(0;\epsilon_i)\vert_{\epsilon_i=0}=1$.
Let $t''=\Psi'(t')$ denote the inverse of $t'=\Psi''(t'')$. We thus have
$t''=\Psi'(t')=\Psi'(\Phi(t))$, where we remind that $t'=\Phi(t)$ provides
the equivalence between the strong tube $M$ and the algebraic CR
generic $M'$.

Since $\Psi'$ is algebraic, the image $M'':=\Psi'(M')$ is also algebraic. Let
$r_j'(t',\bar t')=0$, $j=1,\dots,d$, be defining equations for
$M'$. Then $r_j''(t'',\bar t''):=
r_j'(\Psi''(t''),\overline{\Psi''(t'')})=0$ are defining equations for
$M''$. By assumption, for $\epsilon_i:=e_i\in\R$ real, the family of
algebraic biholomorphisms $G_i'(t';\epsilon_i)$ maps a small piece of
$M'$ through the origin into $M'$. It follows trivially that
$G_i''(t'';\epsilon_i)\equiv\Psi' (G_i'(\Psi''(t'');\epsilon_i))$ maps a
small piece of $M''$ through the origin into $M''$. Furthermore, since
$d\Psi''(0)={\rm Id}$, it follows that if we denote
$X_i'':=\Psi_*'(X_i')$, then $X_i''\vert_0= \partial_{t_i''}\vert_0$.

Next, thanks to the
specific form~\thetag{3.17}, 
by differentiating $\partial_{\epsilon_i}G_i''(t'';
\epsilon_i)\vert_{\epsilon_i=0}$, we get $n$ vector fields of the form
$Z_i''=c_i''(t_i'')\, \partial_{t_i''}$.  By construction, the
functions $c_i''(t_i'')$ are algebraic and satisfy $c_i''(0)=1$.
Differentiating with respect to $e_i$ the identity
$r_j''(G_i''(t'';e_i), \overline{G_i''(t'';e_i)})=0$ for
$r_j''(t'',\bar t'')=0$, {\it i.e.} for $t''\in M''$, we see that
$Z_i''$ is tangent to $M''$, {\it i.e.} we see that $Z_i''$ is an
infinitesimal CR automorphism of $M''$. Consequently, there exist real
constants $\lambda_{i,l}$ such that $Z_i''=\sum_{l=1}^n\,\lambda_{i,l}\,
X_l''$. Since $Z_i''\vert_0=X_i''\vert_0=\partial_{t_i''}\vert_0$, we
have in fact $\lambda_{i,l}=1$ for $i=l$ and $\lambda_{i,l}=0$ for
$i\neq l$. So $Z_i''=X_i''$ and we have shown that
\def\theequation{3.18}\begin{equation}
(\Psi'\circ\Phi)_*(X_i)=X_i''=Z_i''=c_i''(t_i'')\,\partial_{t_i''}.
\end{equation}

We shall call a CR generic manifold $M''$ having infinitesimal CR
automorphisms of the form $X_i''=c_i''(t_i'')\,\partial_{t_i''}$ with
$c_i''(0)\neq 0$ a {\it pseudotube}. Such a pseudotube is not in
general a product by $\R^n$. In fact, there is no hope to tubify all
algebraic peudotubes in algebraic coordinates, as shows the elementary
example ${\rm Im}\, w=\vert z+1\vert^2+\vert z+1\vert^6-2$ having
infinitesimal CR automorphisms $\partial_w$ and $i(z+1)\partial_z$,
since the only change of coordinates for which
$\Phi_*(\partial_w)=\partial_{w'}$ and
$\Phi_*(i(z+1)\partial_z)=\partial_{z'}$ is $z+1=e^{iz'}$, $w=w'$,
which transforms $M$ into $M'$ of nonalgebraic defining equation ${\rm
Im}\, w'=e^{-2y'}+e^{-6y'}$.

The constructions of this paragraph
may be represented by the following symbolic picture.

\bigskip
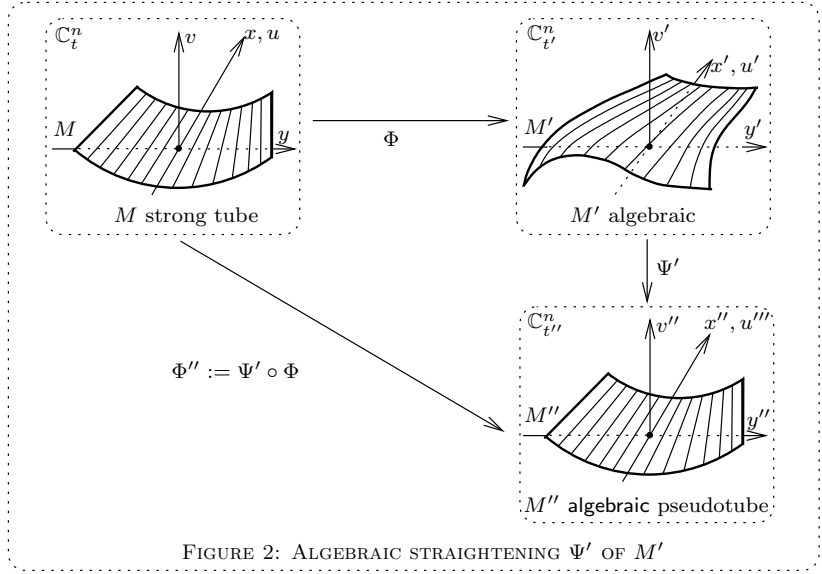
\begin{center}
\input redressement.pstex_t
\end{center}
\bigskip

\smallskip
\noindent
{\it Summary and conclusion of the first step.}  To conclude, let us
denote for simplicity $M''$ again by $M'$, the coordinates $t''$ again
by $t'$ and $t''=\Psi'\circ\Phi(t)$ by $t'=\Phi(t)$. After the above
straightenings, we have shown that the infinitesimal CR automorphisms
$X_i':=\Phi_*(X_i)$ of the algebraic generic manifold $M'$ are of the
sympathetic form $X_i'=c_i'(t_i')\,\partial_{t_i'}$, $i=1,\dots,n$,
with algebraic coefficients $c_i'(t_i')$ satisfying $c_i'(0)=1$.

\subsection*{3.2.~Proof of the second step}
We characterize first finite nondegeneracy for 
tubes of codimension $d$ in $\C^n$.

\def\thelemma{3.2}\begin{lemma}
Let $M$ be a tube of codimension $d$ in $\C^n$ equipped with
coordinates $(z,w)=(x+iy,u+iv)\in\C^m\times \C^d$ given by the
equations $v_j=\varphi_j(y)$, $j=1,\dots,d$, where $\varphi_j(0)=0$.
Then $M$ is finitely nondegenerate at the origin if and only if there
exist $m$ multi-indices $\beta_*^1,\dots,\beta_*^m\in\N^m$ with
$\vert\beta_*^k\vert\geq 1$ and integers $j_*^1,\dots,j_*^m$ with
$1\leq j_*^k\leq d$ such that the real mapping
\def\theequation{3.19}\begin{equation}
\psi(y):=\left(
{\partial^{\vert \beta_*^1\vert}\varphi_{j_*^1}(y)
\over \partial y^{\beta_*^1}},\ldots,
{\partial^{\vert \beta_*^m\vert}\varphi_{j_*^m}(y)
\over \partial y^{\beta_*^m}}
\right)=:y'\in\R^m
\end{equation}
is of rank $m$ at the origin in $\R^m$.
\end{lemma}

\proof
We follow the definition of finite nondegeneracy given in \S1.2. Let
$r_j(t,\bar t):=v_j-\varphi_j(y)=0$ be the defining equations of
$M$. Let $\overline{L}_k:= \partial_{\bar z_k}+\sum_{j=1}^d\,
\varphi_{j,\bar z_k}\,\partial_{\bar w_j}$, $k=1,\dots,m$, be a basis
of $(1,0)$-vector fields tangent to $M$.  We write the first order terms
in the Taylor series of $\varphi_j(y)$ as $\varphi_j(y)=
\sum_{l=1}^n\, \lambda_{j,l} \,y_l+{\rm O}(\vert y\vert^2)$.  Then the
holomorphic gradient of $r_j$ is given by
\def\theequation{3.20}\begin{equation}
\left\{
\aligned
\nabla_t(r_j) 
& \ 
=(\partial_{z_1}r_j,\dots,\partial_{z_m}r_j,
\partial_{w_1}r_j,\dots,\partial_{w_d}r_j)\\
& \
=i2^{-1}(\partial_{y_1}\varphi_{j},\dots,\partial_{y_m}\varphi_{j},0,
\dots,0,-1,0,\dots,0)\\
& \ 
=i2^{-1}(\lambda_{j,1},\cdots,\lambda_{j,m},0,\dots,0,-1,0,\dots,0),
\ \ \text{\rm at the origin}.
\endaligned\right.
\end{equation}
On the other hand, since for $\beta=(\beta_1,\dots,\beta_m)
\in\N^m$ with $\vert\beta\vert \geq 1$ the order $\vert\beta\vert$
derivation
$\overline{L}^\beta:=\overline{L}_1^{\beta_1}\cdots
\overline{L}_m^{\beta_m}$ acts on functions of $y$ as the
operator $(2i)^{-\vert\beta\vert}\, \partial_y^\beta$, 
we can compute
\def\theequation{3.21}\begin{equation}
\left\{
\aligned
\overline{L}^\beta(\nabla_t(r_j)) 
& \
=(\overline{L}^\beta\partial_{z_1}\varphi_{j},\dots,
\overline{L}^\beta\partial_{z_m}\varphi_{j},0,\dots,0,\dots,0)\\
& \
=i^{-\vert\beta\vert+1}2^{-\vert\beta\vert-1}\,
(\partial_y^\beta\partial_{y_1}\varphi_{j},\dots,
\partial_y^\beta\partial_{y_m}\varphi_{j},0,\dots,0,\dots,0).
\endaligned\right.
\end{equation}
By
inspecting the expressions~\thetag{3.20} and~\thetag{3.21},
we see that ${\rm Span}\{(\overline{L}^\beta(\nabla_t (r_j)))(0):
\beta\in\N^m, \, j=1,\dots,d\}=\C^n$ if and only if ${\rm Span}
\{(\partial_y^\beta\partial_{y_1}\varphi_{j}(0),\dots,
\partial_y^\beta\partial_{y_m}\varphi_{j}(0)): \beta\in\N^m, \,
\vert\beta\vert\geq 1, \, j=1,\dots,d\}=\R^m$. This last condition is
clearly equivalent to the one
stated in Lemma~3.2.
\endproof

We can prove now that
the inverse mapping $\psi'(y')$ of the mapping $\psi(y)$ defined
by~\thetag{1.1} (or~\thetag{3.19}) has algebraic first order
derivatives. By Step~1, there exists a biholomorphic transformation
$\Phi$ mapping the strong tube $M$ onto the algebraic pseudotube $M'$
with the property that
$\Phi_*(\partial_{t_i})=c_i'(t_i')\,\partial_{t_i'}$.  Writing
$\Phi(t)=(h_1(t),\dots,h_n(t))$, we have
$\Phi_*(\partial_{t_i})=\sum_{l=1}^n\,
h_{l,t_i}(t)\,\partial_{t_l'}=c_i'(t_i')\,\partial_{t_i'}$, so
$h_i(t)$ depends only on $t_i$ which yields
$\Phi(t)=(h_1(t_1),\dots,h_n(t_n))$. We shall use the convenient
notation $t_i'=h_i(t_i)$ and $t_i=h_i'(t_i')$ for the inverse
$h_i':=h_i^{-1}$, $i=1,\dots,n$.  If accordingly, $\Phi'(t')=t$
denotes the inverse of $\Phi(t)=t'$, we have
$\Phi_*'(c_i'(t_i')\,\partial_{t_i'})= c_i'(t_i')\, h_{i,t_i'}'(t_i')\,
\partial_{t_i} =\partial_{t_i}$, which shows that
$c_i'(t_i')\,h_{i,t_i'}'(t_i')\equiv 1$.  Since $c_i'(0)=1$, we see that
$h_i'(t_i')=\int_0^{t_i'} 1/[c_i'(\sigma)]\, d\sigma$ is the complex
primitive of an algebraic function. This observation will be
important.

After a permutation of the coordinates, we may assume that $M'$ is
given in the coordinates $t'=(z',w')\in\C^m\times\C^d$ by the real
defining equations ${\rm Im} \, w_j'=\varphi_j'(z',\bar z',{\rm Re}\, w')$,
$j=1,\dots,d$, where the functions $\varphi_j'$ are algebraic and
vanish at the origin. Solving in terms of $w'$ by means of the
algebraic implicit function theorem, we can represent $M'$ by the
algebraic complex defining equations
\def\theequation{3.22}\begin{equation}
w_j'=\overline{\Theta}_j'(z',\bar z',\bar w'),
\ \ \ \ \
j=1,\dots,d,
\end{equation}
where $\overline{\Theta}'$ satisfies the vectorial functional equation
$w'\equiv\overline{\Theta}'(z',\bar z',\Theta'(\bar z',
z',w'))$ (which we shall not use). 
According to the splitting $(z',w')$ of coordinates, it is
convenient to modify our previous notation by writing
$z_k=f_k'(z_k')$, $k=1,\dots,m$ and $w_j=g_j'(w_j')$, $j=1,\dots,d$
instead of $t_i=h_i'(t_i')$, $i=1,\dots,n$, and also
\def\theequation{3.23}\begin{equation}
\left\{
\aligned
X_k'= & \ a_k'(z_k')\,\partial_{z_k'}, \ \ \ \ \
k=1,\dots,m, \ \ \ \ \
a_k'(0) = 1,\\
Y_j'= & \ b_j'(w_j')\,\partial_{w_j'}, \ \ \ \ \
j=1,\dots,d, \ \ \ \ \ \
b_j'(0) = 1,
\endaligned\right.
\end{equation}
instead of $X_i'=c_i'(t_i')\,\partial_{t_i'}$. The relation
$c_i'(t_i')\, h_{i,t_i'}'(t_i')\equiv 1$ rewrites down in the form
\def\theequation{3.24}\begin{equation}
\left\{
\aligned
a_k'(z_k')\,f_{k,z_k'}'(z_k')\equiv 
& \ 1,\\
b_j'(w_j')\,g_{j,w_j'}'(w_j')\equiv
& \ 1.
\endaligned\right.
\end{equation}
We remind that the derivatives
of the $f_k'$ and of the $g_j'$ are algebraic.
Let now $t'=(z',w')\in M'$, thus satisfying~\thetag{3.22}. Then 
$h'(t')=(f'(z'),g'(w'))$ belongs to $M$, namely  we have for
$j=1,\dots,d$:
\def\theequation{3.25}\begin{equation}
{g_j'(w_j')-\bar g_j'(\bar w_j')\over
2i}=
\varphi_j\left(
{f_1'(z_1')-\bar f_1'(\bar z_1')\over 2i},\dots,
{f_m'(z_m')-\bar f_m'(\bar z_m')\over 2i}\right),
\end{equation}
where $i=\sqrt{-1}$ here.
Replacing $w_j'$ by $\overline{\Theta}_j'(z',\bar z',\bar w')$ in the
left hand side, we get the following identity between converging power series
of the $2m+d$ complex variables $(z',\bar z',\bar w')$:
\def\theequation{3.26}\begin{equation}
{g_j'(\overline{\Theta}_j'(z',\bar z',\bar w'))-\bar g_j'(\bar w_j')
\over
2i}\equiv
\varphi_j\left(
{f_1'(z_1')-\bar f_1'(\bar z_1')\over 2i},\dots,
{f_m'(z_m')-\bar f_m'(\bar z_m')\over 2i}\right).
\end{equation}
Let us differentiate this identity with respect to 
$z_k'$, for $k=1,\dots,m$. Taking into account the
relations~\thetag{3.24}, we obtain
\def\theequation{3.27}\begin{equation}
{a_k'(z_k')\,
\overline{\Theta}_{j,z_k'}'(z',\bar z',\bar w')\over
b_j'(\overline{\Theta}_j'(z',\bar z',\bar w'))}\equiv
{\partial\varphi_j\over\partial y_k}
\left(
{f_1'(z_1')-\bar f_1'(\bar z_1')\over 2i},\dots,
{f_m'(z_m')-\bar f_m'(\bar z_m')\over 2i}
\right).
\end{equation}
Clearly, the left hand side is an algebraic function 
$\mathcal{A}_{j,k}'(z',\bar z',\bar w')$. Then 
differentiating again with respect to the variables
$z_k'$ the relations~\thetag{3.27}, we see that
for every multi-index $\beta\in\N^m$ with 
$\vert\beta\vert\geq 1$, and every 
$j=1,\dots,d$, there exists an algebraic
function $\mathcal{A}_{j,\beta}'(z',\bar z',\bar w')$ such that
the following identity holds:
\def\theequation{3.28}\begin{equation}
\mathcal{A}_{j,\beta}'(z',\bar z',\bar w')\equiv
{\partial^{\beta_1+\cdots+\beta_m}\varphi_j\over
\partial y_1^{\beta_1}\cdots\partial_{y_m}^{\beta_m}}
\left(
{f_1'(z_1')-\bar f_1'(\bar z_1')\over 2i},\dots,
{f_m'(z_m')-\bar f_m'(\bar z_m')\over 2i}
\right).
\end{equation}
Differentiating~\thetag{3.28} with respect to 
$\bar w'$, we see immediately that $\mathcal{A}_{j,\beta}'$ is
in fact independent of $\bar w'$. Furthermore, we see that
$\mathcal{A}_{j,\beta}'$ is real, namely 
$\mathcal{A}_{j,\beta}'(z',\bar z')\equiv
\overline{\mathcal{A}}_{j,\beta}'(\bar z',z')$.
Now we extract from~\thetag{3.28} the $m$ identities written
for $\beta:=\beta_*^k$, $j:=j_*^k$, $k=1,\dots,m$ and we use
the invertibility of the mapping $\psi$ defined in~\thetag{3.19}
(recall that $\psi'(y')=y$ denotes the inverse of $y'=\psi(y)$), 
which yields
\def\theequation{3.29}\begin{equation}
{f_k'(z_k')-\bar f_k'(\bar z_k')\over
2i}\equiv
\psi_k'(\mathcal{A}_{j_*^1,\beta_*^1}'(z',\bar z'),\dots,
\mathcal{A}_{j_*^m,\beta_*^m}'(z',\bar z')),
\end{equation}
for $k=1,\dots,m$. For simplicity, 
we shall write $\mathcal{A}_k'(z',\bar z')$ instead
of $\mathcal{A}_{j_*^k,\beta_*^k}'(z',\bar z')$. 
Finally,
we differentiate~\thetag{3.29} with 
respect to $z_k'$, which yields, taking
into account~\thetag{3.24}:
\def\theequation{3.30}\begin{equation}
\left\{
\aligned
{1\over 2i\, a_k'(z_k')}\equiv 
& \ 
\sum_{l=1}^m\, 
{\partial \psi_k'\over\partial y_l'}
(\mathcal{A}_1'(z',\bar z'),\dots,
\mathcal{A}_m'(z',\bar z'))\,
{\partial \mathcal{A}_l'
\over \partial z_k'}
(z',\bar z'),\\
0 \equiv
& \ 
\sum_{l=1}^m\, {\partial\psi_k'\over\partial y_l'}
(\mathcal{A}_1'(z',\bar z'),\dots,\mathcal{A}_m'(z',\bar z'))\,
{\partial \mathcal{A}_l'\over\partial z_{\widetilde{k}}'}(z',\bar z'), \ \ 
\ \ \
\widetilde{k}\neq k.
\endaligned\right.
\end{equation}
It follows from these relations~\thetag{3.30} viewed in matrix form
that the constant matrix $({\partial \mathcal{A}_l'\over\partial
z_k'}(0,0))_{1\leq l,k\leq m}$ is invertible, because the diagonal
matrix $(\delta_k^{\widetilde{k}} \, [2ia_k'(z_k')]^{-1})_{1\leq
k,\widetilde{k}\leq m}$ is evidently invertible at $z_k'=0$ (recall
$a_k'(0)=1$). Consequently, there exist algebraic functions
$\mathcal{B}_{k,l}'(z',\bar z')$ so that
\def\theequation{3.31}\begin{equation}
{\partial \psi_k'\over\partial y_l'}(
\mathcal{A}_1'(z',\bar z'),\dots,\mathcal{A}_m'(z',\bar z'))\equiv
\mathcal{B}_{k,l}'(z',\bar z').
\end{equation}
Next, setting $\widetilde{y}_k'=\mathcal{A}_k'(iy',-iy')$,
$k=1,\dots,m$ we see, from the invertibility of the matrix $({\partial
\mathcal{A}_k'\over\partial z_l'}(0,0))_{1\leq k,l\leq m}$ and from the 
reality of $\mathcal{A}_k'(z',\bar z')$, that the
Jacobian determinant at the origin of the mapping $y'\mapsto
\mathcal{A}'(iy',-iy')=\widetilde{y}_k'$ 
is nonzero. Thus there are real algebraic functions $\mathcal{C}_k'$  so 
that we can express $y'$ in terms of
$\widetilde{y}'$ as $y_k'=\mathcal{C}_k'(\widetilde{y}')$. Finally, we get
\def\theequation{3.32}\begin{equation}
{\partial\psi_k'\over\partial y_l'}(\widetilde{y}_1',\dots,\widetilde{y}_m')=
\mathcal{B}_{k,l}'(i\mathcal{C}'(\widetilde{y}'),
-i\mathcal{C}'(\widetilde{y}')), 
\end{equation}
where the right hand sides are algebraic; this shows that the partial
derivatives ${\partial_{y'l} \psi_k'}$ are
algebraic functions of $\tilde{y}'$.  

To obtain the equivalent formulation of Theorem~1.1, we observe
the following.

\def\thelemma{3.3}\begin{lemma}
For every $k,l=1,\dots,m$,
the functions $\partial_{y_k'}\psi_l'(y')$ are algebraic functions
of $y'$ if and only if for every $k_1,k_2=1,\dots,m$, 
the second derivative $\partial^2_{y_{k_1}y_{k_2}}(y)$ is an algebraic
function of $\psi(y)=(\partial_{y_1}\varphi(y),\dots,\partial_{y_m}
\varphi(y))$. 
\end{lemma}

\proof
Differentiating the identities $y_k\equiv \psi_k'(\psi(y))$, 
$k=1,\dots,m$, with respect to $y_l$, we get
\def\theequation{3.33}\begin{equation}
\delta_k^l\equiv \sum_{j=1}^{m}\, 
\partial_{y'_j}\psi_{k}'(\psi(y))\, \partial_{y_j}\psi_{j}(y)\equiv
\sum_{j=1}^{m}\, \partial_{y'_j}\psi_{k}'(y')\, 
\partial^2_{y_my_l}\varphi(y).
\end{equation}
Applying Cramer's rule, we see that there exist universal rational 
functions $R_{k,l}$ and $S_{k,l}$ such that
\def\theequation{3.34}\begin{equation}
\left\{
\aligned
\partial^2_{y_ky_l}\varphi(y)\equiv R_{k,l}(\{
\partial_{y'_{k_2}}\psi_{k_1}'(y')\}_{1\leq k_1,k_2\leq m}\}),\\
\partial_{y'_l}\psi_{k}'(y')\equiv S_{k,l}(\{
\partial^2_{y_{k_1}y_{k_2}}\varphi(y)\}_{1\leq k_1,k_2\leq m}\}).
\endaligned\right.
\end{equation}
This implies the equivalence of Lemma~3.3.
\endproof
In conclusion, taking
Theorem~2.1 for granted, the proof of Theorem~1.1 is now complete.
\qed
\subsection*{3.3.~Proof of Theorem~1.5}
Let $M:v=\varphi(z,\bar z)$ be a rigid Levi nondegenerate hypersurface
in $\C^n$ passing through the origin. We may assume that
$v=\sum_{k=1}^{n-1}\, \varepsilon_k\, \vert
z_k\vert^2+\varphi^3(z,\bar z)$, where $\varepsilon_k =\pm 1$ and we
may write $\varphi^3(z,\bar z)= \sum_{k=1}^{n-1}\,[\bar z_k
\,\varphi_k^3(z)+z_k\bar \varphi_k^3(\bar z)] +\varphi^4(z,\bar z)$,
with $\varphi^4(0,\bar z)\equiv \varphi_{z_k}^4(0,\bar z)\equiv 0$ and
$\varphi_k^3={\rm O}(2)$. After making the change of coordinates
$z_k':=z_k+\varepsilon_k\, \varphi_k^3(z)$, $w':=w$, we come to the
simple equation $v'=\sum_{k=1}^{n-1}\, \varepsilon_k\, \vert
z_k'\vert^2+\chi'(z',\bar z')$, where $\chi'(0,\bar z')\equiv
\chi_{z_k'}(0,\bar z')\equiv 0$, considered in Theorem 1.5.

Assume that $M$ is strongly rigid, locally algebraizable and let $M'$
be an algebraic equivalent of $M$. Let $t'=h(t)$ be such an
equivalence, or in our previous notation $z'=f(z,w)$ and
$w'=g(z,w)$. We note $z=f'(z',w')$ and $w=g'(z',w')$ the inverse
equivalence. Since $M$ is strongly rigid, namely $\mathfrak{Hol}(M)$
is generated by the single vector field $X_1:=\partial_w$, it follows
that $\mathfrak{Hol}(M')$ is also one-dimensional, generated by the
single vector field $X_1':=h_*(X_1)$. Taking again Theorem~2.1 for
granted and proceeding as in the first step of the proof of
Proposition~3.1, we may algebraically straighten the complex foliation
induced by $X_1'$ to the ``vertical'' foliation by
$w'$-lines. Equivalently, we may assume that
$X_1'=b'(z',w')\,\partial_{w'}$ with $b'$ algebraic and $b'(0)=1$.
The assumption $h_*(\partial_w)=b'(z',w')\, \partial_{w'}$ yields that
$f'(z',w')$ is independant of $w'$ and that $b'(z',w')\,
g_{w'}'(z',w')\equiv 1$, so that as in~\thetag{3.24} above, the
derivative $g_{w'}'$ is algebraic.  Let $w'=\overline{\Theta}'
(z',\bar z',\bar w')$ be the complex defining equation of $M'$ in
these coordinates. The assumption $h'(M')=M$ yields the following
power series identity
\def\theequation{3.35}\begin{equation}
g'(z',\overline{\Theta}'(z',\bar z',\bar w'))-
\bar g'(\bar z',\bar w')\equiv 2i\,\varphi(f'(z'),\bar f'(\bar z')).
\end{equation}
By differentiating this identity with respect to $z_k'$, we get
\def\theequation{3.36}\begin{equation}
\partial_{z_k'}g'(z',\overline{\Theta}'(z',\bar z',\bar w'))+
{\partial_{z_k'}\overline{\Theta}'(z',\bar z',\bar w')\over
b'(\overline{\Theta}'(z',\bar z',\bar w'))}\equiv
2i\, \sum_{l=1}^{n-1}\, \partial_{z_l}\varphi(f'(z'),\bar f'(\bar z'))\, 
\partial_{z_k'}f_l'(z').
\end{equation}
We notice that the second term in the left hand side of~\thetag{3.36}
is algebraic. By differentiating in turn~\thetag{3.36} with respect to
$\bar z_k'$ and using the algebraicity of $\partial^2_{z_k'w'} 
g'(z',\bar\Theta'(z',\bar z',\bar
w'))$, we obtain that there exist algebraic functions
$\mathcal{A}_{k_1,k_2}'(z',\bar z')$ such that
\def\theequation{3.37}\begin{equation}
\mathcal{A}_{k_1,k_2}'(z',\bar z')\equiv 
\sum_{l_1,l_2=1}^{n-1}\, 
\partial^2_{z_{l_1}\bar z_{l_2}}\varphi(f'(z'),\bar f'(\bar z'))\, 
\partial_{z'_{k_1}}f_{l_1}'(z')\, 
\partial_{\bar z_{k_2}'}\bar f_{l_2}'(\bar z').
\end{equation}
Without loss of generality, we may assume that $h'$ is tangent to the
identity map at $t'=0$. Then setting $\bar z':=0$ in~\thetag{3.37} and
using the fact that $\partial^2_{z_{l_1}\bar z_{l_2}}\varphi(z,0)=
\delta_{l_1}^{l_2}\varepsilon_{l_1}+ \partial^2_{z_{l_1}\bar z_{l_2}}\chi(z,0)
\equiv \delta_{l_1}^{l_2}\varepsilon_{l_1}$ by the
properties of $\chi$ in Theorem~1.5 we get, since $\partial_{\bar z_{k_2}'}
\bar f_{l_2}'(0)=\delta_{l_2}^{k_2}$~:
\def\theequation{3.38}\begin{equation}
\mathcal{A}_{k_1,k_2}'(z',0)\equiv
\varepsilon_{k_2}\,\partial_{z_{k_1}'} f_{k_2}'(z'),
\end{equation}
which shows that all the first order derivatives $\partial_{z_k'}f_{l}'(z')$
are algebraic. 

Next, since the canonical transformation to normalizing coordinates is
algebraic and preserves the ``horizontal'' coordinates $z'$ ({\it
cf.}~[CM1974]), hence does not perturb the complex foliation induced
by $X_1'$, we may also assume that $M'$ is given in normal
coordinates, namely that the function $\Theta'$ satisfies
$\Theta'(0,\bar z',\bar w')\equiv \Theta'(z',0,\bar w')\equiv \bar
w'$. Since the coordinates are normal for both $M$ and $M'$, it
follows by setting $\bar z':=0$ and $\bar w':=0$ in~\thetag{3.35} that
$g'(z',0)\equiv 0$. Consequently, $\partial_{z_k'}g'(z',0)\equiv
0$. Finally, by setting $z':=0$ and $\bar w':=0$ in~\thetag{3.36}, we
see that the first term in the left hand side vanishes and that the
second term is algebraic with respect to $\bar z'$, so we obtain that
there exist algebraic functions $\overline{\mathcal{B}}'_k (\bar z')$
such that
\def\theequation{3.39}\begin{equation}
\overline{\mathcal{B}}_k'(\bar z')\equiv 
\sum_{l=1}^{n-1}\, \partial_{z_l}\varphi(0,\bar f'(\bar z'))\, 
\partial_{z_k'}f_{l}'(0)\equiv \varepsilon_k \, \bar f_k'(\bar z').
\end{equation}
We have proved that the components $f_k'(z')$ are 
all algebraic.

Finally, coming back to the relation~\thetag{3.36}, we want to
prove that the derivatives $\partial_{z_l}\varphi(f'(z'),\bar f'(\bar z'))$
are all algebraic. However, the first term of~\thetag{3.36} is not
algebraic in general. Fortunately, using the fact that
$\overline{\Theta}'=\bar w'+{\rm O}(2)$, we see that there exists a
unique algebraic solution $\bar w'=\overline{\Lambda}'(z',\bar z')$ of
the implicit equation $\overline{\Theta}'(z',\bar z',\bar w')=0$,
namely satisfying $\overline{\Theta}'(z',\bar z',\overline{\Lambda}'(z',\bar
z'))\equiv 0$. Then by replacing $\bar w'$ by $\overline{\Lambda}'$ 
in~\thetag{3.36}, we get that there exist algebraic functions 
$\mathcal{C}_k'(z',\bar z')$ such that
\def\theequation{3.40}\begin{equation}
\mathcal{C}_k'(z',\bar z')\equiv 
\sum_{l=1}^{n-1}\, \partial_{z_l}\varphi(f'(z'),\bar f'(\bar z'))\, 
\partial_{z_k'}f_{l}'(z').
\end{equation}
Since $f'$ is tangent to the identity map, we can solve by Cramer's
rule this linear system for the derivatives $\partial_{z_l}\varphi$, which
yields that the $\partial_{z_l}\varphi(f'(z'),\bar f'(\bar z'))$ are all
algebraic. Since $f'(z')$ is also algebraic, we obtain in sum that the
derivatives $\partial_{z_l}\varphi(z,\bar z)$ are all algebraic.  In
conclusion, taking Theorem~2.1 for granted, the proof of Theorem~1.5
is complete. \qed

\subsection*{3.4.~Proof of Theorem~1.4}
Let $M: v=\varphi(z\bar z)$ in $\C^2$ with $\mathfrak{Hol}(M)$ generated by
$\partial_w$ and $iz\partial_z$. Without loss of generality, we can
assume that $\varphi(r)=r+{\rm O}(r^2)$. Let $M'$ be an algebraic
equivalent of $M$. Let $t=h'(t')$, or $z=f'(z',w')$, $w=g'(z',w')$ be
a local holomorphic equivalence satisfying $h'(M')=M$. Let $t'=h(t)$ be its
inverse. Then $\mathfrak{Hol}(M')$ is two-dimensional and generated by
$h_*(\partial_w)$ and $h_*(iz\partial_z)$. First of all, using the
algebraicity of the CR automorphism group of $M'$ and proceeding as in
the proof of Proposition~3.1, we can prove that there exist two
generators of $\mathfrak{Hol}(M')$ of the form
$X_1'=b'(w')\,\partial_{w'}$ and $X_2'=a'(z')\,\partial_{z'}$ where
$b'$ and $a'$ are algebraic and satisfy $b'(0)=1$ and $a'(z')=iz'+{\rm
O}({z'}^2)$. Furthermore, we may assume that $h'$ is tangent to the
identity map and that $h_*'(b'(w')\,\partial_{w'})=\partial_w$ and
$h_*'(a'(z')\,\partial_{z'})=iz\partial_z$. As in~\thetag{3.24}, it
follows that $b'(w')\, g_{w'}'(w')\equiv 1$ and $a'(z')\,
f_{z'}'(z')\equiv i f'(z')$.  Let $w'=\overline{\Theta}' (z',\bar
z',\bar w')$ be the complex algebraic equation of $M'$.  Then we get
the following power series identity:
\def\theequation{3.41}\begin{equation}
g'(\overline{\Theta}'(z',\bar z',\bar w'))-
\bar g'(\bar w')\equiv 2i\, \varphi(f'(z')\bar f'(\bar z')),
\end{equation}
which yields after differentiating with respect to $z'$:
\def\theequation{3.42}\begin{equation}
\left\{
\aligned
\overline{\Theta}_{z'}'(z',\bar z',\bar w')/
[b'(\overline{\Theta}'(z',\bar z',\bar w'))]\equiv 
& \
2i\, \partial_{z'}f(z')\, \bar f'(\bar z')\, 
\partial_r\varphi(f'(z')\bar f'(\bar z'))\\
\equiv & \ 
-2 \, f'(z')\bar f'(\bar z') \, \partial_r\varphi(f'(z')\bar
f'(\bar z'))/[a'(z')].
\endaligned
\right.
\end{equation}
Here, we consider the function $\varphi$ as a function $\varphi(r)$ of
the real variable $r\in\R$. Since the left hand side is an algebraic 
function and $a'(z')$
is also algebraic, there exists an algebraic function
$\mathcal{A}'(z',\bar z')$ such that we can write
\def\theequation{3.43}\begin{equation}
\mathcal{A}'(z',\bar z')\equiv
f'(z')\bar f'(\bar z') \, \partial_r\varphi(f'(z')\bar
f'(\bar z')).
\end{equation}
Next, using the property $\varphi(r)=r+{\rm O}(r^2)$,
differentiating~\thetag{3.43} with respect to $\bar z'$ at $\bar
z'=0$, we obtain that $f'(z')$ is algebraic. Coming back
to~\thetag{3.43}, this yields that $\partial_r\varphi(f'(z')\bar f'(\bar z'))$
is algebraic. Since $f'(z')$ is also algebraic, we finally obtain that
$\partial_r\varphi(r)$ is algebraic. Excepting the examples which will be 
treated in \S7.5, the proof of Theorem~1.4 is complete. 
\qed

The next three sections are devoted to the statement of Theorem 4.1,
which implies directly Theorem 2.1 (\S4), and to its proof (\S\S5-6).

\section*{\S4.~Local Lie group structure for the CR automorphism group}

\subsection*{4.1.~Local representation of a real 
algebraic generic submanifold}  
We consider a connected real algebraic (or more generally, real analytic)
generic submanifold $M$ in $\C^n$ of codimension $d\geq 1$ and CR
dimension $m=n-d\geq 1$. Pick a point $p\in M$ and consider some
holomorphic coordinates $t=(t_1,\dots,t_n)=
(z_1,\dots,z_m,w_1,\dots,w_d)\in\C^m\times\C^d$ vanishing at $p$ in
which $T_0M=\{{\rm Im}\, w=0\}$. If we denote $w=u+iv$, it follows
that there exists (Nash) real algebraic power series $\varphi_j(z,\bar
z,u)$ with $\varphi_j(0)=0$ and $d\varphi_j(0)=0$ such that the
defining equations of $M$ are of the form $v_j=\varphi_j(z,\bar z,u)$,
$j=1,\dots,d$ in a neighborhood of the origin. By means of the
algebraic implicit function theorem, we can solve with respect to
$\bar w$ the equations $w_j-\bar w_j=2i\,\varphi_j(z,\bar z,(w+\bar
w)/2)$, $j=1,\dots, d$, which yields $\bar w_j=\Theta_j(\bar z,z,w)$
for some power series $\Theta_j$ which are complex algebraic with
respect to their $2m+d$ variables. Here, we have $\Theta_j=w_j+{\rm
O}(2)$, since $T_0M=\{{\rm Im}\, w=0\}$. Without
loss of generality, we shall assume that the
coordinates are normal, namely the functions $\Theta_j(\bar z,z,w)$
satisfy $\Theta_j(0,z,w)\equiv w_j$ and $\Theta_j(\bar z, 0, w)\equiv
w_j$. It may be shown that the
power series $\Theta_j=w_j+{\rm O}(2)$ satisfy the vectorial
functional equation $\Theta(\bar z,z,\overline{\Theta}(z,\bar z,\bar
w))\equiv \bar w$ in $\C\{z,\bar z,\bar w\}^d$ and conversely that to
every such power series mapping satisfying this vectorial functional
equation, there corresponds a unique real algebraic generic manifold
$M$ ({\it cf.}~for instance the manuscript
[GM2001c] for the details). So we can equivalently take $\bar
w_j=\Theta_j(\bar z,z,w)$ or $w_j=\overline{\Theta}_j(z,\bar z,\bar
w)$ as complex defining equations for $M$.

For arbitrary $\rho>0$, we shall
often consider the open polydisc $\Delta_n(\rho):=\{t\in\C^n:\, \vert
t\vert < \rho\}$ where we denote by $\vert t\vert:=\max_{1\leq i\leq
n}\, \vert t_i\vert$ the usual polydisc norm.  Without loss of
generality, we may assume that the power series $\Theta_j$ converge
normally in the polydic $\Delta_{2m+d}(2\rho_1)$, where $\rho_1>0$. In
fact, we shall successively introduce some other positive constants
(radii) $0<\rho_5<\rho_4<\rho_3<\rho_2<\rho_1$
afterwards. Finally, we define $M$ as:
\def\theequation{4.1}\begin{equation}
M=\{(z,w)\in \Delta_n(\rho_1):\,
\bar w_j=\Theta_j(\bar z,z,w), \ j=1,\dots, d\}.
\end{equation}
Next, let $\rho_2$ arbitrary with $0<\rho_2<\rho_1$. For 
$h',h\in\mathcal{O}(\Delta_n(\rho_1),\C^n)$, we define
\def\theequation{4.2}\begin{equation}
\vert\vert h'-h\vert\vert_{\rho_2}:=
\sup \, \{\vert h'(t)-h(t)\vert : \, t\in
\Delta_n(\rho_2)\}.
\end{equation}
For $k\in\N$, we shall also consider the $\mathcal{C}^k$ norms
\def\theequation{4.3}\begin{equation}
\vert\vert J^kh'-J^kh\vert\vert_{\rho_2}:=
\sup \, \{\vert\partial_t^\alpha h'(t)-\partial_t^\alpha h(t)\vert : \,
t\in\Delta_n(\rho_2), \, \alpha\in\N^n, \,
\vert\alpha\vert \leq k\}.
\end{equation}

For $k\in\N$ and $t\in\Delta_n(\rho_1)$, we denote by $J^kh(t)$ the
collection of partial derivatives $(\partial_t^\alpha h_i(t))_{1\leq
i\leq n, \, \vert \alpha\vert \leq k}$ of length $\leq k$ of the
components $h_1,\dots,h_n$, so $J^kh(t)\in\C^{N_{n,k}}$, where
$N_{n,k}:=n{(n+k)!\over n!\ k!}$.  In particular, the expression $J^k
h(0)=(\partial_t^\alpha h_i(0))_{1\leq i\leq n, \, \vert \alpha\vert
\leq k}$ denotes the $k$-jet of $h$ at $0$. So, the space of $k$-jets
at the origin of holomorphic mappings
$h\in\mathcal{O}(\Delta_n(\rho_1),\C^n)$ may be identified with the
complex linear space $\C^{N_{n,k}}$.  We denote the natural
coordinates on $\C^{N_{n,k}}$ by $(J_i^\alpha)_{ 1\leq i\leq n,\,
\vert\alpha\vert \leq k}$. Sometimes, we abbreviate this collection of
coordinates by $J^k\equiv(J_i^\alpha)_{ 1\leq i\leq n,\,
\vert\alpha\vert \leq k}$. Finally, we denote by $J_{\rm Id}^k$ the
$k$-jet at the origin of the identity mapping. We introduce the
important set of holomorphic self-mappings of $M$ defined by
\def\theequation{4.4}\begin{equation}
\left\{
\aligned
\mathcal{H}_{M,k,\varepsilon}^{\rho_2,\rho_1}:=
\{h\in\mathcal{O}(\Delta_n(\rho_1),\C^n): \ & \ 
\vert\vert 
J^kh -J_{\rm Id}^k \vert\vert_{\rho_2} < \varepsilon, \\
& \ \,
h(M\cap \Delta_n(\rho_2))\subset
M\cap \Delta_n(\rho_1)\}.
\endaligned\right.
\end{equation}
Here, $k\in\N$ and $\varepsilon>0$ is a small positive
number that we shall shrink many times in the sequel.

\bigskip
\begin{center}
\input polydiscs.pstex_t
\end{center}
\bigskip

We may now state the main theorem of \S4, \S5 and \S6, namely
Theorem~4.1, which provides a complete parametrized
description of the set $\mathcal{H}_{M,k,\varepsilon}^{\rho_2,\rho_1}$
of local biholomorphic self-mappings of $M$, with $k$ equal to an
integer $\kappa_0$ depending on $M$. During the course of the (rather
long) proof, for technical reasons, we shall have to introduce first a
third positive radius $\rho_3$ with $0<\rho_3<\rho_2<\rho_1$ which is
related to the finite nondegeneracy of $M$, and then afterwards a
fourth positive radius $\rho_4$ with $0<\rho_4<\rho_3<\rho_2<\rho_1$,
which is related to the minimality of $M$.  This is why the radius
notation ``$\rho_4$'' appears after ``$\rho_2$'' and ``$\rho_1$''
without mention of ``$\rho_3$'' ({\it cf.}~{\sc Figure~3}).

\def\thetheorem{4.1}\begin{theorem}
Assume that the real algebraic generic submanifold $M$ defined
by~\thetag{4.1} is minimal and finitely nondegenerate at the
origin. As above, fix two radii $\rho_1$ and $\rho_2$ with
$0<\rho_2<\rho_1$. Then there exists an {\rm even} integer
$\kappa_0\in\N_*$ which depends only on the local geometry of $M$ near
the origin, there exists $\varepsilon>0$, there exists $\rho_4>0$ with
$\rho_4< \rho_2$, there exists a complex algebraic $\C^n$-valued
mapping $H(t,J^{\kappa_0})$ which is defined for $t\in\C^n$ with
$\vert t\vert < \rho_4$ and for $J^{\kappa_0}\in\C^{N_{n,\kappa_0}}$
$($where $N_{n,\kappa_0}=n{(n+\kappa_0)!\over n!\ \kappa_0!})$
with $\vert J^{\kappa_0} - J_{\rm Id}^{\kappa_0}\vert < \varepsilon$
and which depends only on $M$ and there exists a geometrically smooth
real algebraic totally real submanifold $E$ of
$\C^{N_{n,\kappa_0}}$ passing through the identity jet $J_{\rm
Id}^{\kappa_0}$ which depends only on $M$, which is defined by
\def\theequation{4.5}\begin{equation}
E=\{J^{\kappa_0}:\, \vert
J^{\kappa_0}-J_{\rm Id}^{\kappa_0}\vert < \varepsilon, \, 
C_l(J^{\kappa_0},\overline{J^{\kappa_0}})=0, \, 
l=1,\dots,\upsilon\},
\end{equation}
where the $C_l(J^{\kappa_0}, \overline{J^{\kappa_0}})$,
$l=1,\dots,\upsilon$, are real algebraic functions defined on the
polydisc $\{\vert J^{\kappa_0} - J_{\rm Id}^{\kappa_0}\vert <
\varepsilon\}$, and which can be constructed
algorithmically by means only of the
defining equations of $M$, such that the following six statements hold{\rm :}
\begin{itemize}
\item[{\bf (1)}]
Every local biholomorphic self-mapping $h\in\mathcal{H}_{M,\kappa_0,
\varepsilon}^{\rho_2,\rho_1}$ of $M$ $($which is defined on the
large polydisc $\Delta_n(\rho_1)$$)$ is represented by
\def\theequation{4.6}\begin{equation}
h(t)=H(t,J^{\kappa_0} h(0)), 
\end{equation}
 on the smallest polydisc $\Delta_n(\rho_4)$. In
particular, each $h\in \mathcal{H}_{M, \kappa_0, \varepsilon}^{
\rho_2,\rho_1}$ is a complex algebraic biholomorphic
mapping. Furthermore, the $\kappa_0$-jet of $h$ at the origin belongs
to the real algebraic submanifold $E$, namely we have
$C_l(J^{\kappa_0} h(0), J^{\kappa_0} \bar h(0))=0$,
$l=1,\dots,\upsilon$.
\item[{\bf (2)}] 
Conversely, shrinking $\varepsilon$ if necessary, given an arbitrary
jet $J^{\kappa_0}$ in $E$ there exists a smaller positive radius
$\rho_5<\rho_4$ such that the mapping defined by
$h(t):=H(t,J^{\kappa_0})$ for $\vert t\vert <\rho_5$ sends $M\cap
\Delta_n(\rho_5)$ CR-diffeomorphically onto its image which is
contained in $M\cap \Delta_n(\rho_4)$. We may therefore say that the
set $\mathcal{H}_{M,\kappa_0, \varepsilon}^{\rho_2,\rho_1}$ of local
biholomorphic self-mappings of $M$ is parametrized by the real
algebraic submanifold $E$.
\item[{\bf (3)}] 
For every choice of two smaller positive radii 
$\widetilde{\rho}_1\leq \rho_1$ and $\widetilde{\rho}_2\leq \rho_2$ with 
$\widetilde{\rho}_2 < \widetilde{\rho}_1$, 
there exists a positive radius $\widetilde{\rho}_4\leq\rho_4$ 
with $\widetilde{\rho}_4<\widetilde{\rho}_2$, and a
positive $\widetilde{\varepsilon}\leq \varepsilon$ such that the
same complex algebraic mapping $H(t,J^{\kappa_0})$ as in statement
{\bf (1)} above represents all local biholomorphic
self-mappings $\widetilde{h}\in \mathcal{H}_{M, \kappa_0,
\widetilde{\varepsilon}}^{\widetilde{\rho}_2,\widetilde{\rho}_1}$ of
$M$, namely we have $\widetilde{h}(t)=
H(t,J^{\kappa_0}\widetilde{h}(0))$ for all $\vert t\vert <
\widetilde{\rho}_4$ as in~\thetag{4.6}. 
Furthermore, the corresponding real algebraic
totally real submanifold $\widetilde{E}$ coincides with $E$ in the
polydisc $\{\vert J^{\kappa_0}- J_{\rm Id}^{\kappa_0}\vert <
\widetilde{\varepsilon}\}$ and it is defined by the same real
algebraic equations $C_l(J^{\kappa_0}, \overline{J^{\kappa_0}})=0$,
$l=1,\dots,\upsilon$, as in equation~\thetag{4.5}. In fact, the
algebraic mapping $H(t,J^{\kappa_0})$ and the real algebraic totally
real submanifold $E$ depend only on the local geometry of $M$ in
a neighborhood of the origin, namely on the germ of $M$ at $0$.
\item[{\bf (4)}] 
The set $\mathcal{H}_{M,\kappa_0,\varepsilon}^{\rho_2,\rho_1}$,
equipped with
the law of composition of holomorphic mappings, is a real
algebraic local Lie group. More precisely,
let the positive integer $c_0$ denote the real
dimension of $E$, which is independent
of $\rho_1, \,\rho_2$ and consider a parametrization
\def\theequation{4.7}\begin{equation}
\R^{c_0}\ni e=(e_1,\dots,e_{c_0})\mapsto j_{\kappa_0}(e)\in E\subset
\C^{N_{n,\kappa_0}}
\end{equation} 
of the real algebraic totally real submanifold $E$. Then there exist a
real algebraic associative local multiplication mapping $(e,e')\mapsto
\mu(e,e')$ and a real algebraic local inversion mapping $e\mapsto
\iota(e)$ such that if we define $H(t;e):=H(t,j_{\kappa_0}(e))$,
then $H(H(t;e);e')\equiv H(t;\mu(e,e'))$ and $H(t;e)^{-1}\equiv
H(t;\iota(e))$, with the local Lie transformation group 
axioms, as defined in {\rm \S2.3}, being satisfied by $H$, $\mu$
and $\iota$.
\item[{\bf (5)}]
For $i=1,\dots,c_0$, consider the one-parameter families of
transformations defined by $H(t;0,\dots,0,e_i,0,\dots,0)=:
H_i(t;e_i)=:H_{i,e_i}(t)$. Then for each $i=1,\dots,c_0$, the vector
field $X_i\vert_{(t;e_i)}:= [\partial_{i,e_i}
H_{e_i}(t')]_{t'=H_{e_i}^{-1}(t)}$, is defined for
$t\in\Delta_n(\rho_5)$ and $\vert e_i\vert < \varepsilon$, has
algebraic coefficients depending on the ``time'' parameter $e_i$, and
has an algebraic flow, since this coincides with the algebraic
mapping $(t,e_i)\mapsto H_{i,e_i}(t)$.
\item[{\bf (6)}]
Let $\rho_5$ be as in statement {\bf (2)}. Then the dimension $c_0$ of
the real Lie algebra $\mathfrak{Hol}(M,\Delta_n(\rho_5))$ is finite, 
bounded by the fixed integer $N_{n,\kappa_0}:= n{(n+\kappa_0)!\over
n!  \ \kappa_0!}$. Furthermore, each vector field $X\in
\mathfrak{Hol}(M,\Delta_n(\rho_5))$ has complex algebraic
coefficients.
\end{itemize}
If $M$ is real analytic, the same theorem holds
with the word ``algebraic'' replaced everywhere by the word
``analytic''.
\end{theorem}

We shall explain below how the integer $\kappa_0$ is related
to the minimality and to the finite nondegeneracy of $M$ at the origin.
The next \S5 and \S6 are devoted to the proof Theorem 4.1, namely
the existence
of the mapping $H(t,J^{\kappa_0})$, the existence of
the real algebraic totally real submanifold $E$ and the completion of 
the proof of properties {\bf (1-6)}.

\section*{\S5.~Minimality and finite nondegeneracy}

\subsection*{5.1.~Local CR geometry of complexified
real analytic generic submanifolds} Let $\zeta\in\C^m$ and
$\xi\in\C^d$ denote some independent coordinates corresponging to
the complexification of the variables $\bar z$ and $\bar w$, which we
denote symbolically by $\zeta:=(\bar z)^c$ and $\xi:= (\bar w)^c$,
where the letter ``c'' stands for the word ``complexified''. We also
write $\tau:=(\bar t)^c$, so $\tau=(\zeta,\xi)\in\C^n$. The {\it
extrinsic complexification}\, $\mathcal{M}:=(M)^c$ of $M$ is the
complex submanifold of codimension $d$ defined by
\def\theequation{5.1}\begin{equation}
\mathcal{M}:=\{(z,w,\zeta,\xi)\in\Delta_n(\rho_1)\times
\Delta_n(\rho_1):\, 
\xi=\Theta(\zeta,z,w)\}.
\end{equation}  
If $M$ is (real, Nash) algebraic, so is $\mathcal{M}$.
 As remarked, we can choose the
equivalent defining equation $w=\overline{\Theta}(z,\zeta,\xi)$ for
$\mathcal{M}$.  In the remainder of \S5, we shall essentially deal
with $\mathcal{M}$ instead of $M$. In fact, $M$ clearly
imbeds in $\mathcal{M}$ as the intersection of $\mathcal{M}$ with the
antiholomorphic diagonal
$\underline{\Lambda}:=\{(t,\tau)\in\C^n\times \C^n: \, \tau=\bar t\}$.

Following [Me1998], [Me2001], we shall complexify a conjugate pair of
generating families of CR vector fields tangent to $M$, namely
$L_1,\dots,L_m$ of type $(1,0)$ and their conjugates
$\overline{L}_1,\dots,\overline{L}_m$ which are of type $(0,1)$.
Here, we can explicitely choose the generators
$L_k=\partial/\partial z_k+\sum_{j=1}^d\, [\partial
\overline{\Theta}_j /\partial z_k (z, \bar z, \bar w)] \,
\partial/\partial w_j$ for $k=1,\dots,m$. Then their complexification
yields a pair of collections of $m$ vector fields defined over
$\Delta_n(\rho_1)\times \Delta_n(\rho_1)$ by
\def\theequation{5.2}\begin{equation}
\left\{
\aligned
\mathcal{L}_k:= & \ 
{\partial \over\partial z_k}+
\sum_{j=1}^d\,
{\partial\overline{\Theta}_j\over \partial z_k}
(z,\zeta,\xi)\,
\, {\partial\over\partial w_j}, \ \ \ \ 
k=1,\dots,m,\\
\underline{\mathcal{L}}_k:= & \
{\partial\over\partial \zeta_k}+\sum_{j=1}^d\,
{\partial \Theta_j\over \partial \zeta_k}(\zeta,z,w)\,
{\partial\over\partial \xi_j}, \ \ \ \
k=1,\dots,m.
\endaligned\right.
\end{equation}
The reader may check directly that $\mathcal{L}_k(w_j-
\overline{\Theta}_j (z,\zeta,\xi)) \equiv 0$, which shows that the
vector fields $\mathcal{L}_k$ are tangent to $\mathcal{M}$.
Similarly, $\underline{\mathcal{L}}_k (\xi_j- \Theta_j(\zeta,z,w))
\equiv 0$, so the vector fields $\underline{\mathcal{L}}_k$ are also
tangent to $\mathcal{M}$. Furthemore, we may check the commutation
relations $[\mathcal{L}_k,\mathcal{L}_{k'}]=0$ and
$[\underline{\mathcal{L}}_k,\underline{\mathcal{L}}_{k'}]=0$ for all
$k, k'=1,\dots,m$. It follows from the Frobenius theorem that the two
$m$-dimensional distributions spanned by each of these two collections
of $m$ vector fields has the integral manifold property. This is not
surprising since the vector fields $\mathcal{L}_k$
are the vector fields tangent to the intersection of
$\mathcal{M}$ with the sets $\{\tau=\tau_p=ct.\}$, which are clearly
$m$-dimensional complex integral manifolds. Following [Me1998], [Me2001], we
denote these manifolds by $\mathcal{S}_{\tau_p}:= \{(t,\tau_p): \,
w=\overline{\Theta}(z,\zeta_p,\xi_p)\}$, where $\tau_p$ is a constant,
and we call them {\it complexified Segre varieties}.  Similarly, the
integral manifolds of the vector fields $\underline{\mathcal{L}}_k$
are the {\it conjugate complexified Segre varieties}\,
$\underline{\mathcal{S}}_{t_p}:=\{(t_p,\tau):\,
\xi=\Theta(\zeta,z_p,w_p)\}$, where $t_p$ is fixed.  The union of
the manifolds $\mathcal{S}_{\tau_p}$ induces a local complex algebraic 
foliation
$\mathcal{F}$ of $\mathcal{M}$ by $m$-dimensional leaves. Similarly,
there is a second foliation $\underline{\mathcal{F}}$ whose leaves are
the $\underline{\mathcal{S}}_{t_p}$.

The following symbolic picture summarizes our constructions. However,
we warn the reader that the codimension $d\geq 1$ of the union
of the two foliations $\mathcal{F}$ and $\underline{\mathcal{F}}$ in 
$\mathcal{M}$ is not visible in this two-dimensional figure.

\bigskip
\begin{center}
\input complexification.pstex_t
\end{center}
\bigskip

Now, we introduce the ``multiple'' flows of the two collections of
conjugate vector fields $(\mathcal{L}_k)_{1\leq k\leq m}$ and
$(\underline{\mathcal{L}}_k)_{1\leq k\leq m}$. For an
arbitrary point $p=(w_p,z_p,\zeta_p,\xi_p)\in\mathcal{M}$ and for an
arbitrary complex ``multitime'' parameter $z_1=(z_{1,1},\dots,
z_{1,m})\in\C^m$, we define
\def\theequation{5.3}\begin{equation}
\left\{
\aligned
\mathcal{L}_{z_1}(z_p,w_p,\zeta_p,\xi_p) & \ :=
\exp(z_1\mathcal{L})(p):=
\exp(z_{1,1}\mathcal{L}_1(\cdots(\exp(z_{1,m}\mathcal{L}_m(
p)))\cdots)):=\\
& \
:=(z_p+z_1,\overline{\Theta}(z_p+z_1,\zeta_p,\xi_p),
\zeta_p,\xi_p).
\endaligned\right.
\end{equation}
With this formal definition, there exists a maximal connected open
subset $\Omega$ of $\mathcal{M}\times \C^m$ containing
$\mathcal{M}\times\{0\}$ such that
$\mathcal{L}_{z_1}(p)\in\mathcal{M}$ for all
$(z_1,p)\in\Omega$. Analogously, for $(\zeta_1,p)$ running in a
similar open subset $\underline{\Omega}$, we may also define
the map
\def\theequation{5.4}\begin{equation}
\underline{\mathcal{L}}_{\zeta_1}(z_p,w_p,\zeta_p,\xi_p):=
(z_p,w_p,\zeta_p+\zeta_1,\Theta(\zeta_p+\zeta_1,z_p,w_p)).
\end{equation}
We notice that the two maps given by ~\thetag{5.3} and~\thetag{5.4} are
holomorphic in their variables. Since $M$ is real algebraic, they are
moreover complex algebraic.  

\subsection*{5.2.~Segre chains}
Let us start from the point $p$ being the
origin and let us move alternately in the direction of $\mathcal{S}$ or of
$\underline{\mathcal{S}}$, namely we consider the two maps
$\Gamma_1(z_1):=\mathcal{L}_{z_1}(0)$ and
$\underline{\Gamma}_1(z_1):=
\underline{\mathcal{L}}_{z_1}(0)$. Next, we start from these
endpoints and we move in the other direction, namely, we consider the
two maps
\def\theequation{5.5}\begin{equation}
\Gamma_2(z_1,z_2):=\underline{\mathcal{L}}_{z_2}(\mathcal{L}_{z_1}(0)), 
\ \ \ \ \
\underline{\Gamma}_2(z_1,z_2):=
\mathcal{L}_{z_2}(\underline{\mathcal{L}}_{z_1}(0)),
\end{equation}
where $z_1,\,z_2\in\C^m$. Also, we define $\Gamma_3(z_1,z_2,z_3):=
\mathcal{L}_{z_3}(\underline{\mathcal{L}}_{z_2}(\mathcal{L}_{z_1}(0)))$,
{\it etc.} By induction, for every positive integer $k$, we obtain two maps 
$\Gamma_k(z_1,\dots,z_k)$ and
$\underline{\Gamma}_k(z_1,\dots,z_k)$. In the sequel, we shall often
use the notation $z_{(k)}:=(z_1,\dots,z_k)\in\C^{mk}$.  Since
$\Gamma_k(0)=\underline{\Gamma}_k(0)=0$, for every $k\in\N_*$, there
exists a sufficiently small open polydisc $\Delta_{mk}(\delta_k)$
centered at the origin in $\C^{mk}$ with $\delta_k>0$ such that
$\Gamma_k(z_{(k)})$ and $\underline{\Gamma}_k(z_{(k)})$ belong to
$\mathcal{M}$ for all $z_{(k)} \in \Delta_{mk}(\delta_k)$.

We also exhibit a simple link between the maps $\Gamma_k$ and
$\underline{\Gamma}_k$. Let $\sigma$ be the antiholomorphic involution
defined by $\sigma(t,\tau):=(\bar\tau,\bar t)$.  Since
$w=\overline{\Theta}(z,\zeta,\xi)$ if and only if
$\xi=\Theta(\zeta,z,w)$, this involution maps $\mathcal{M}$ onto
$\mathcal{M}$ and it also fixes the antidiagonal $\underline{\Lambda}$
pointwise.  Using the definitions~\thetag{5.3} and~\thetag{5.4}, we
see readily that $\sigma(\mathcal{L}_{z_1}(0)) =
\underline{\mathcal{L}}_{\bar z_1}(0)$.  It follows generally that $\sigma(
\Gamma_k(z_{(k)}))= \underline{\Gamma}_k( \overline{z_{(k)}})$.

Next, we observe that $\Gamma_{k+1}(z_{(k)},0)=\Gamma_k(z_{(k)})$,
since $\mathcal{L}_0$ and $\underline{\mathcal{L}}_0$ coincide with
the identity map. So the ranks at the origin of
the maps $\Gamma_k$ increase with $k$.

\def\thedefinition{5.1}\begin{definition}
{\rm The real analytic generic manifold $M$ is said to 
be {\it minimal}\, at $p$ if the maps $\Gamma_k$ are 
of (maximal possible) rank equal to $2m+d={\rm dim}_\C\,\mathcal{M}$ at
the origin in $\Delta_{mk}(\delta_k)$ for all $k$ large enough.} 
\end{definition}

The following fundamental properties are established in [Me1998],
[Me2001].

\def\thetheorem{5.2}\begin{theorem}
The minimality of $M$ at $0$ is a biholomorphically invariant property.
It depends neither on the choice of a defining equation for $M$ nor on
the choice of a system of generating complexified CR vector fields
$(\mathcal{L}_k)_{1\leq k\leq m}$ and
$(\underline{\mathcal{L}}_k)_{1\leq k\leq m}$. Also, minimality is
equivalent to the fact that the Lie algebra generated by the
complexified CR vector fields $(\mathcal{L}_k)_{1\leq k\leq m}$ and
$(\underline{\mathcal{L}}_k)_{1\leq k\leq m}$ spans $T\mathcal{M}$ in
a neighborhood of $0$. Furthermore, there exists an invariant integer
$\nu_{0}$, called the {\rm Segre type}\, of $M$ at $0$ satisfying
$\nu_{0}\leq d+1$ which is the smallest integer such that the mappings
$\Gamma_k$ and $\underline{\Gamma}_k$ are of generic rank equal
to $2m+d$ over $\Delta_{mk}(\delta_k)$ for all $k\geq
\nu_{0}+1$. Finally, with this integer $\nu_0$, the odd integer
$\mu_0:=2\nu_0+1$, called the {\rm Segre type} $\mathcal{M}$ at $0$ is
the smallest integer such that the mappings $\Gamma_k$ and
$\underline{\Gamma}_k$ are of rank equal to $2m+d$ at the origin
in $\Delta_{mk}(\delta_k)$.
\end{theorem}

Let $\mu_{0}:=2\nu_{0}+1$ be the Segre type of $\mathcal{M}$ at
$0$ (notice that this is always odd). In the remainder of this section,
we assume that $M$ is minimal at $0$ and we exploit
the rank condition on $\Gamma_k$. More precisely we choose a
positive $\eta$ with $0<\eta\leq \delta_{\mu_{0}}$ such that
$\Gamma_{\mu_{0}}$ has rank $2m+d$ at every point of the polydisc
$\Delta_{m\mu_{0}}(\eta)$. Without loss of generality, we can also
assume that $\Gamma_{\mu_{0}}(\Delta_{m\mu_{0}}(\eta))$ contains
$\mathcal{M}\cap (\Delta_n(\rho_4)\times \Delta_n(\rho_4))$. Simple
examples in the hypersurface case show that $\rho_4<<\rho_1$ and in
fact, one has necessarily 
an inequality of the form $\rho_4\leq (\rho_1)^N$, where $N$ is a
certain integer depending on the vector fields $(\mathcal{L}_k)_{1\leq
k\leq m}$ and $(\underline{\mathcal{L}}_k)_{1\leq k\leq m}$ ({\it
cf.}~[Be1996]).

\subsection*{5.3.~Finite nondegeneracy}
The last ingredient for Theorem~4.1 consists in developing
the equations of $M$ in powers of $\bar z$ as follows
\def\theequation{5.6}\begin{equation}
\bar w_j=\sum_{\beta\in\N^m}\, (\bar z)^\beta\,
\Theta_{j,\beta}(t), \ \ \ \ \
j=1,\dots,d,
\end{equation}
where the functions $\Theta_{j,\beta}(t)$ are holomorphic in the
polydisc $\Delta_n(2\rho_1)$.  So we may introduce the holomorphic
maps $\psi_k(t):= (\Theta_{j,\beta}(t))_{1\leq j\leq d, \,
\vert\beta\vert \leq k}$ with values in $\C^{d{(m+k)!\over m!\ k!}}$.
Obviously, the ranks at the origin of the $\psi_k$ increase with $k$.

\def\thedefinition{5.3}\begin{definition}
{\rm The generic manifold $M$ is said to be {\it finitely nondegenerate}\,
at $0$ if there exists a positive integer $k$ such that the 
rank at the origin of the map $\psi_k$ is equal to $n$.}
\end{definition}

It may be checked that this definition depends neither on the system
of coordinates nor on the choice of a collection of $d$ defining
equations for $M$ and that it coincides with the definition given
in \S1.2. If $M$ is finitely nondegenerate at $0$ we denote by $\ell_{0}$ 
the
smallest integer $k$ given by definition 5.3 and we say that
$M$ is $\ell_{0}$-nondegenerate at the origin.

Finite nondegeneracy is interesting for the following reason.
In the sequel, we shall have to consider an infinite collection of
equations of the form
\def\theequation{5.7}\begin{equation}
\Theta_{j,\beta}(t)+\sum_{\gamma\in\N_*^m}\,
(\zeta)^\gamma\, {(\beta+\gamma)!\over \beta!\ \gamma!}\,
\Theta_{j,\beta+\gamma}(t)=\omega_{j,\beta},
\end{equation}
where $\N_*^m:=\N^m\backslash \{0\}$, where $j$ runs from $1$ to $d$,
where $\beta$ runs in $\N^m$ and where the right hand sides
$\omega_{j,\beta}$ are independent complex variables. For $\beta=0$,
the equations~\thetag{5.7} write simply
$\Theta_j(\zeta,t)=\omega_{j,0}$. By definition, if $M$ is
$\ell_{0}$-nondegenerate at $0$, there exists $n$ integers
$j_*^1,\dots,j_*^n$ with $1\leq j_*^i\leq d$ and $n$ multi-indices
$\beta_*^1,\dots,\beta_*^n\in\N^m$ with $\vert\beta_*^i\vert\leq
\ell_{0}$ such that the local holomorphic self-mapping $t\mapsto
(\Theta_{j_*^k,\beta_*^k}(t))_{1\leq k\leq n}$ of $\C^n$ is of rank
$n$ at the origin. Considering the equations~\thetag{5.7} for
$j=j_*^1,\dots,j_*^n$ and $\beta=\beta_*^1,\dots,\beta_*^n$ and
applying the implicit function theorem, we observe that we can solve
$t$ in terms of
$(\zeta,\omega_{ j_*^1,\beta_*^1}, \dots,\omega_{j_*^n, \beta_*^n})$ by
means of a holomorphic mapping, namely
\def\theequation{5.8}\begin{equation}
t=\Psi(\tau,\omega_{j_*^1,\beta_*^1},\dots,\omega_{j_*^n,\beta_*^n}).
\end{equation}
Without loss of generality, we may assume that $\Psi$ is 
holomorphic for $\vert \zeta\vert < \widetilde{\rho}_3$ 
and $\vert \omega_{j_*^i,
\beta_*^i} \vert < \widetilde{\rho}_3$, 
where $0<\widetilde{\rho}_3<\rho_2<\rho_1$.

\section*{\S6.~Algebraicity of local CR automorphism groups}

\subsection*{6.1.~Fundamental reflection identity for the mapping}
So $M$ is $\ell_0$-nondegenerate at the origin. Recall that $\mu_0=2\nu_0+1$
is the Segre type of $\mathcal{M}$ and introduce the new integer
$\kappa_0:=\ell_0(\mu_0+1)$. Notice that
$\kappa_0$ is even. Let us 
take an arbitrary local holomorphic self map $h$ of $M$ close to the
identity in the set $\mathcal{H}_{M,\kappa_0, \varepsilon}^{
\rho_2,\rho_1}$, {\it i.e.} with $k:=\kappa_0$ in the
definition~\thetag{4.4}. We denote the map $h$ by
$(h_1,\dots,h_n)=(f_1,\dots,f_m,g_1,\dots,g_d)$, according to the
splitting $t=(z,w)$ of the coordinates. The complexification
$h^c:=(h,\bar h)$ induces a local holomorphic self map of the
complexification $\mathcal{M}$. More precisely, for all
$(t,\tau)\in\mathcal{M}$ with $\vert t\vert, \, \vert \tau\vert \leq
\rho_2$, we have $(h(t),\bar h(\tau))\in\mathcal{M}$ and $\vert
h(t)\vert,\, \vert h(\tau)\vert <\rho_1$, so we can write
\def\theequation{6.1}\begin{equation}
\bar g_j(\tau)=\Theta_j(\bar f(\tau),h(t)),
\end{equation}
for $j=1,\dots,d$.  Since $h$ is a biholomorphism and
$T_0^cM=\{w=0\}$, it follows that the determinant
\def\theequation{6.2}\begin{equation}
{\rm det}\, (\underline{\mathcal{L}}_k \bar f_l(\tau))_{1\leq k,l\leq n},
\end{equation}
which is a $\K$-analytic function of $(t,\tau)\in\mathcal{M}$, does not
vanish at the origin. Shrinking $\varepsilon$ if necessary, we can
assume that for every holomorphic map $h\in\mathcal{H}_{M,\kappa_0,\varepsilon
}^{\rho_2,\rho_1}$, the determinant~\thetag{6.2} does
not vanish for all $\vert t\vert, \, \vert\tau\vert <\rho_2$. We now
differentiate~\thetag{6.2} by applying the vector fields
$\underline{\mathcal{L}}_1,\dots, \underline{\mathcal{L}}_m$, which
gives
\def\theequation{6.3}\begin{equation}
\underline{\mathcal{L}}_k \bar g_j(\tau)=
\sum_{l=1}^m\,
{\partial\Theta_j\over
\partial \zeta_l}(\bar f(\tau),h(t))\,
\underline{\mathcal{L}}_k \bar f_l(\tau),
\end{equation}
for $k=1,\dots,m$ and $j=1,\dots,d$. For fixed $j$, we consider the
$m$ equations~\thetag{6.3} as an affine system satisfied by
the partial derivatives $\partial\Theta_j/\partial \zeta_l$.
By Cramer's rule, there exists universal polynomials
$\Omega_{j,k}$ in their variables such that
\def\theequation{6.4}\begin{equation}
{\partial\Theta_j\over
\partial \zeta_k}(\bar f(\tau),h(t))=
{\Omega_{j,k}(\{\underline{\mathcal{L}}_l\,
\bar h(\tau)\}_{1\leq l\leq m})
\over
{\rm det}\, (\underline{\mathcal{L}}_k 
\bar f_l(\tau))_{1\leq k,l\leq n}}
\end{equation}
for all $(t,\tau)\in\mathcal{M}$ with $\vert t\vert,\,\vert\tau\vert <
\rho_2$ and for $k=1,\dots,m$, $j=1,\dots,d$.  

Applying the derivations $\underline{\mathcal{L}}_k$
to~\thetag{6.4} we see by induction that for every multi-index
$\beta\in\N_*^m$ and for every $j=1,\dots,d$, there exists a universal
polynomial $\Omega_{j,\beta}$ in its variables such that
\def\theequation{6.5}\begin{equation}
{1\over\beta!}\,
{\partial^{\vert\beta\vert}\Theta_j\over
\partial \zeta^\beta}(\bar f(\tau),h(t))=
{\Omega_{j,\beta}(\{\underline{\mathcal{L}}^\gamma\,
\bar h(\tau)\}_{\vert\gamma\vert\leq\vert\beta\vert})
\over
[{\rm det}\, (\underline{\mathcal{L}}_k 
\bar f_l(\tau))_{1\leq k,l\leq n}]^{2\vert\beta\vert-1}},
\end{equation}
for all $(t,\tau)\in\mathcal{M}$ with $\vert t\vert,\,\vert\tau\vert <
\rho_2$.  Here, for $\gamma=(\gamma_1,\dots, \gamma_m)\in\N^m$, we
denote by $\underline{\mathcal{L}}^\gamma$ the derivation
$(\underline{\mathcal{L}}_1)^{\gamma_1}\dots
(\underline{\mathcal{L}}_m)^{\gamma_m}$.
Next, denoting by $\omega_{j,\beta}(t,\tau)$
the right hand side of~\thetag{6.5} and developing the left hand side
in power series using~\thetag{5.7}, we may write
\def\theequation{6.6}\begin{equation}
\Theta_{j,\beta}(h(t))+\sum_{\gamma\in\N_*^m}\,
(\bar f(\tau))^\gamma\,\Theta_{j,\beta+\gamma}(h(t))=
\omega_{j,\beta}(t,\tau).
\end{equation}
Recall that $M$ is $\ell_{0}$-nondegenerate at $0$.
Using~\thetag{5.8}, we can solve $h(t)$ in terms
of the derivatives of $\bar h(\tau)$, namely
\def\theequation{6.7}\begin{equation}
\left\{
\aligned
h(t)= 
& \
\Psi\left(\bar f(\tau),
{\Omega_{j_*^1,\beta_*^1}(\{\underline{\mathcal{L}}^\gamma\,
\bar h(\tau)\}_{\vert\gamma\vert\leq\vert\beta_*^1\vert})\over
[{\rm det}\, (\underline{\mathcal{L}}_k 
\bar f_l(\tau))_{1\leq k,l\leq n}]^{2\vert\beta_*^1\vert-1}},
\dots\right.\\
& \ \ \ \ \ \ \ \ \ \ \ \ \ \ 
\left.
\dots,
{\Omega_{j_*^n,\beta_*^n}(\{\underline{\mathcal{L}}^\gamma\,
\bar h(\tau)\}_{\vert\gamma\vert\leq\vert\beta_*^n\vert})\over
[{\rm det}\, (\underline{\mathcal{L}}_k 
\bar f_l(\tau))_{1\leq k,l\leq n}]^{2\vert\beta_*^n\vert-1}}\right)=\\
= & \ \Psi(\bar f(\tau), \omega_{j_*^1,\beta_*^1}(t,\tau), \dots,
\omega_{j_*^n,\beta_*^n}(t,\tau)).
\endaligned\right.
\end{equation} 

\noindent
Here, the maximal length of the multi-indices
$\beta_*^1,\dots,\beta_*^n$ is equal to $\ell_{0}$. According
to~\thetag{5.8}, the representation~\thetag{6.7} of $h(t)$ holds
provided $\vert \bar g(\tau)\vert < \widetilde{\rho}_3$ and $\vert
\omega_{j_*^i, \beta_*^i}\vert < \widetilde{\rho}_3$.  Since the
coordinates are normal, we have $\Theta_j(\bar z,0,0)\equiv 0$, or
equivalently $\Theta_{j,\beta}(0)=0$ for all $j=1,\dots,d$ and all
$\beta\in\N^m$.  It follows from~\thetag{6.6} and from $h(0)=0$ that
$\omega_{j,\beta}(0)=0$, for all $j=1,\dots,d$ and all
$\beta\in\N^m$. Consequently, there exists a radius $\rho_3\sim
\widetilde{\rho}_3$ with $0<\rho_3<\rho_2<\rho_1$ such that $\vert
\omega_{j_*^i,\beta_*^i}(t,\tau)\vert < \widetilde{\rho}_3$, $i=1,\dots,n$ and
such that $\vert \bar g(\tau)\vert < \widetilde{\rho}_3$ for all
$h\in\mathcal{H}_{M,\kappa_0, \varepsilon}^{\rho_2,\rho_1}$ and for
all $(t,\tau)\in\mathcal{M}$ with $\vert t\vert, \, \vert \tau\vert
<\rho_3$.

In conclusion, the relation~\thetag{6.7} holds for all
$h\in\mathcal{H}_{M,\kappa_0, \varepsilon}^{\rho_2,\rho_1}$ and for
all $(t,\tau)\in\mathcal{M}$ with $\vert t\vert, \, \vert \tau\vert
<\rho_3$.

Next, using the explicit expressions of
the vector fields $\underline{\mathcal{L}}_k$ given in~\thetag{5.2},
we may develop the higher order derivatives
$\underline{\mathcal{L}}^\gamma\bar h(\tau)$ as polynomials in the
$\vert\gamma\vert$-jet $(\partial_\tau^{\gamma'}\bar
h(\tau))_{\vert\gamma'\vert\leq \vert\gamma\vert}$ of $\bar h(\tau)$
with coefficients being certain holomorphic functions of $(t,\tau)$
obtained as certain polynomials with respect to the partial
derivatives of the functions $\Theta_j(\zeta,t)$.

To be more explicit in this desired new representation of~\thetag{6.7}, 
we remind first our jet notation. For each
$i=1,\dots,n$ and each $\alpha\in\N^n$, we introduced a new {\it
independent}\, coordinate $J_i^\alpha$ corresponding to the
partial derivative $\partial_\tau^\alpha\bar h_i(\tau)$ (or
$\partial_t^\alpha h_i(t)$). The space of
$k$-jets of holomorphic mappings $\bar h(\tau)$ is then the complex space
$\C^{n{(n+k)!\over n!\ k!}}$ with coordinates
$(J_i^\alpha)_{1\leq i\leq n,\,\vert\alpha\vert\leq k}$. It will
be convenient to use the abbreviations
$J^k:=(J_i^\alpha)_{1\leq i\leq n,\,\vert\alpha\vert\leq
k}$ and $J^k\bar h(\tau):= (\partial_\tau^\alpha\bar h_i(\tau))_{
1\leq i\leq n,\,\vert\alpha\vert\leq k}$. 

So pursuing with~\thetag{6.7}, we argue that for
every $\gamma\in\N^m$, there exists a polynomial in the jet
$J^{\vert\gamma\vert}\bar h(\tau)$ with holomorphic cooeficients
depending only on $\Theta$ such that
\def\theequation{6.8}\begin{equation}
\underline{\mathcal{L}}^\gamma\,\bar h(\tau)\equiv
P_\gamma(t,\tau,J^{\vert\gamma\vert}\bar h(\tau)).
\end{equation}
Putting all these expressions in~\thetag{6.7}, we obtain an important
relation between $h$ and the $\ell_0$-jet
of $\bar h$ which we may now summarize.
At first, as $\kappa_0=\ell_0(\mu_0+1)\geq \ell_0$,
observe that for every
$h\in\mathcal{H}_{M,\kappa_0,\varepsilon}^{\rho_2,\rho_1}$, we have
$\vert\vert J^{\ell_0}h-J_{\rm Id}^{\ell_0}\vert\vert_{\rho_2}\leq
\vert\vert J^{\kappa_0}h-J_{\rm Id}^{\kappa_0}\vert\vert_{\rho_2}\leq
\varepsilon$. 
Shrinking $\varepsilon$ if
necessary, we have proved the following lemma.

\def\thelemma{6.1}\begin{lemma}
There exists a complex algebraic $\C^n$-valued mapping $\Pi(t,\tau,
J^{\ell_{0}})$ defined for $\vert t\vert,\, \vert\tau\vert < \rho_3$
and for $\vert J^{\ell_{0}}- J_{\rm Id}^{\ell_{0}} \vert<\varepsilon$
which depends only on the defining functions
$\xi_j-\Theta_j(\zeta,t)$ of $\mathcal{M}$, such that for every
local holomorphic self-mapping
$h\in \mathcal{H}_{ M,\kappa_0, \varepsilon}^{ \rho_2,\rho_1}$ of $M$
$($hence satisfying $\vert\vert J^{\ell_{0}}h- J_{\rm
Id}^{\ell_{0}} \vert\vert_{\rho_2}< \varepsilon)$, the relation
\def\theequation{6.9}\begin{equation}
h(t)=\Pi(t,\tau,J^{\ell_{0}}\bar h(\tau))
\end{equation}
holds for all $(t,\tau)\in\mathcal{M}$ with $\vert t\vert,\,
\vert \tau\vert <\rho_3$.
\end{lemma}

\subsection*{6.2.~Reflection identity for arbitrary jets}
Let now $\Upsilon_j$ and $\underline{\Upsilon}_j$
be the vector fields tangent to 
$\mathcal{M}$ defined by 
\def\theequation{6.10}\begin{equation}
\Upsilon_j:={\partial\over\partial w_j}+
\sum_{l=1}^d\,\Theta_{l,w_j}(\zeta,t)\,{\partial\over
\partial \xi_l}, \ \ \ \ \
\underline{\Upsilon}_j:=
{\partial\over\partial\xi_j}+\sum_{l=1}^d\,
\overline{\Theta}_{l,\xi_j}(z,\tau)\,
{\partial\over\partial w_l},
\end{equation}
for $j=1,\ldots,d$.
We observe that the collection of $2m+d$ vector fields
$\mathcal{L}_k$, $\underline{\mathcal{L}}_k$, $\Upsilon_j$ span
$T\mathcal{M}$. The same holds for the collection $\mathcal{L}_k$,
$\underline{\mathcal{L}}_k$, $\underline{\Upsilon}_j$. We also
have the commutation relations
$[\Upsilon_j,\underline{\mathcal{L}}_k]=0$ and
$[\underline{\Upsilon}_j,\mathcal{L}_k]=0$.
We observe that $\Upsilon^\gamma h(t)=
\partial_w^\gamma h(t)$ for all $\gamma\in\N^d$.
Let $\alpha=(\beta,\gamma)\in\N^m\times\N^d$. By expanding
$\mathcal{L}^\beta\,\Upsilon^\gamma \, h(t)$ using the explicit
expressions~\thetag{5.2}, we obtain a polynomial 
$Q_{\beta,\gamma}(t,\tau,(\partial_t^{\alpha'}h(t))_{\vert
\alpha'\vert\leq\vert\alpha\vert})$, where $Q_{\beta,\gamma}$ 
is a polynomial in its last variables with coefficients
depending on $\overline{\Theta}$ and its partial derivatives.
Conversely, 
since $\mathcal{L}_k\vert_0=\partial_{z_k}$ at the origin, 
we can invert these formulas, so
there exist polynomials $P_\alpha$ in their last 
variables with coefficients depending only on 
$\overline{\Theta}$ such that
\def\theequation{6.11}\begin{equation}
\partial_t^\alpha h(t)=
P_\alpha(t,\tau,(
\mathcal{L}^{\beta'}\, 
\Upsilon^{\gamma'}\, 
h(t))_{\vert\beta'\vert\leq\vert\beta\vert, \,
\vert\gamma'\vert\leq\vert\gamma\vert}).
\end{equation}

\def\thelemma{6.2}\begin{lemma}
For every $\ell\in\N$, there exists a complex algebraic mapping
$\Pi_\ell$ with values in $\C^{N_{n,\ell}}$ defined for $\vert t\vert, \,
\vert\tau\vert<\rho_3$ and $\vert J^{\ell_0}-J_{\rm
Id}^{\ell_0}\vert <\varepsilon$ which is relatively polynomial with respect
to the higher order jets $J_i^\alpha$ with
$\vert\alpha\vert\geq \ell_0+1$, $i=1,\dots,n$, such that for
every local holomorphic self-mapping $h\in\mathcal{H}_{M,\kappa_0,
\varepsilon}^{\rho_2,\rho_1}$, the two conjugate relations
\def\theequation{6.12}\begin{equation}
\left\{
\aligned
J^\ell h(t)= & \
\Pi_\ell(t,\tau,J^{\ell_0+\ell}\bar h(\tau)),\\
J^\ell\bar h(\tau)= & \
\overline{\Pi}_\ell(\tau,t,J^{\ell_0+\ell} h(t)).
\endaligned\right.
\end{equation}
hold for all
$(t,\tau)\in\mathcal{M}$ with $\vert t\vert, \, 
\vert\tau\vert < \rho_3$.
\end{lemma}

\proof
Applying the derivations $\mathcal{L}^\beta\,
\Upsilon^\gamma$ to~\thetag{6.9}, and using the chain rule, 
we obtain 
\def\theequation{6.13}\begin{equation}
\mathcal{L}^\beta\,\Upsilon^\gamma \, h(t)=
\Pi_{\beta,\gamma}(t,\tau,J^{\ell_0+\vert\beta\vert
+\vert\gamma\vert}\bar h(\tau)), 
\end{equation}
where the function $\Pi_{\beta,\gamma}$ (as the function $\Pi$) is
holomorphic for $\vert t\vert, \, \vert \tau\vert < \rho_3$ and $\vert
J^{\ell_0}-J_{\rm Id}^{\ell_0}\vert < \varepsilon$ and
relatively polynomial with respect to the jets $J_i^\alpha$
with $\vert \alpha \vert \geq \ell_0+1$.  Applying~\thetag{6.11}, we
obtain the function $\Pi_\ell$, which completes the proof.
\endproof

\subsection*{6.3.~Substitutions of reflection identities}
Let $\pi_t(t,\tau):=t$ and $\pi_\tau(t,\tau):=\tau$ denote the two
canonical projections. We write $h^c(t,\tau):=(h(t),\bar h(\tau))$.
We make the following slight abuse of notation: instead of rigorously
writing $h(\pi_t(t,\tau))$, we write $h(t,\tau)=h(t)$ and $\bar
h(t,\tau)= \bar h(\tau)$.

Let $x\in\C^\nu$ and let $\mathcal{Q}(x)=(\mathcal{Q}_1(x),\dots,
\mathcal{Q}_{2n}(x))\in\C\{x\}^{2n}$. As the multiple flow of
$\underline{\mathcal{L}}$ given by~\thetag{5.3} does not act on the
$(z,w)$ variables, we have the trivial but important property
$h(\underline{\mathcal{L}}_{z_1}(\mathcal{Q}(x)))= h(\mathcal{Q}(x))$.
More generally, for every multi-index $\alpha\in\N^n$, we have
$\partial_t^\alpha h(\underline{\mathcal{L}}_{z_1}(\mathcal{Q}(x)))=
\partial_t^\alpha h(\mathcal{Q}(x))$. Analogously, we have
$\partial_\tau^\alpha \bar
h(\mathcal{L}_{z_1}(\mathcal{Q}(x)))=\partial_\tau^\alpha
\bar h(\mathcal{Q}(x))$. Since for $k$ even, we have
$\Gamma_k(z_{(k)})=
\underline{\mathcal{L}}_{z_k}(\Gamma_{k-1}(z_{(k-1)}))$, the
following two properties hold:
\def\theequation{6.14}\begin{equation}
\left\{
\aligned
J^\ell h(\Gamma_k(z_{(k)}))= & \
J^\ell h(\Gamma_{k-1}(z_{(k-1)})), \ \ \ \ \ 
{\sf if} \ k \ {\sf is} \ {\sf even;}\\
J^\ell \bar h(\Gamma_k(z_{(k)}))= & \
J^\ell \bar h(\Gamma_{k-1}(z_{(k-1)})), \ \ \ \ \
{\sf if} \ k \ {\sf is} \ {\sf odd}.
\endaligned\right.
\end{equation}

Let now $\kappa_0:=\ell_0(\mu_0+1)$ be the product of the Levi type
with the Segre type of $\mathcal{M}$ plus $1$ and consider the open
subset of the $\kappa_0$-order jet space $\C^{N_{n,\kappa_0}}$ defined
by the inequality $\vert J^{\kappa_0}-J_{\rm Id}^{\kappa_0}\vert<
\varepsilon$. Let $z_{(k)}\in\Delta_{mk}$ as in \S5.6 above. Since the
maps $\Gamma_k$ are holomorphic and satisfy $\Gamma_k(0)=0$, we may
choose $\delta>0$ sufficiently small in order that the following two
conditions are satisfied for every $k\leq \mu_0$ and for and for every
$\vert z_{(k)}\vert <\delta$:
\def\theequation{6.15}\begin{equation}
\vert\Gamma_k(z_{(k)})\vert <\rho_3 
\ \ \ \ \
{\rm and}
\ \ \ \ \ 
\vert J^{\kappa_0}h(\Gamma_k(z_{(k)}))-
J_{\rm Id}^{\kappa_0}\vert < \varepsilon.
\end{equation}
This choice of $\delta$ is convenient to make several 
susbtitutions by means of formulas~\thetag{6.12}.
The formulas~\thetag{6.16} that we will obtain below strongly differ from the
previous formulas~\thetag{6.12}, because they depend on
the jet of $h$ at the origin only.

\def\thelemma{6.3}\begin{lemma}
Shrinking $\varepsilon$ if necessary, for every integer $k\leq
\mu_0+1$ and for every integer $\ell\geq 0$, there exists a complex
algebraic mapping $\Pi_{\ell,k}$ with values in $\C^{N_{n,\ell}}$
defined for $\vert t\vert, \, \vert\tau\vert <\rho_3$ and for $\vert
J^{k\ell_0}-J_{\rm Id}^{k\ell_0}\vert < \varepsilon$, which is
relatively polynomial with respect to the higher order jets
$J_i^\alpha$ with $\vert\alpha\vert \geq k\ell_0+1$, $i=1,\ldots,n$,
and which depends only on the defining functions
$\xi_j-\Theta_j(\zeta,t)$ of $\mathcal{M}$, such that the following
two families of conjugate identities are satisfied
\def\theequation{6.16}\begin{equation}
\left\{
\aligned
J^\ell h(\Gamma_k(z_{(k)}))= & \ 
\Pi_{\ell,k}(\Gamma_k(z_{(k)}),
J^{k\ell_0+\ell}\bar h(0)), \ \ \ \ \
{\rm if} \ k \ {\sf is} \ {\sf odd};\\
J^\ell \bar h(\Gamma_k(z_{(k)}))= & \
\overline{\Pi_{\ell,k}}(\Gamma_k(z_{(k)}),
J^{k\ell_0+\ell}\bar h(0)), \ \ \ \ \
{\rm if} \ k \ {\sf is} \ {\sf even}.
\endaligned\right.
\end{equation}
\end{lemma}

\proof
For $k=1$, replacing $(t,\tau)$ by $\Gamma_1(z_{(1)})$ in
the first relation~\thetag{6.12}
and using the second property~\thetag{6.14}, we get
\def\theequation{6.17}\begin{equation}
\left\{
\aligned
J^\ell h(\Gamma_1(z_{(1)}))= & \
\Pi_\ell(\Gamma_1(z_{(1)}),
J^{\ell_0+\ell}\bar h(\Gamma_1(z_{(1)})))=\\
= & \ \Pi_\ell(\Gamma_1(z_{(1)}),J^{\ell_0+\ell}\bar h(0)),
\endaligned\right.
\end{equation}
so the lemma holds true for $k=1$ if we simply choose
$\Pi_{\ell,1}:=\Pi_\ell$. By induction, suppose that the lemma holds
true for $k\leq \mu_0$.  To fix the ideas, let us assume that this
$k$ is even (the odd case is completely similar).  Then replacing the
arguments $(t,\tau)$ in the first relation~\thetag{6.12} 
by $\Gamma_{k+1}(z_{(k+1)})$,
using again the second property~\thetag{6.14}, and using the induction
assumption, namely using the conjugate of the second relation~\thetag{6.16}
with $\ell$ replaced by $\ell_0+\ell$, we get
\def\theequation{6.18}\begin{equation}
\left\{
\aligned
J^\ell h(\Gamma_{k+1}(z_{(k+1)})) & \ =
\Pi_\ell(\Gamma_{k+1}(z_{(k+1)}),
J^{\ell_0+\ell}\bar h(\Gamma_{k+1}(z_{(k+1)}))) =\\
& \ =
\Pi_\ell(\Gamma_{k+1}(z_{(k+1)}),J^{\ell_0+\ell}\bar  h
(\Gamma_k(z_{(k)})))= \\
& \ =
\Pi_\ell(\Gamma_{k+1}(z_{(k+1)}),\overline{
\Pi_{\ell_0+\ell,k}}(
\Gamma_k(z_{(k)}),J^{k\ell_0+\ell_0+\ell} \bar h (0)))=: \\
& \ =:
\Pi_{\ell,k+1}(\Gamma_{k+1}(z_{(k+1)}),
J^{(k+1)\ell_0+\ell} \bar h(0)), 
\endaligned\right.
\end{equation}
which yields the desired formula at level $k+1$.
For the above formal composition formulas to be correct, we possibly
have to shrink $\varepsilon$. Finally, a direct inspection of relative
polynomialness shows that $\Pi_{\ell,k+1}$ is polynomial with respect
to the jet variables $J_i^\alpha$ with $\vert\alpha\vert \geq
(k+1)\ell_0+1$, $i=1,\dots,n$. The proof of Lemma~6.21 is complete.
\endproof

\subsection*{6.4.~Algebraic parameterization of 
CR mappings by their jet at the origin} Finally, as in the paragraph
after Theorem~5.2, we choose $\rho_4>0$ sufficiently small such that
$\Gamma_{\mu_{0}}$ maps the polydisc $\Delta_{m\mu_{0}}(\eta)$
submersively onto an open neighborhood of the origin in $\mathcal{M}$
which contains the open subset $\mathcal{M}\cap
(\Delta_n(\rho_4)\times \Delta_n(\rho_4))$.  From the relation
$\Gamma_{\mu_0+1}(z_{(\mu_0)},0)\equiv \Gamma_{\mu_0}(z_{(\mu_0)})$,
it follows trivially that $\Gamma_{\mu_0+1}$ also induces a submersion
from $\Delta_{m(\mu_0+1)}(\eta)$ onto $\mathcal{M}\cap
(\Delta_n(\rho_4)\times \Delta_n(\rho_4))$. It follows that the
composition $\pi_t\circ \Gamma_{\mu_0+1}$ also maps submersively the
polydisc $\Delta_{m(\mu_{0}+1)}(\eta)$ onto an open neighborhood of
the origin in $\C^n$ which contains $\Delta_n(\rho_4)$. Consequently,
in the representation obtained in Lemma~6.21 with $\ell=0$
and $k:=\mu_0+1=2\nu_0+2$ (which is even), namely in the representation
\def\theequation{6.19}\begin{equation}
\bar h(\Gamma_{\mu_0+1}(z_{(\mu_0+1)}))=
\overline{\Pi_{0,\mu_0+1}}(\Gamma_{\mu_0+1}(z_{(\mu_0+1)}),
J^{(\mu_0+1)\ell_0} \bar h(0)), 
\end{equation} 
we can write an arbitrary $t\in\Delta_n(\rho_4)$ in the form
$\Gamma_{\mu_0+1}(z_{(\mu_0+1)})$, and finally, conjugating $(6.19)$,
we obtain a complex algebraic mapping $H$ with the property that
$h(t)=H(t,J^{(\mu_0+1)\ell_0} h(0))$. We may now summarize what we
have proved so far.

\def\thetheorem{6.4}\begin{theorem}
Let $M$ be a real algebraic generic submanifold in $\C^n$ passing
through the origin, of codimension $d\geq 1$ and of CR
dimension $m=n-d\geq 1$.  Assume that $M$ is $\ell_{0}$-nondegenerate
at $0$.  Assume that $M$ is minimal at $0$, let $\nu_{0}$ be the Segre
type of $M$ at $0$ and let $\mu_{0}:=2\nu_{0}+1$ be the Segre type of
$\mathcal{M}$ at $0$. Let $\kappa_0:=(\mu_0+1)\ell_0$. Let
$t=(z,w)\in\C^m\times \C^d$ be holomorphic coordinates vanishing at
$0$ with $T_0M=\{{\rm Im}\, w=0\}$ and let $\rho_1>0$ be such that $M$
is represented by the complex analytic defining equations 
$\xi_j=\Theta_j(\zeta,t)$, $j=1,\dots,d$ in $\Delta_n(\rho_1)$. Then
there exist $\varepsilon>0$, $\rho_4>0$ and there exists a complex
algebraic $\C^n$-valued mapping $H(t,J^{\kappa_0})$ defined for $\vert
t\vert <\rho_4$ and for $\vert J^{\kappa_0}- J_{\rm
Id}^{\kappa_0}\vert< \varepsilon$ which satisfies $H(t,J_{\rm
Id}^{\kappa_0})\equiv t$ and which depends only on the defining
functions $\bar w_j-\Theta_j(\bar z,t)$ of $\mathcal{M}$, such that for
every local holomorphic self-mapping $h$ of $M$ belonging to
$\mathcal{H}_{M,\kappa_0,\varepsilon}^{\rho_2,\rho_1}$, we have
the representation formula
\def\theequation{6.20}\begin{equation}
h(t)=H(t,J^{\kappa_0} h(0)),
\end{equation}
for all $t\in\C^n$ with $\vert t\vert < \rho_4$. Furthermore
the mapping $H$  depends neither on the
choice of smaller radii $\widetilde{\rho}_1\leq \rho_1$,
$\widetilde{\rho}_2\leq \rho_2$, $\widetilde{\rho}_3\leq \rho_3$ and
$\widetilde{\rho}_4\leq \rho_4$ satisfying $0< \widetilde{\rho}_4 <
\widetilde{\rho}_3 < \widetilde{\rho}_2 < \widetilde{\rho}_1$ nor on
the choice of a smaller constant $\widetilde{\varepsilon}< \varepsilon$,
so that the first sentence of property {\bf (3)} in Theorem~4.1 holds
true. Finally, if $M$ is real analytic, the same statement holds with
the word ``algebraic'' everywhere replaced by the word ``analytic''.
\end{theorem}

It remains now to construct the submanifold $E$ whose
existence is stated in Theorem~4.1 and to establish that $\mathcal{H}_{M,
\kappa_0,\varepsilon}^{\rho_2,\rho_1}$ may be endowed with the 
structure of a local real algebraic Lie group.

\subsection*{6.5.~Local real algebraic Lie group structure}
In order to construct this submanifold $E$, we introduce the
$\kappa_0$-th jet mapping $\mathcal{J}^{\kappa_0}:
\mathcal{H}_{M,\kappa_0,\varepsilon}^{\rho_2,\rho_1} \to \C^{N_{n,\kappa_0}}$
defined by $\mathcal{J}^{\kappa_0}(h):=(\partial_t^\alpha h(0))_{\vert
\alpha\vert \leq \kappa_0}=J^{\kappa_0}h(0)$. The following lemma is crucial.

\def\thelemma{6.5}\begin{lemma}
Shrinking $\varepsilon$ if necessary,
the set 
\def\theequation{6.21}\begin{equation}
E:=\mathcal{J}^{\kappa_0}(\mathcal{H}_{M,\kappa_0,\varepsilon}^{
\rho_2,\rho_1})=
\{J^{\kappa_0}h(0): \, 
h\in \mathcal{H}_{M,\kappa_0,\varepsilon}^{\rho_2,\rho_1}\}
\end{equation}
is a real algebraic totally real submanifold of
the polydisc $\{J^{\kappa_0}\in\C^{N_{n,\kappa_0}}: \, 
\vert J^{\kappa_0}- J_{\rm Id}^{\kappa_0} \vert <\varepsilon\}$.
\end{lemma}

\proof
Let $h\in \mathcal{H}_{M,\kappa_0,\varepsilon}^{\rho_2,\rho_1} $.
Substituting the representation formula $h(t)=H(t,J^{\kappa_0} h(0))$ given by
Theorem~6.4 in the defining equations of $M$, we get
\def\theequation{6.22}\begin{equation}
r_j(H(t,J^{\kappa_0} h(0)),
\overline{H}(\tau,J^{\kappa_0}\bar h(0)))=0,
\end{equation}
for $j=1,\ldots,d$ and $(t,\tau)\in\mathcal{M}$
with $\vert t\vert, \, \vert\tau\vert<\rho_4$. As $(t,\tau)\in\mathcal{M}$, 
we replace $\xi$ by $\Theta(\zeta,t)$ and we use the
$2m+d$ coordinates $(t,\zeta)$ on $\mathcal{M}$. So, by expanding 
the functions~\thetag{6.22} in power series with respect 
to $(t,\zeta)$, we can write
\def\theequation{6.23}\begin{equation}
r_j(H(t,J^{\kappa_0}),\overline{H}(\zeta,\Theta(\zeta,t),
\overline{J^{\kappa_0}}))=\sum_{\alpha\in\N^n,\,
\beta\in\N^m}\,
t^\alpha\, \zeta^\beta\, C_{j,\alpha,\beta}(J^{\kappa_0},
\overline{J^{\kappa_0}}). 
\end{equation}
Here,
we obtain an infinite collection of complex-valued real algebraic functions
$C_{j,\alpha,\beta}$ defined in $\{\vert
J^{\kappa_0}-J_{\rm Id}^{\kappa_0}\vert <\varepsilon\}$
with the property that a mapping $H(t,J^{\kappa_0})$
sends $M\cap \Delta_n(\rho_4)$ into $M$
if and only if
\def\theequation{6.24}\begin{equation}
C_{j,\alpha,\beta}(J^{\kappa_0},
\overline{J^{\kappa_0}})=0, \ \ \ \ \
\forall \ j, \, \alpha,\, \beta.
\end{equation}
Consequently, the set $E$ defined by the vanishing of all the
equations~\thetag{6.24} is a real algebraic subset.

It follows from the representation formula~\thetag{6.20}
that the mapping $\mathcal{J}^{\kappa_0}$ is injective and from the Cauchy 
integral formula that $\mathcal{J}^{\kappa_0}$
is continuous on its domain of definition
$\mathcal{H}_{M,\kappa_0,\varepsilon}^{\rho_2,\rho_1}$ endowed with the
topology of uniform convergence on compact sets.

On the reverse side, let $J^{\kappa_0}\in E$. Then the mapping
$h(t):=H(t,J^{\kappa_0})$ defined for $\vert t\vert <\rho_4$
maps $M\cap \Delta_n(\rho_4)$ into $M$. Applying Theorem~6.4 to this
mapping $h(t)$, with $\rho_1$ replaced by $\rho_4$, we deduce that
there exists a radius $\rho_6<\rho_4$ such that we can represent
$h(t)=H(t,J^{\kappa_0}h(0))$ for $\vert t\vert < \rho_6$, with
the same mapping $H$, as stated in the end of Theorem~6.4.
By differentiating this representation with respect to $t$
at $t=0$, we deduce that $J^{\kappa_0}h(0)=([\partial_t^\alpha
H(t,J^{ \kappa_0}h(0))]_{t=0})_{ \vert \alpha \vert \leq
\kappa_0}$. Consequently, since $h(t)=H(t,J^{\kappa_0})$ by
definition, we get $J^{\kappa_0}=([\partial_t^\alpha
H(t,J^{\kappa_0})]_{t=0})_{ \vert \alpha \vert \leq \kappa_0}$.
In conclusion, we proved that
$\mathcal{J}^{\kappa_0}(H(t,J^{\kappa_0}))=J^{\kappa_0}$
for every $J^{\kappa_0}\in E$, so $\mathcal{J}^{\kappa_0}$ has a
continuous local inverse on $E$, formally defined by
$H(t,J^{\kappa_0})$.

It follows from the above two paragraphs that the mapping
$\mathcal{J}^{\kappa_0}$ is a local homeomorphism from a neighborhood
of the identity in $\mathcal{H}_{M,\kappa_0,\varepsilon}^{\rho_2,\rho_1} $
onto its image $E$.

Furthermore, we claim that the real algebraic
subset $E$ is in fact geometrically smooth at every point,
namely it is a real algebraic submanifold. Indeed, let
$J_1^{\kappa_0}$ be a regular point of $E$ where
$E$ is of maximal geometrical dimension $c_0$, with $J_1^{\kappa_0}$
arbitrarily close to the identity jet $J_{\rm Id}^{\kappa_0}$. 
Let $h_1\in
\mathcal{H}_{M,\kappa_0,\varepsilon}^{\rho_2, \rho_1}$ such that 
$J_1^{\kappa_0}=
\mathcal{J}^{\kappa_0}(h_1)$. Let $\mathcal{U}_1$ be a small
neighborhood of $J_1^{\kappa_0}$ in $\C^{N_{n,\kappa_0}}$ in
which $E\cap \mathcal{U}_1$ is a regular $c_0$-dimensional real
algebraic submanifold and consider the complex algebraic mapping
defined over $\mathcal{U}_1$ by
\def\theequation{6.25}\begin{equation}
\mathcal{F}_1(J^{\kappa_0}):= ([\partial_t^\alpha\, (h_1^{-1}(
H(t,J^{\kappa_0})))]_{t=0})_{\vert\alpha\vert\leq \kappa_0}
\in\C^{N_{n, \kappa_0}}.
\end{equation}
We have $\mathcal{F}_1(J_1^{\kappa_0})=J_{\rm
Id}^{\kappa_0}$ and the restriction of $\mathcal{F}_1$ to $E\cap
\mathcal{U}_1$ induces a homeomorphism onto its image, which is a
neighborhood of $J_{\rm Id}^{\kappa_0}$ in $E$.  We remind that
the mapping $J^{\kappa_0}\to ([\partial_t^\alpha\, (
H(t,J^{\kappa_0}))]_{t=0})_{\vert\alpha\vert\leq \kappa_0}$
restricted to $E\cap \mathcal{U}_1$ is the identity and consequently
of constant rank equal to $c_0$. As $h_1$ is invertible, it
follows from the chain rule by developing~\thetag{6.25} that
$\mathcal{F}_1\vert_{E\cap \mathcal{U}_1}$ is also of locally constant
rank equal to $c_0$. This proves that $E$ is a $c_0$-dimensional real
algebraic submanifold in $\C^{N_{n,\kappa_0}}$ through $J_{\rm
Id}^{\kappa_0}$. More generally, this reasoning shows that $E$ is
geometrically smooth at every point.

Finally, applying Lemma~6.3 with the odd integer
$k=\mu_0=2\nu_0+1$ (instead of $k=\mu_0+1$), we get a new, different
representation formula $h(t)=\widetilde{H}(t,J^{\ell_0\mu_0}\bar
h(0))$ (notice $\bar h(0)$). Accordingly, we can define a real
algebraic submanifold $\widetilde{E}$. It is clear that we can
identify $E$ and $\widetilde{E}$, since they both parametrize the
local biholomorphic self-mappings of $M$, so they are algebraically
equivalent by means of the natural projection from the
$\ell_0(\mu_0+1)$-th jet space onto the $\ell_0\mu_0$-th jet
space. Next, we see by differentiating
$h(t)=\widetilde{H}(t,J^{\ell_0\mu_0}\bar h(0))$ with respect to $t$
that
\def\theequation{6.26}\begin{equation} 
J^{\ell_0\mu_0} h(0)=
([\partial_t^\alpha \widetilde{H}(t, 
J^{\ell_0\mu_0}\bar h(0))]_{t=0})_{\vert
\alpha\vert\leq \ell_0\mu_0}.
\end{equation}
Consequently, if $K$ is the holomorphic map defined by
\def\theequation{6.27}\begin{equation}
K(J^{\ell_0\mu_0}):=([\partial_t^\alpha \widetilde{H}
(t,J^{\ell_0\mu_0}]_{t=0})_{
\vert\alpha\vert \leq \ell_0\mu_0}), 
\end{equation}
we get the equality $J^{\ell_0\mu_0}=K(\overline{J^{\ell_0\mu_0}})$
for every $J^{\ell_0\mu_0}\in \widetilde{E}$, which proves that
$\widetilde{E}$ is totally real. It follows that $E$ is
totally real, which completes the proof.
\endproof

\def\thelemma{6.6}\begin{lemma}
The submanifold $E$ is naturally equipped with a local 
real algebraic Lie group structure in a neighborhood of
$J_{\rm Id}^{\kappa_0}$.
\end{lemma}

\proof
Indeed, let us parametrize $E$ by a real algebraic mapping
\def\theequation{6.28}\begin{equation}
\R^{c_0}\ni(e_1,\dots,e_{c_0})\longmapsto 
j_{\kappa_0}(e)\in\C^{N_{n,\kappa_0}},
\end{equation}
where $c_0$ is the dimension of $E$. Here, to avoid excessive formal
complexity, we shall avoid to mention all the polydiscs of
variation of the variables. For $e\in E$, we shall use the notation
\def\theequation{6.29}\begin{equation}
H(t;e):=H(t,j_{\kappa_0}(e)).
\end{equation}
Let $e\in E$ and $e'\in E$, set
$J^{\kappa_0}:=j_{\kappa_0}(e)$ and
$'\!J^{\kappa_0}:=j_{\kappa_0}(e')$.  Then we can define
the Lie group multiplication $\mu_J$ by
\def\theequation{6.30}\begin{equation}
\mu_{J}('\! J^{\kappa_0},
J^{\kappa_0}):=
([\partial_t^\alpha (H(H(t, J^{\kappa_0}),
\, '\!J^{\kappa_0}))]_{
t=0})_{\vert\alpha\vert\leq \kappa_0}.
\end{equation}
Accordingly, in terms of the coordinates $(e_1,\dots,e_{c_0})$ on $E$,
the Lie group multiplication $\mu$ is defined by
\def\theequation{6.31}\begin{equation}
\mu(e,e'):=(j_{\kappa_0})^{-1}(
\mu_J(j_{\kappa_0}(e'), j_{\kappa_0}(e)))\in\R^{c_0}
\end{equation} 
It follows from the algebraicity of the mappings $H$ and
$j_{\kappa_0}$ that the mappings
$\mu_J$ and $\mu$ are algebraic. 

We must check the associativity of $\mu$, namely
$\mu(\mu(e,e'),e'')=\mu(e,\mu(e',e''))$.  So we set
$h(t):=H(t,j_{\kappa_0}(e))$,
$h'(t):=H(t,j_{\kappa_0}(e'))$ and
$h''(t):=H(t,j_{\kappa_0}(e''))$. By the
definition~\thetag{6.30}, we have
$\mu_J(j_{\kappa_0}(e),j_{\kappa_0}(e'))=
J^{\kappa_0} (h\circ h')(0)$. Applying then Theorem~6.4, we get
$H(t,J^{\kappa_0}(h\circ h')(0))\equiv (h\circ
h')(t)$. Consequently, using again~\thetag{6.30} and the associativity
of the composition of mappings, we may compute
\def\theequation{6.32}\begin{equation}
\left\{
\aligned
\mu_J(\mu_J(j_{\kappa_0}(e),j_{\kappa_0}(e')),
j_{\kappa_0}(e''))= & \
\mu_J(J^{\kappa_0}((h\circ h')(0), j_{\kappa_0}(e'')) \\
= & \ 
J^{\kappa_0}((h\circ h')\circ h'')(0) \\
= & \ 
J^{\kappa_0}(h\circ(h'\circ h''))(0) \\
= & \
\mu_J(j_{\kappa_0}(e),J^{\kappa_0}(h'\circ h'')(0)) \\
= & \ 
\mu_J(j_{\kappa_0}(e), 
\mu_J(j_{\kappa_0}(e),j_{\kappa_0}(e'))), 
\endaligned\right.
\end{equation}
which proves the associativity. 

Finally, we may define an algebraic inversion mapping $\iota$ as
follows. First of all, for $J^{\kappa_0}$ close to $J_{\rm
Id}^{\kappa_0}$, the mapping
$h(t):=H(t,J^{\kappa_0})=t+\sum_{\alpha\in\N^n}\,
t^\alpha \, H_\alpha(J^{\kappa_0})$ is an invertible algebraic
biholomorphic mapping. Here, the coefficients
$H_\alpha(J^{\kappa_0})$ are algebraic functions of
$J^{\kappa_0}$ which vanish at $J_{\rm Id}^{\kappa_0}$
(since $H(t,J_{\rm Id}^{\kappa_0})\equiv t$ in
Theorem~6.4). From the algebraic implicit function theorem, it
follows that the local inverse $h^{-1}(t)$ writes uniquely in the form
$h^{-1}(t)= t+\sum_{\alpha\in\N^n}\, t^\alpha\,
\widetilde{H}_\alpha(J^{\kappa_0})=:
\widetilde{H}(t,J^{\kappa_0})$, where the
$\widetilde{H}_\alpha(J^{\kappa_0})$ are algebraic functions of
$J^{\kappa_0}$ also satisfying $\widetilde{H}_\alpha (J_{\rm
Id}^{\kappa_0})=0$. Consequently, choosing $e\in E$ such that
$J^{\kappa_0}=j_{\kappa_0}(e)$, we can 
define
\def\theequation{6.33}\begin{equation}
\iota_J(J^{\kappa_0}):=
([\partial_t^\alpha \widetilde{H}(t,J^{\kappa_0})]_{t=0})_{\vert
\alpha\vert \leq \kappa_0}.
\end{equation}
Accordingly, in terms of the coordinates $(e_1,\dots,e_{c_0})$ on 
$E$, the Lie group inverse mapping is defined by 
\def\theequation{6.34}\begin{equation}
\iota(e):= (j^{\kappa_0})^{-1}(
i_J(j_{\kappa_0}(e))).
\end{equation}
Of course, with this definition we have $\iota_J(J_{\rm
Id}^{\kappa_0})= J_{\rm Id}^{\kappa_0}$. Finally, we leave to
the reader to verify that $\mu_J(j_{\kappa_0}(e),
i_J(j_{\kappa_0}(e)))=J_{\rm Id}^{\kappa_0}$.  This
completes the proof of property {\bf (4)} of Theorem~4.1.
\endproof

\smallskip
\noindent
{\it End of proof of Theorem~4.1.}  We notice that statement {\bf (5)}
does not need to be proved.  Furthermore that the dimensional
inequality $c_0\leq {(n+\kappa_0)!\over n! \ \kappa_0!}$ in {\bf (6)}
follows from the fact each local biholomorphic mapping in the local
Lie group $\mathcal{H}_{M,\kappa_0,\varepsilon}^{\rho_2,\rho_1}\cong
E$ writes uniquely as $h(t)=H(t,J^{\kappa_0} h(0))$, so the complex
dimension of the local Lie group $E$ is $\leq {(n+\kappa_0)!\over n!\
\kappa_0!}$, the dimension of the $\kappa_0$-th jet space.  As $E$ is
totally real, the real dimension of $E$ is also $\leq
{(n+\kappa_0)!\over n!\ \kappa_0!}$. Finally, it follows that the real
local Lie algebra of vector fields
$\mathfrak{Hol}(M,\Delta_n(\rho_5))$ is of dimension $\leq
{(n+\kappa_0)!\over n!\ \kappa_0!}$. The proof of Theorem~4.1 is
complete.
\qed

\section*{\S7.~Description of explicit families of strong tubes in $\C^n$}

\subsection*{7.1.~Introduction}
Theorems~1.1, 1.4 and 1.5 provide sufficient conditions for some real 
analytic real
submanifold in $\mathbb C^n$ to be not locally algebraizable. 
For the sake
of completeness, we exhibit explicit examples of such
nonalgebraizable submanifolds which are effectively strong tubes and
effectively nonalgebraizable, proving corollaries 1.2, 1.3, 1.6 and 1.7.
Consequently we will deal with the two following families of nonalgebraizable 
real
analytic Levi nondegenerate hypersurfaces in $\mathbb C^n$ $(n \geq
2)$ : the Levi nondegenerate strong tube
hypersurfaces in $\mathbb C^n$ and the strongly rigid hypersurfaces in 
$\mathbb C^n$.
For heuristic reasons, we
shall sometimes start with the case $n=2$ and treat the general case $n\geq 2$
afterwards. In fact, our goal will be to construct infinite families
of pairwise non biholomorphically equivalent and non locally
algebraizable hypersurfaces. Our computations for the construction of
families of manifolds with a control on the structure of their
automorphism group are all based on the Lie theory of symmetries of
differential equations. For the convenience of the reader, we recall
briefly the procedure ({\it see}\, [Su2001a,b], [GM2001a,b,c] for more
details).

\subsection*{7.2.~Hypersurfaces and differential equations}
Let $M$ be a real analytic hypersurface in $\mathbb C^n$. Assume
that $M$ is Levi nondegenerate at one of its points 
$p$. Then there exist some local
holomorphic coordinates $(z,w)=(z,u+iv)\in\C^{n-1}\times \C$
vanishing at $p$ such that $M$ is
given by the real analytic equation 
\def\theequation{7.1}\begin{equation}
v=\varphi(z,\bar z,u)=
\varepsilon_1\vert z_1\vert^2+\cdots+
\varepsilon_{n-1}\vert z_{n-1}\vert^2+
\psi(z,\bar z,u),
\end{equation} 
where $\varepsilon_k=\pm 1$, $k=1,\dots,n-1$ and where $\psi={\rm
O}(3)$. Passing to the extrinsic
complexification $\mathcal{M}$ of $M$, we may consider the variables
$\bar z$ and $\bar w$ as independent complex parameters
$\zeta\in\C^{n-1}$ and $\xi\in\C$. Then the associated complex defining 
equation is of the form
\def\theequation{7.2}\begin{equation}
w=\overline{\Theta}(z,\zeta,\xi)=\xi+2i(
\varepsilon_1z_1\zeta_1+\cdots+\varepsilon_{n-1} z_{n-1}\zeta_{n-1}+
\overline{\Xi}(z,\zeta,\xi)),
\end{equation}
where $\overline{\Xi}={\rm O}(3)$.   By [Me1998] ({\it cf.}~\S5.1
above), for $\tau_p=(\zeta_p,\xi_p)$ fixed, the family of complexified
Segre varieties $\mathcal{S}_{\tau_p}:=\{(t,\tau_p):
w=\overline{\Theta}(z,\tau_p)\}$ is invariantly and biholomorphically
attached to $M$.  

Following [Se1931] and [Su2001a,b], we may consider
this family as a family of graphs of the solutions of a second order
completely integrable system of partial differential equations as
follows. By differentiating the left and the right hand sides
of~\thetag{7.2} with respect to $z_k$, we get
\def\theequation{7.3}\begin{equation}
\partial_{z_k}w=\partial_{z_k}\overline{\Theta}(z,\tau)=2i
(\varepsilon_k \zeta_k+
\partial_{z_k}\overline{\Xi}(z,\tau)),
\end{equation}
for $k=1,\dots,n-1$. Here, we consider $w$ as a function of $z$. 
Using the analytic implicit function theorem to
solve $\tau$ in the $1+(n-1)=n$ equations~\thetag{7.2}
and~\thetag{7.3}, we may express $\tau$ in terms of $w$, of $z$ and of
the first order derivative $w_{z_l}$, which yields
\def\theequation{7.4}\begin{equation}
\tau=\Pi(z,w,(\partial_{z_l}w)_{1\leq l\leq n-1}),
\end{equation}
where $\Pi$ is holomorphic in its variables.
If we take the second derivative
$w_{z_{k_1}z_{k_2}}$ of $w$
and replace the value of $\tau$, we get the desired
system of partial differential equations:
\def\theequation{7.5}\begin{equation}
\aligned
\partial^2_{z_{k_1}z_{k_2}}w & \ =
\partial^2_{z_{k_1}z_{k_2}}\overline{\Theta}(z,\tau)
=\partial^2_{z_{k_1}z_{k_2}}\overline{\Theta}(z,
\Pi(z,w,(\partial_{z_l}w)_{1\leq l\leq n-1}))=:\\ 
& \ =: F_{k_1,k_2}(z,w,(\partial_{z_l}w)_{1\leq l\leq n-1}).
\endaligned
\end{equation}
Here, $k_1,k_2=1,\dots,n-1$ and the $F_{k_1,k_2}\equiv F_{k_2,k_1}$
are holomorphic in their variables.  We denote by $\mathcal{E}_M$ this
system of partial differential equations (here, to construct
$\mathcal{E}_M$, we have used the Levi nondegeneracy of $M$ but we
note that if $M$ were finitely nondegenerate the same conclusion would
be true, by considering some derivatives of $w$ of larger order).
Since the solutions of $\mathcal E_M$ are precisely the complexified
Segre varieties $\mathcal{S}_\tau$, the system $\mathcal{E}_M$ is
completely integrable.

To study the local geometry of $M$, we may consider on one hand the
real Lie algebra of infinitesimal CR automorphisms of $M$ ({\it
cf.}~\S2.2), namely $\mathfrak{Aut}_{CR}(M)=2\, {\rm Re}\,
\mathfrak{Hol}(M)$. On the other hand, following the general ideas of
Lie ({\it cf.}~the modern restitution by Olver in [Ol1986, Ch~2]), we
may consider the Lie algebra of infinitesimal generators of the local
symmetry group of the system of partial differential
equations~$\mathcal{E}_M$, which we shall denote by
$\mathfrak{Sym}(\mathcal E_M)$. By definition,
$\mathfrak{Sym}(\mathcal E_M)$ consists of holomorphic vector fields
in the $(z,w)$-space whose local flow transforms the graph of every
solution of $\mathcal{E}_M$ (namely a complexified Segre variety) into
the graph of another solution of $\mathcal{E}_M$ (namely into another
complexified Segre variety). The link between $\mathfrak{Aut}_{CR}(M)$ and
$\mathfrak{Sym}(\mathcal{E}_M)$ is as follows: by [Ca1932, p.~30--32],
one can prove that $\mathfrak{Aut}_{CR}(M)$ is a maximally real subspace
of $\mathfrak{Sym}(\mathcal E_M)$ ({\it see} also [Su2001a,b],
[GM2001a,b,c]).

The computation of explicit generators of $\mathfrak{Sym}(\mathcal{E}_M)$ may
be performed using the Lie theory of symmetries of differential
equations. By inspecting some examples, it appears that dealing with
$\mathfrak{Sym}(\mathcal{E}_M)$ generally shortens the complexity of
the computation of $\mathfrak{Aut}_{CR}(M)$ by at least one half.

The Lie procedure to compute $\mathfrak{Sym}(\mathcal{E}_M)$ is as
follows. Let $J^2_{n-1,1}(\mathbb C)$ denote the space of second order
jets of a function $w(z_1,\dots,z_{n-1})$ of $(n-1)$ complex
variables, equipped with independent coordinates
$(z,w,W_l^1,W_{k_1,k_2}^2)$ corresponding to
$(z,w,w_{z_l},w_{z_{k_1}z_{k_2}})$, where $l=1,\dots,n-1$, where
$k_1,k_2=1,\dots,n-1$, and where we of course identify $W_{k_1,k_2}^2$
with $W_{k_2,k_1}^2$. To the system $\mathcal{E}_M$, we associate the
complex submanifold of $J_{n-1,1}^2(\C)$ defined by replacing the
derivatives of $w$ by the independent jet variables in the system
$\mathcal{E}_M$, which yields ({\it cf.}~\thetag{7.5}):
\def\theequation{7.6}\begin{equation}
W_{k_1,k_2}^2=F_{k_1,k_2}(z,w,(W_l^1)_{1\leq l\leq n-1}), 
\end{equation}
for $k_1,k_2=1,\dots,n-1$.  Let $\Delta_M$ denote this submanifold. By
Lie's theory, every vector field $X=\sum_{k=1}^{n-1}\, Q^k(z,w)\,
\partial_{z_k}+ R(z,w)\, \partial_w$ defined in a neighborhood of the
origin in $\mathbb C^n$ can be uniquely lifted to a vector field
$X^{(2)}$ in $J^2_{n-1,1}(\mathbb C)$, which is called the {\it second
prolongation}\, of $X$ (by definition, the lift $X^{(2)}$ shows how
the flow of $X$ transforms second order jets of graphs of functions
$w(z)$). The coefficients $R_l^1$ and $R_{k_1,k_2}^2$ of the second
prolongation
\def\theequation{7.7}\begin{equation}
X^{(2)}=\sum_{k=1}^{n-1}\, Q^k\, {\partial\over\partial z_k}+
R\, {\partial \over \partial w}
+\sum_{l=1}^{n-1}\, R_l^1\, {\partial\over\partial W_l^1}+
\sum_{k_1,k_2=1}^{n-1}\, R_{k_1,k_2}^2\, 
{\partial\over \partial W_{k_1,k_2}^2},
\end{equation}
are completely determined by the following universal formulas
({\it cf.}~[Ol1986], [Su2001a,b], [GM2001a]):
\def\theequation{7.8}\begin{equation}
\left\{
\aligned
R_{l}^1= & \ \partial_{z_{l}}R+\sum_{m_1} \,
[\delta_{l}^{m_1} \, \partial_wR - \partial_{z_{l}}Q^{m_1}]
\, W_{m_1}^1+\sum_{m_1,m_2} \left[ -\delta_{l}^{m_1} \, 
\partial_wQ^{m_2}\right] \, W_{m_1}^1 \, W_{m_2}^1.\\
R_{k_1,k_2}^2= & \ \partial^2_{z_{k_1}z_{k_2}}R+\sum_{m_1} \, 
\left[\delta_{k_1}^{m_1} \,\partial^2_{z_{k_2}w}R+\delta_{k_2}^{m_1} \, 
\partial^2_{z_{k_1}w}R\,-\partial^2_{z_{k_1}z_{k_2}}Q^{m_1}\right] \, 
W_{m_1}^1+\\
+& \ 
\sum_{m_1,m_2}\,\left[
\delta_{k_1,\,k_2}^{m_1,m_2}\,\partial^2_{w^2}R-\delta_{k_1}^{m_1}\,
\partial^2_{z_{k_2}w}Q^{m_2}-\delta_{k_2}^{m_1}\,\partial^2_{z_{k_1}w}Q^{m_2}
\right]\, 
W_{m_1}^1\,W_{m_2}^1+\\
+&\ \sum_{m_1,m_2,m_3}\,\left[-\delta_{k_1,\,k_2}^{m_1,m_2}\,
\partial^2_{w^2}Q^{m_3}\right]\,W_{m_1}^1\,W_{m_2}^1\,W_{m_3}^1+ \\
+ & \
\sum_{m_1,m_2}\,\left[\delta_{k_1,\,k_2}^{m_1,m_2}\,\partial_wR-\delta_{k_1}^
{m_1}\,
\partial_{z_{k_2}}Q^{m_2}-\delta_{k_2}^{m_1}\,
\partial_{z_{k_1}}Q^{m_2}\right]\,W_{m_1,m_2}^2+\\
+&\ \sum_{m_1,m_2,m_3}\,\left[-\delta_{k_1,\,k_2}^{m_1,m_2}
\, \partial_wQ^{m_3}-\delta_{k_1,\,k_2}^{m_2,m_3}\,\partial_wQ^{m_1}-
\delta_{k_1,\,k_2}^{m_3,m_1}\,\partial_wQ^{m_2}\right]\,W_{m_1}^1\,
W_{m_2,m_3}^2.
\endaligned\right.
\end{equation}
In these formulas, by $\delta_l^m$ we denote the Kronecker symbol
equal to $1$ if $l=m$ and to $0$ otherwise. The multiple Kronecker
symbol $\delta_{l_1,\ l_2}^{m_1,m_2}$ is defined to be the product
$\delta_{l_1}^{m_1}\cdot \delta_{l_2}^{m_2}$.  Finally, in the sums
$\sum_{m_1}$, $\sum_{m_1,m_2}$ and $\sum_{m_1,m_2,m_3}$, the integers
$m_1,m_2,m_3$ run from $1$ to $n-1$. We would like to mention that in
[GM2001a], we also provide some explicit expression of the $k$-th
prolongation $X^{(k)}$ for $k\geq 3$.

Then the {\it Lie criterion}\, states that {\it a holomorphic vector
field $X$ belongs to $\mathfrak{Sym}(\mathcal E_M)$ if and only if its
second prolongation $X^{(2)}$ is tangent to $\Delta_M$} ([Ol1986,
Ch~2]). This gives the following equations:
\def\theequation{7.9}\begin{equation}
R_{k_1,k_2}^2-\sum_{k=1}^{n-1}\, 
Q^k\, \partial_{z_k}F_{k_1,k_2}-
R\, \partial_wF_{k_1,k_2}-
\sum_{l=1}^{n-1}\, 
R_l^1\, \partial_{W_l^1}F_{k_1,k_2}\equiv 0,
\end{equation}
where $1\leq k_1,k_2\leq n-1$ and where each occurence of 
$W_{l_1,l_2}^2$ is replaced by its value $F_{l_1,l_2}$
on $\Delta_M$. By developping~\thetag{7.9}
in power series with respect to 
the variables $W_l^1$, we get an expression of the form
\def\theequation{7.10}\begin{equation}
\sum_{l_1,\dots,l_{n-1}\geq 0}\, 
W_{l_1}^1\cdots W_{l_{n-1}}^1\,
\Phi_{l_1,\dots,l_{n-1}}\equiv 0, 
\end{equation} 
where each term $\Phi_{l_1,\dots,l_{n-1}}$ is a certain linear
partial differential expression involving the derivatives of
$Q^1,\dots,Q^{n-1},R$ up to order two with coefficients being
holomorphic functions of $(z,w)$.  The determination of a system of
generators $X_1,\dots,X_c$ of $\mathfrak{Sym}(\mathcal{E}_M)$ is obtained by
solving the infinite collection of these linear partial differential
equations $\Phi_{l_1,\dots,l_{n-1}}=0$ ({\it cf.}~[Ol1986],
[Su2001a,b], [GM2001a,b,c]).  We shall apply this general procedure to
provide different families of nonalgebraizable real analytic
hypersurfaces in $\mathbb C^n$.

\subsection*{7.3.~Hypersurfaces in $\mathbb C^2$ with control of their CR 
automorphism group} The goal of this paragraph is to construct some
classes of strong tubes, namely tubes having the smallest possible CR
automorphism group.  We start with the case $n=2$ and study afterwards
the case $n\geq 3$ in the next subparagraph.  Let $M_\chi$ be the strong
tube hypersurface in $\C^2$ defined by the equation
\def\theequation{7.11}\begin{equation}
M_\chi: \ \ \ \ \
v=\varphi(y):=y^2+y^6+y^9+y^{10}\,\chi(y).
\end{equation}
where $\chi$ is a real analytic function defined in a
neighborhood of the origin in $\mathbb R$.

\def\thelemma{7.1}\begin{lemma}
The hypersurfaces $M_\chi$ are pairwise not biholomorphically 
equivalent strong tubes.
\end{lemma}

\proof
To check that $M_\chi$ is a strong tube, it suffices to 
show that every hypersurface of the form $v=y^2+y^6+{\rm O}(y^9)$ is
a strong tube (the term $y^9$ will be used afterwards).
Writing $v=(w-\bar w)/2i$ and $y=(z-\bar z)/2i$, 
considering $w$ as a function of $z$ and 
$\bar w$, $\bar z$ as constants, the differentiation of $w$ with
respect to $z$ in~\thetag{7.11} yields:
\def\theequation{7.12}\begin{equation}
\partial_zw=2y+6y^5+{\rm O}(y^8).
\end{equation}
The implicit function theorem yields:
\def\theequation{7.13}\begin{equation}
y=(1/2)\partial_zw
-(3/2^5)(\partial_zw)^5+{\rm O}((\partial_zw)^8).  
\end{equation}
One further differentiation of 
equation~\thetag{7.12} with respect to $z$ gives:
\def\theequation{7.14}\begin{equation}
\partial^2_{zz}w= -i-(15i)\, y^4+
{\rm O}(y^7).
\end{equation}
Replacing $y$ in this equation by its value obtained in~\thetag{7.13},
we obtain the following second order ordinary equation
$\mathcal{E}_M$ satisfied by $\partial_zw$ and $\partial^2{_zz}w$:
\def\theequation{7.15}\begin{equation}
\partial^2_{zz}w = -i-(15i/2^4)(\partial_zw)^4+
{\rm O}((\partial_zw)^7).
\end{equation}
In the four dimensional jet space $J^2_{1,1}(\mathbb C)$ equipped with
the coordinates $(z,w,W^1,W^2)$ the equation of the corresponding
complex hypersurface $\Delta_M$ is of course:
\def\theequation{7.16}\begin{equation}
W^2=-i-(15i/2^4)(W^1)^4+
{\rm O}((W^1)^7).
\end{equation}
Then the Lie criterion states that a holomorphic vector field
$X=Q\,\partial_z + R\, \partial_w$ belongs to $\mathfrak{Sym}(
\mathcal{E}_M)$ if and only if its second prolongation $X^{(2)}=
Q\,\partial_z+R\,\partial_w+R^1\, \partial_{W^1}+ R^2\,
\partial_{W^2}$ is tangent to $\Delta_M$, where the coefficients $R^1$
and $R^2$ are given by the formulas~\thetag{7.8} specified for $n=2$, namely:
\def\theequation{7.17}\begin{equation}
\left\{
\aligned
R^1 =& \ \partial_zR+[\partial_wR-\partial_zQ]\,W^1 -\partial_wQ\, (W^1)^2.\\
R^2 =& \ \partial^2_{zz}R+[2\partial^2_{zw}R-\partial^2_{zz}Q]\, W^1+
[\partial^2_{ww}R-2\partial^2_{zw}Q]\, (W^1)^2-
\partial^2_{ww}Q\, (W^1)^3+\\
& \ \ \ \ \ +[\partial_wR-2\partial_zQ]\, W^2-3\partial_wQ \, W^1W^2.
\endaligned\right.
\end{equation}
The tangency condition yields the 
following equation which is satisfied on $\Delta_M$, {\it i.e.}
after replacing $W^2$ by its value given by~\thetag{7.16}:
\def\theequation{7.18}\begin{equation}
\aligned
R^2+(15i/2^2)R^1(W^1)^3+{\rm O}((W^1)^6)=0.
\endaligned
\end{equation}
By expanding equation~\thetag{7.18} in powers of $W^1$ up to order
five, we obtain the following system of six linear partial
differential equations which must be satisfied by the derivatives of
$Q$ and $R$ up to order two:
\def\theequation{7.19}\begin{equation}
\left\{\begin{array}{lll}
(e_0): & &\partial^2_{zz}R-i(\partial_wR-2\partial_zQ) \equiv 0.\\
\\
(e_1): & & 2\partial^2_{zw}R-\partial^2_{zz}Q \equiv 0.\\
\\
(e_2):& & \partial^2_{ww}R-2\partial^2_{zw}Q \equiv 0.\\
\\
(e_3): & &-\partial^2_{ww}Q+\frac{15i}{2^2}\partial_zR \equiv 0.\\
\\
(e_4): & & -\frac{15i}{2^4}(\partial_wR-2\partial_zQ)+
\frac{15i}{2^2}(\partial_wR-\partial_zQ)\equiv 0.\\
\\
(e_5): & & -\frac{15i}{2^4}(-3\partial_wQ)-
\frac{15i}{2^2}(\partial_wQ)\equiv 0.
\end{array}\right.
\end{equation}
It follows from the equation $(e_5)$ that $\partial_wQ \equiv 0$ which implies
$\partial^2_{ww}Q \equiv 0$. Then by equation $(e_3)$ we obtain $\partial_zR 
\equiv 0$,
implying $\partial^2_{zz}R \equiv 0$. From equation $(e_0)$ we get $\partial_w
R \equiv
2\partial_zQ$ and, from equation $(e_4)$, we get $\partial_wR \equiv
\partial_zQ$. Consequently $\partial_zR \equiv \partial_wR \equiv \partial_zQ 
\equiv \partial_wQ \equiv
0$. Since the two vector fields $\partial_z$ and $\partial_w$
evidently belong to $\mathfrak{Sym}(\mathcal{E}_M)$, it follows that
${\rm dim}_\C \mathfrak{Sym}(\mathcal{E}_M)=2$. Finally, this implies
that $\dim_{\mathbb R}\mathfrak{Aut}_{CR}(M)=2$ and that
$\mathfrak{Aut}_{CR}(M)$ is generated by $\partial_w+\partial_{\bar
w}$ and $\partial_z+\partial_{\bar z}$. 

Next, let $\chi(y)$ and $\chi'(y')$ be two real analytic functions,
and assume that $M_\chi$ and $M_{\chi'}'$ are biholomorphically
equivalent.  Let $t'=h(t)$ be such an equivalence. Reasoning as in \S4 and
taking into account that both are strong tubes, we see that
$h_*(\partial_z)$ and $h_*(\partial_w)$ must be linear combinations of
$\partial_{z'}$ and $\partial_{w'}$ with real coefficients. It follows
that $h$ must be linear, of the form $z'=az+bw$, $w'=cz+dw$, where
$a$, $b$, $c$ and $d$ are real. Since $T_0M_\chi=\{v=0\}$
and $T_0M_{\chi'}'=\{v'=0\}$, we have $c=0$.
Next, in the equation
\def\theequation{7.20}\begin{equation}
\left\{
\aligned
{} & \
d(y^2+y^6+y^9+y^{10}\chi(y))\equiv
[ay+b(y^2+y^6+y^9+y^{10}\chi(y))]^2+\\
& \
+[ay+b(y^2+y^6+y^9+y^{10}\chi(y))]^6
+[ay+b(y^2+y^6+y^9+y^{10}\chi(y))]^9+\\
& \
+[ay+b(y^2+y^6+y^9+y^{10}\chi(y))]^{10}
\chi'(ay+b(y^2+y^6+y^9+y^{10}\chi(y))),
\endaligned\right.
\end{equation}
we firstly see that $b=0$, and then from
\def\theequation{7.21}\begin{equation} 
d(y^2+y^6+y^9+y^{10}\chi(y))\equiv
a^2y^2+a^6y^6+a^9y^9+a^{10}y^{10}\chi'(ay),
\end{equation}
we see that $a=d=1$. In other words, $h={\rm Id}$, whence $y'=y$ and
$\chi'(y')\equiv \chi(y)$. This proves Lemma~7.1.
\endproof

In the remainder of \S7, we shall exhibit other classes of 
hypersurfaces with a control on their CR automorphism
group. Since the computations are generally similar, we shall
summarize them. 

\subsection*{7.4.~Some classes of strong tube hypersurfaces in 
$\C^n$} Generalizing Lemma~7.1, we may state:
\def\thelemma{7.2}\begin{lemma}
The real analytic hypersurfaces $M_{\chi_1,\dots,\chi_{n-1}}
\subset \C^n$ of equation 
\def\theequation{7.22}\begin{equation}
v=\sum_{k=1}^{n-1}\,[ \varepsilon_k\, y_k^2+
y_k^6+y_k^9y_1\cdots y_{k-1}+
y_k^{n+8}\chi_k(y_1,\dots,y_{n-1})], 
\end{equation}
where $\varepsilon_k=\pm 1$,
are pairwise not biholomorphically 
equivalent strong tubes.
\end{lemma}

\proof
The associated system of partial differential equations is of the form
\def\theequation{7.23}\begin{equation}
\left\{
\aligned
\partial^2_{z_kz_k}w= 
& \
-i\varepsilon_k-(15i/2^4)\, (\partial_{z_k}w)^4+
{\rm O}((\partial_{z_1}w)^7)+\cdots+{\rm O}((\partial_{z_{n-1}}w)^7)
,\\
\partial^2_{z_{k_1}z_{k_2}}w=
& \
0, \ \ \ \ \ \ \ \ \ \ \ \ \ \ \ {\rm for} \ \ 
k_1\neq k_2.
\endaligned\right.
\end{equation}
Using the formulas~\thetag{7.8} and inspecting the coefficients
of the monomials in the $W_l^1$ up to order five in the $(n-1)$ equations
extracted from the set of Lie equations
\def\theequation{7.24}\begin{equation}
\left\{
\aligned
{} & 
R_{k,k}^2+(15i/2^2) (W_k^2)\, R_k^1+{\rm O}((W_1^1)^6)+
\cdots+{\rm O}((W_{n-1}^1)^6)=0,\\
&
R_{k_1,k_2}^2=0, \ \ \ \ \ \ \ \ \ \ \ \ \ \ \ {\rm for} \ \ 
k_1\neq k_2,
\endaligned\right.
\end{equation}
we get $\partial_{z_l}R\equiv \partial_wR\equiv \partial_{z_l}Q^k\equiv 
\partial_wQ^k\equiv 0$ for
$l,k=1,\dots,n-1$. Thus, $M_{\chi_1,\dots,\chi_{n-1}}$
is a strong tube.

Next, reasoning as in the end of the proof of Lemma~7.1, we see first
that an equivalence between $M_{\chi_1,\dots,\chi_{n-1}}$ and
$M_{\chi_1',\dots,\chi_{n-1}'}'$ must be of the form
$z_k'=\sum_{l=1}^{n-1}\, \lambda_k^l\, z_l$, $w'=\mu\, w$, where
$\lambda_k^l$, $1\leq l,k\leq n-1$ and $\mu$ are real. Inspecting the
terms of degree $9, 10,\dots,n+7$, we get $\lambda_k^l=0$ if $k\neq
l$, {\it i.e.} $y_k'=\lambda_k^k\,y_k$ and $w'=\mu\, w$. Finally,
$\lambda_k^k=1$ and $\mu=1$, which completes the proof.
\endproof

\subsection*{7.5.~Families of strongly rigid hypersurfaces}
Alongside the same recipe, we can study some classes of hypersurfaces
of the form $v=\varphi(z\bar z)$.

\def\thelemma{7.3}\begin{lemma}
The Lie algebra $\mathfrak{Hol}(M_\chi)$ of the rigid real analytic
hypersurfaces $M_\chi$ in $\C^2$ of equation $v=\varphi(z\bar
z)=z\bar z+z^5\bar z^5+z^7\bar z^7+ z^8\bar z^8\chi(z\bar z)$ is
two-dimensional and generated by $\partial_w$ and
$iz\partial_z$. Furthermore, $M_\chi$ is biholomorphically equivalent
to $M_{\chi'}'$ if and only if $\chi=\chi'$.
\end{lemma}

\proof
The associated differential equation is of the form
\def\theequation{7.25}\begin{equation}
\partial^2_{zz}w=
[5z^3/4](\partial_zw)^5-
[21z^5/32](\partial_zw)^7+{\rm O}((\partial_zw)^9).
\end{equation}
Extracting from the associated Lie equations~\thetag{7.10} the
coefficients of the monomials $(W^1)^4$, $(W^1)^5$, $(W^1)^6$ and
$(W^1)^7$, we obtain four equations which are solved by $z\partial_zQ-Q\equiv
0$, $\partial_wQ\equiv 0$, $\partial_zR\equiv 0$ and $\partial_wR\equiv 0$. 
 Next, if $M_\chi$
and $M_{\chi'}'$ are biholomorphically equivalent, reasoning as in
\S4, taking into account that $h_*(iz\partial_z)$ and
$h_*(\partial_w)$ are linear combinations of $iz'\partial_{z'}$ and
$\partial_{w'}$ with real coefficients, we see first that $z'=\lambda\,
z\, e^{\gamma w/2i}$ and $w'=\mu\,w$ for some three real constants
$\gamma$, $\lambda\neq 0$ and $\mu\neq 0$.  Replacing $z'$ and $w'$ in
the equation of $M_{\chi'}'$, we get $\gamma=0$, $\mu=1$ and $\lambda\pm
1$. In other words, $z'=\pm z$ and $w'=w$, which entails $\chi'(z'\bar
z')\equiv \chi(z\bar z)$, as claimed.
\endproof

Perturbing this family we may exhibit other strongly rigid hypersurfaces~:

\def\thelemma{7.4}\begin{lemma}
The Lie algebra $\mathfrak{Hol}(M_\chi)$ of the real analytic
hypersurfaces $M_\chi$ in $\C^2$ of equation $v=\varphi(z,\bar
z)=z\bar z+z^5\bar z^5+z^7\bar z^7+ z^8\bar z^8(z+\bar z)+ z^{10}\bar
z^{10}\chi(z,\bar z)$ is one-dimensional and generated by $\partial_w$.
Furthermore, $M_\chi$ is biholomorphically
equivalent to $M_{\chi'}'$ if and only if $\chi=\chi'$.
\end{lemma}

\proof
We already know that $z\partial_zQ-Q\equiv \partial_wQ\equiv \partial_zR 
\equiv \partial_wR\equiv 0$.
Extracting from the associated Lie equations~\thetag{7.10} the
coefficient of the monomials $(W^1)^8$, we also get $Q\equiv 0$. Next,
let $M_\chi$ and $M_{\chi'}'$ be biholomorphically equivalent.  Let
$t'=h(t)$ be such an equivalence. Using $h_*(\partial_w)= \mu\,
\partial_{w'}$, where $\mu\in\R$ is nonzero, we get $z'=f(z)$ and
$w'=\mu w+g(z)$. Next from the equation
\def\theequation{7.26}\begin{equation}
\left\{
\aligned
{}
&
\mu(z\bar z+z^5\bar z^5+z^7\bar z^7+
z^8\bar z^8(z+\bar z)+{\rm O}(z^9\bar z^9))+
[g(z)-\bar g(\bar z)]/2i\equiv \\
& \
\equiv f(z)\bar f(\bar z)+f(z)^5\bar f(\bar z)^5+
f(z)^7\bar f(\bar z)^7+
f(z)^8\bar f(\bar z)^8(f(z)+\bar f(\bar z))+
{\rm O}(z^9\bar z^9),
\endaligned\right.
\end{equation}
we get firstly $f(z)=\sqrt{\vert \mu\vert}e^{i\theta}z$ by differentiating
with respect to $\bar z$ at $\bar z=0$ and secondly 
$\mu=e^{i\theta}=1$, which completes the proof.
\endproof

We provide a second family of strongly rigid hypersurfaces in $\mathbb
C^2$ with a one-dimensional Lie algebra\,:

\def\thelemma{7.5}\begin{lemma}
The Lie algebra $\mathfrak{Hol}(M_\chi)$ of the real analytic
hypersurfaces $M_\chi\subset \mathbb C^2$ of equation
$v=z\bar{z}+z^5\bar{z}^5(z+\bar{z})+z^{10}\bar{z}^{10}\chi(z,\bar{z})$ is
one-dimensional and generated by $\partial_w$. Furthermore $M_\chi$ is
biholomorphically equivalent to $M'_{\chi'}$ if and only if
$\chi=\chi'$.
\end{lemma}

\proof 
The derivatives $\partial_zw$ and $\partial^2_{z^2}w$ of $w$ with
respect to $z$ are given by\,:
\def\theequation{7.27}\begin{equation}
\left\{
\aligned
\partial_zw = & \
2i
\bar{z}+12iz^5\bar{z}^5+10iz^4\bar{z}^6+
{\rm O}(\bar{z}^{10})),\\
\partial^2_{zz}w= 
& \ 
60i z^4\bar{z}^5+40iz^3\bar{z}^6+
{\rm O}(\bar{z}^{10}).
\endaligned\right.
\end{equation}
Replacing $\bar{z}$ in the second equation by its expression given by
the first equation we obtain the following second order differential
equation, interpreted in the jet space\,:
\def\theequation{7.28}\begin{equation}
W_2=[15z^4/8]\,(W^1)^5-
[5iz^3/8]\,(W^1)^6
-[225z^9/64]\,(W^1)^9+
{\rm O}((W^1)^{10}).
\end{equation}
Solving the partial differential equations involving $Q$, $\partial_zQ$,
$\partial_wQ$, $\partial_zR$ and $\partial_wR$ given in the coefficients of 
$(W^1)^4$,
$(W^1)^5$, $(W^1)^6$, $(W^1)^7$ and $(W^1)^9$ we obtain $Q\equiv
\partial_zQ\equiv \partial_wQ\equiv \partial_zR\equiv \partial_wR\equiv 0$ 
which is the desired
information. Finally, proceeding exactly as in the end of 
the proof of Lemma~7.4, we see that the $M_\chi$ are
pairwise biholomorphically not equivalent.
\endproof

 The dimension of $\mathfrak{Hol}(M)$ for
the five examples of Corollary~1.3, for the seven examples of
Theorem~1.4, for the seven examples of Corollary~1.7 and for 
the hypersurface $v=e^{z\bar z}-1$ at a point $p$ with $z_p\neq 0$ was
 computed with the package {\tt
diffalg} of Maple Release~6. Since at a point $p$ with $z_p\neq
0$ the hypersurface $v=e^{z\bar z}-1$ is biholomorphically equivalent
to the hypersurface $M_a$ of equation $v=\varphi^a(y):=e^{a(e^y-1)}-1$
with $a=\vert z_p\vert^2$, this defines a strong tube.
Applying Theorem~1.1 and Lemma~3.3, we see that $M_a$ is not locally
algebraizable at the origin, because $\varphi^a_{yy}(y)=
ae^ye^{a(e^y-1)} + a^2e^{2y}e^{a(e^y-1)}$ and $\varphi^a_y(y)=
ae^ye^{a(e^y-1)}$ are algebraically independent.  Finally all the examples 
of Corollary~1.3, Theorem~1.4 and
Corollary~1.6 are not locally algebraic since they satisfy the required 
transcendence conditions.

\section*{\S8.~Analyticity versus algebraicity}
Intuitively there seems to be  much
more analytic mappings, manifolds and varieties than algebraic
ones. Our goal is to elaborate a precise statement about this.
 By complexification, every local real analytic object yields a
local complex analytic object, so we shall only work in the
holomorphic category. Let $\Delta_n$ be the complex polydisc of radius
one in $\C^n$ and $\overline{\Delta}_n$ its closure. Let $k\in\N$. We
consider the space $\mathcal{O}^k (\overline{\Delta}_n):=
\mathcal{O}(\Delta_n)\cap \mathcal{C}^k (\overline{\Delta}_n)$ of
holomorphic functions extending up to
the boundary as a function of class $\mathcal{C}^k$ . This is a Banach 
space for the $\mathcal{C}^k$ norm $\vert \vert \varphi \vert
\vert_k:=\sum_{l=0}^k \sup_{z\in \overline{\Delta}_n}\, \vert
\varphi_{z^l}(z)\vert$. 
The last statements of Corollaries~1.2 and~1.6 are a direct
consequence of the following lemma.

\def\thelemma{8.1}\begin{lemma}
The set of holomorphic functions
$\varphi\in\mathcal{O}^k(\overline{\Delta}_n)$ such that there exists
a polynomial $P$ such that
\def\theequation{8.1}\begin{equation}
P(z,j^k\varphi(z))\equiv 0,
\end{equation}
is of first category, namely it can be represented as the countable
union of nowhere dense closed subsets. Conversely, the set of functions
$\varphi\in\mathcal{O}^k(\overline{\Delta}_n)$ such that there is no
algebraic dependence relation like~\thetag{8.1} is generic in the
sense of Baire, namely it can be represented as the countable
intersection of everywhere dense open subsets.
\end{lemma}

\proof
Let $N\in\N$. Consider the set $F_N$ of functions $\varphi$ such that
there exists a polynomial of degree $N$  satisfying~\thetag{8.1}.  It
suffices to show that $F_N$ is closed and that its complement is
everywhere dense.  Suppose that a sequence $(\varphi^{(m)})_{m\in\N}$
converges to $\varphi\in\mathcal{O}^k(\overline{\Delta}_n)$. Let the
zero-set of a degree $N$ polynomial $P_N^{(m)}(z,J_k)$ contain the
graph of the $k$-jet of $\varphi^{(m)}$. The coefficients of
$P_N^{(m)}$ belong to a certain complex projective space $P_A(\C)$,
where the integer $A=A(n,k)$ is independent of $m$. By compactness of
$P_A(\C)$, passing to a subsequence if necessary, the $P_N^{(m)}$
converge to a nonzero polynomial $P_N$. By continuity,
$P_N(z,j^k\varphi(z))= 0$ for all $z\in
\mathcal{O}(\overline{\Delta}_n)$. We claim that the complement of the
union of the $F_N$ is dense in $\mathcal{O}^k(\overline{\Delta}_n)$.
Indeed, let $\varphi(z)$ be such that there exists a degree $N$
polynomial $P$ satisfying~\thetag{8.1}. Fix $z_0\in\Delta_n$ having
rational real and imaginary parts.  Then the complex numbers $z_0$,
$\partial_z^\alpha\varphi(z_0)$, $\vert\alpha\vert \leq k$, are
algebraically dependent. By a Cantorian argument, there exists complex
numbers $\chi_0^\alpha$ arbitrarily close to
$\partial_z^\alpha\varphi(z_0)$ such that $z_0$, $\chi_0^\alpha$ are
algebraically independent. Let $\chi(z)$ be a polynomial with
$\partial_z^\alpha(z_0)=\chi_0^\alpha-\partial_t^\alpha
\varphi(z_0)$. We can choose $\chi$ to be arbitrarily close
to zero in the $\mathcal{C}^k(\overline{\Delta}_n)$ norm. Then the function
$\varphi(z)+\chi(z)$ is not Nash algebraic.
\endproof

\vfill
\end{document}

%% file: straightening.pstex_t
\begin{picture}(0,0)%
\epsfig{file=straightening.pstex}%
\end{picture}%
\setlength{\unitlength}{3947sp}%
\begingroup\makeatletter\ifx\SetFigFont\undefined%
\gdef\SetFigFont#1#2#3#4#5{%
  \reset@font\fontsize{#1}{#2pt}%
  \fontfamily{#3}\fontseries{#4}\fontshape{#5}%
  \selectfont}%
\fi\endgroup%
\begin{picture}(5424,1824)(1189,-2173)
\put(1351,-2044){\makebox(0,0)[lb]{\smash{\SetFigFont{9}{10.8}{\familydefault}{\mddefault}{\updefault}{\sc Figure~1: Local algebraic straightening of the orbits of $G_1(t;\epsilon_1)$}}}}
\put(5681,-1517){\makebox(0,0)[lb]{\smash{\SetFigFont{7}{8.4}{\familydefault}{\mddefault}{\updefault}$t_1$}}}
\put(2813,-1625){\makebox(0,0)[lb]{\smash{\SetFigFont{7}{8.4}{\familydefault}{\mddefault}{\updefault}$0$}}}
\put(2817,-516){\makebox(0,0)[lb]{\smash{\SetFigFont{7}{8.4}{\familydefault}{\mddefault}{\updefault}$t_2,\dots,t_n$}}}
\put(5555,-763){\makebox(0,0)[lb]{\smash{\SetFigFont{7}{8.4}{\familydefault}{\mddefault}{\updefault}$\mathcal{C}_{t_2',\dots,t_n'}'$}}}
\put(4938,-981){\makebox(0,0)[lb]{\smash{\SetFigFont{7}{8.4}{\familydefault}{\mddefault}{\updefault}$G_1(G_1(0,t_2',\dots,t_n';t_1');\epsilon_1)$}}}
\put(3910,-1137){\makebox(0,0)[lb]{\smash{\SetFigFont{7}{8.4}{\familydefault}{\mddefault}{\updefault}$G_1(0,t_2',\dots,t_n';t_1')$}}}
\put(2789,-1211){\makebox(0,0)[lb]{\smash{\SetFigFont{7}{8.4}{\familydefault}{\mddefault}{\updefault}$G_1(0,t_2',\dots,t_n';0)$}}}
\put(3808,-701){\makebox(0,0)[lb]{\smash{\SetFigFont{7}{8.4}{\familydefault}{\mddefault}{\updefault}$\sim$}}}
\put(4024,-679){\makebox(0,0)[lb]{\smash{\SetFigFont{7}{8.4}{\familydefault}{\mddefault}{\updefault}$t$}}}
\put(3599,-736){\makebox(0,0)[lb]{\smash{\SetFigFont{7}{8.4}{\familydefault}{\mddefault}{\updefault}$t^*$}}}
\end{picture}

%% file: redressement.pstex_t
\begin{picture}(0,0)%
\epsfig{file=redressement.pstex}%
\end{picture}%
\setlength{\unitlength}{3947sp}%
\begingroup\makeatletter\ifx\SetFigFont\undefined%
\gdef\SetFigFont#1#2#3#4#5{%
  \reset@font\fontsize{#1}{#2pt}%
  \fontfamily{#3}\fontseries{#4}\fontshape{#5}%
  \selectfont}%
\fi\endgroup%
\begin{picture}(5116,3594)(64,-2795)
\put(3327,-24){\makebox(0,0)[lb]{\smash{\SetFigFont{8}{9.6}{\familydefault}{\mddefault}{\updefault}$M'$}}}
\put(4134,563){\makebox(0,0)[lb]{\smash{\SetFigFont{8}{9.6}{\familydefault}{\mddefault}{\updefault}$v'$}}}
\put(3348,565){\makebox(0,0)[lb]{\smash{\SetFigFont{8}{9.6}{\familydefault}{\mddefault}{\updefault}$\C_{t'}^n$}}}
\put(370,-36){\makebox(0,0)[lb]{\smash{\SetFigFont{8}{9.6}{\familydefault}{\mddefault}{\updefault}$M$}}}
\put(383,550){\makebox(0,0)[lb]{\smash{\SetFigFont{8}{9.6}{\familydefault}{\mddefault}{\updefault}$\C_t^n$}}}
\put(1199,544){\makebox(0,0)[lb]{\smash{\SetFigFont{8}{9.6}{\familydefault}{\mddefault}{\updefault}$v$}}}
\put(1570,572){\makebox(0,0)[lb]{\smash{\SetFigFont{8}{9.6}{\familydefault}{\mddefault}{\updefault}$x,u$}}}
\put(4159,-1272){\makebox(0,0)[lb]{\smash{\SetFigFont{8}{9.6}{\familydefault}{\mddefault}{\updefault}$v''$}}}
\put(3329,-1838){\makebox(0,0)[lb]{\smash{\SetFigFont{8}{9.6}{\familydefault}{\mddefault}{\updefault}$M''$}}}
\put(3350,-1244){\makebox(0,0)[lb]{\smash{\SetFigFont{8}{9.6}{\familydefault}{\mddefault}{\updefault}$\C_{t''}^n$}}}
\put(4458,-1249){\makebox(0,0)[lb]{\smash{\SetFigFont{8}{9.6}{\familydefault}{\mddefault}{\updefault}$x'',u'''$}}}
\put(4726,-1857){\makebox(0,0)[lb]{\smash{\SetFigFont{8}{9.6}{\familydefault}{\mddefault}{\updefault}$y''$}}}
\put(4711,-33){\makebox(0,0)[lb]{\smash{\SetFigFont{8}{9.6}{\familydefault}{\mddefault}{\updefault}$y'$}}}
\put(1782,-63){\makebox(0,0)[lb]{\smash{\SetFigFont{8}{9.6}{\familydefault}{\mddefault}{\updefault}$y$}}}
\put(2442,-98){\makebox(0,0)[lb]{\smash{\SetFigFont{8}{9.6}{\familydefault}{\mddefault}{\updefault}$\Phi$}}}
\put(3616,-578){\makebox(0,0)[lb]{\smash{\SetFigFont{8}{9.6}{\familydefault}{\mddefault}{\updefault}$M'$ algebraic}}}
\put(4499,367){\makebox(0,0)[lb]{\smash{\SetFigFont{8}{9.6}{\familydefault}{\mddefault}{\updefault}$x',u'$}}}
\put(3325,-2401){\makebox(0,0)[lb]{\smash{\SetFigFont{8}{9.6}{\familydefault}{\mddefault}{\updefault}$M''$ {\sf algebraic} pseudotube}}}
\put(753,-575){\makebox(0,0)[lb]{\smash{\SetFigFont{8}{9.6}{\familydefault}{\mddefault}{\updefault}$M$ strong tube}}}
\put(1115,-1556){\makebox(0,0)[lb]{\smash{\SetFigFont{8}{9.6}{\familydefault}{\mddefault}{\updefault}$\Phi'':=\Psi'\circ\Phi$}}}
\put(1178,-2695){\makebox(0,0)[lb]{\smash{\SetFigFont{8}{9.6}{\familydefault}{\mddefault}{\updefault}{\sc Figure~2: Algebraic straightening $\Psi'$ of $M'$}}}}
\put(4156,-911){\makebox(0,0)[lb]{\smash{\SetFigFont{8}{9.6}{\familydefault}{\mddefault}{\updefault}$\Psi'$}}}
\end{picture}

%% file: polydiscs.pstex_t
\begin{picture}(0,0)%
\epsfig{file=polydiscs.pstex}%
\end{picture}%
\setlength{\unitlength}{3947sp}%
\begingroup\makeatletter\ifx\SetFigFont\undefined%
\gdef\SetFigFont#1#2#3#4#5{%
  \reset@font\fontsize{#1}{#2pt}%
  \fontfamily{#3}\fontseries{#4}\fontshape{#5}%
  \selectfont}%
\fi\endgroup%
\begin{picture}(6124,2829)(551,-2443)
\put(1623,-1071){\makebox(0,0)[lb]{\smash{\SetFigFont{9}{10.8}{\familydefault}{\mddefault}{\updefault}$M$}}}
\put(4623,-1071){\makebox(0,0)[lb]{\smash{\SetFigFont{9}{10.8}{\familydefault}{\mddefault}{\updefault}$M$}}}
\put(5499,-290){\makebox(0,0)[lb]{\smash{\SetFigFont{9}{10.8}{\familydefault}{\mddefault}{\updefault}nondegeneracy}}}
\put(5502,-171){\makebox(0,0)[lb]{\smash{\SetFigFont{9}{10.8}{\familydefault}{\mddefault}{\updefault}finite}}}
\put(5495,-410){\makebox(0,0)[lb]{\smash{\SetFigFont{9}{10.8}{\familydefault}{\mddefault}{\updefault}assumption}}}
\put(4418,-176){\makebox(0,0)[lb]{\smash{\SetFigFont{7}{8.4}{\familydefault}{\mddefault}{\updefault}$\Delta_n(\rho_3)$}}}
\put(4702,-19){\makebox(0,0)[lb]{\smash{\SetFigFont{7}{8.4}{\familydefault}{\mddefault}{\updefault}$\Delta_n(\rho_2)$}}}
\put(3411,-768){\makebox(0,0)[lb]{\smash{\SetFigFont{8}{9.6}{\familydefault}{\mddefault}{\updefault}$\Delta_n(\rho_4)$}}}
\put(1825,-2349){\makebox(0,0)[lb]{\smash{\SetFigFont{9}{10.8}{\familydefault}{\mddefault}{\updefault}{\sc Figure~3: Nest of polydiscs centered at $0\in M$}}}}
\put(5523,-1728){\makebox(0,0)[lb]{\smash{\SetFigFont{9}{10.8}{\familydefault}{\mddefault}{\updefault}Minimality}}}
\put(5517,-1848){\makebox(0,0)[lb]{\smash{\SetFigFont{9}{10.8}{\familydefault}{\mddefault}{\updefault}assumption}}}
\put(5511,-2006){\makebox(0,0)[lb]{\smash{\SetFigFont{9}{10.8}{\familydefault}{\mddefault}{\updefault}$\rho_4\sim (\rho_1)^ N$}}}
\put(5506,-2169){\makebox(0,0)[lb]{\smash{\SetFigFont{9}{10.8}{\familydefault}{\mddefault}{\updefault}$N>>1$}}}
\put(3492,-1490){\makebox(0,0)[lb]{\smash{\SetFigFont{5}{6.0}{\familydefault}{\mddefault}{\updefault}$\Delta_n(\rho_5)$}}}
\put(5563,-1071){\makebox(0,0)[lb]{\smash{\SetFigFont{9}{10.8}{\familydefault}{\mddefault}{\updefault}${\rm Re}\, w, \, z, \bar z$}}}
\put(4890,126){\makebox(0,0)[lb]{\smash{\SetFigFont{6}{7.2}{\familydefault}{\mddefault}{\updefault}$\Delta_n(\rho_1)$}}}
\put(3292,254){\makebox(0,0)[lb]{\smash{\SetFigFont{9}{10.8}{\familydefault}{\mddefault}{\updefault}${\rm Im}\, w$}}}
\end{picture}

%% file: complexification.pstex_t
\begin{picture}(0,0)%
\epsfig{file=complexification.pstex}%
\end{picture}%
\setlength{\unitlength}{3947sp}%
\begingroup\makeatletter\ifx\SetFigFont\undefined%
\gdef\SetFigFont#1#2#3#4#5{%
  \reset@font\fontsize{#1}{#2pt}%
  \fontfamily{#3}\fontseries{#4}\fontshape{#5}%
  \selectfont}%
\fi\endgroup%
\begin{picture}(5724,2349)(56,-1625)
\put(2887,146){\makebox(0,0)[lb]{\smash{\SetFigFont{7}{8.4}{\familydefault}{\mddefault}{\updefault}The complexification of a real analytic}}}
\put(2887,-389){\makebox(0,0)[lb]{\smash{\SetFigFont{7}{8.4}{\familydefault}{\mddefault}{\updefault}whose leaves coincide with the complexified}}}
\put(2887,-657){\makebox(0,0)[lb]{\smash{\SetFigFont{7}{8.4}{\familydefault}{\mddefault}{\updefault}These leaves also coincide with the intersection}}}
\put(2887, 12){\makebox(0,0)[lb]{\smash{\SetFigFont{7}{8.4}{\familydefault}{\mddefault}{\updefault}CR-generic manifold $M$ carries two }}}
\put(2887,-123){\makebox(0,0)[lb]{\smash{\SetFigFont{7}{8.4}{\familydefault}{\mddefault}{\updefault}complex foliations $\mathcal{F}_{\mathcal{L}}$ and $\mathcal{F}_{\underline{\mathcal{L}}}$ directed}}}
\put(2887,-256){\makebox(0,0)[lb]{\smash{\SetFigFont{7}{8.4}{\familydefault}{\mddefault}{\updefault}by the complefixied CR-vector fields $\mathcal{L}$ and $\underline{\mathcal{L}}$}}}
\put(2887,-524){\makebox(0,0)[lb]{\smash{\SetFigFont{7}{8.4}{\familydefault}{\mddefault}{\updefault}Segre varieties $\mathcal{S}_{\tau_p}$ and $\underline{\mathcal{S}}_{t_p}$.}}}
\put(2887,-792){\makebox(0,0)[lb]{\smash{\SetFigFont{7}{8.4}{\familydefault}{\mddefault}{\updefault}of $\mathcal{M}$ with the horizontal slices $\{\tau=\tau_p\}$}}}
\put(2887,-924){\makebox(0,0)[lb]{\smash{\SetFigFont{7}{8.4}{\familydefault}{\mddefault}{\updefault}and with the vertical slices $\{t=t_p\}$.}}}
\put(866,-1099){\makebox(0,0)[lb]{\smash{\SetFigFont{8}{9.6}{\familydefault}{\mddefault}{\updefault}$0$}}}
\put(1480,417){\makebox(0,0)[lb]{\smash{\SetFigFont{8}{9.6}{\familydefault}{\mddefault}{\updefault}$\{t=t_p\}$}}}
\put(1688,-635){\makebox(0,0)[lb]{\smash{\SetFigFont{8}{9.6}{\familydefault}{\mddefault}{\updefault}$\mathcal{F}$}}}
\put(1617,224){\makebox(0,0)[lb]{\smash{\SetFigFont{8}{9.6}{\familydefault}{\mddefault}{\updefault}$\underline{\mathcal{S}}_{t_p}$}}}
\put(2375,503){\makebox(0,0)[lb]{\smash{\SetFigFont{8}{9.6}{\familydefault}{\mddefault}{\updefault}$\underline{\Lambda}$}}}
\put(2183,-381){\makebox(0,0)[lb]{\smash{\SetFigFont{8}{9.6}{\familydefault}{\mddefault}{\updefault}$\{\tau=\tau_p\}$}}}
\put(1625,-1081){\makebox(0,0)[lb]{\smash{\SetFigFont{8}{9.6}{\familydefault}{\mddefault}{\updefault}$t_p$}}}
\put(2315,-1060){\makebox(0,0)[lb]{\smash{\SetFigFont{8}{9.6}{\familydefault}{\mddefault}{\updefault}$t$}}}
\put(1859,-839){\makebox(0,0)[lb]{\smash{\SetFigFont{8}{9.6}{\familydefault}{\mddefault}{\updefault}$\mathcal{L}$}}}
\put(631,-304){\makebox(0,0)[lb]{\smash{\SetFigFont{8}{9.6}{\familydefault}{\mddefault}{\updefault}$\tau_p$}}}
\put(528,220){\makebox(0,0)[lb]{\smash{\SetFigFont{8}{9.6}{\familydefault}{\mddefault}{\updefault}$\mathcal{M}$}}}
\put(1229,-95){\makebox(0,0)[lb]{\smash{\SetFigFont{8}{9.6}{\familydefault}{\mddefault}{\updefault}$\underline{\mathcal{F}}$}}}
\put(869,123){\makebox(0,0)[lb]{\smash{\SetFigFont{8}{9.6}{\familydefault}{\mddefault}{\updefault}$\underline{\mathcal{L}}$}}}
\put(1738,-267){\makebox(0,0)[lb]{\smash{\SetFigFont{8}{9.6}{\familydefault}{\mddefault}{\updefault}$\mathcal{S}_{\tau_p}$}}}
\put(1074,-1496){\makebox(0,0)[lb]{\smash{\SetFigFont{9}{10.8}{\familydefault}{\mddefault}{\updefault}{\sc Figure~4: Geometry of the complexification $\mathcal{M}$}}}}
\end{picture}

%% file: nonalg.bbl
\begin{thebibliography}{XXXXXX}

\bibitem[Ar1974]{ar1974}
{\sc Arnold}, V.I.:
{\em \'Equations diff\'erentielles ordinaires. Champs de vecteurs, 
groupes \`a un param\`etre, diff\'eomorphismes, 
flots, syst\`emes lin\'eaires, stabilit\'e des positions 
d'\'equilibre, th\'eorie des oscillations, \'equations diff\'erentielles
sur les vari\'et\'es}. Traduit du russe par Djilali Embarek, 
\'Editions Mir, Moscou, 1974. 267pp.

\bibitem[BER1999]{ber1999} 
{\sc Baouendi}, M.S.;
{\sc Ebenfelt}, P.; {\sc Rothschild}, L.P.: 
{\em Rational dependence of smooth and analytic CR mappings on
their jets}. Math. Ann. {\bf 315} (1999), 205--249.

\bibitem[BER2000]{ber2000}
{\sc Baouendi}, M.S.; {\sc Ebenfelt}, P.;
{\sc Rothschild}, L.P.:
{\em Local geometric properties of real submanifolds in complex space}, 
Bull. Amer. Math. Soc. {\bf 37} (2000), no.3, 309--336.

\bibitem[Be1996]{be1996}
{\sc Bella\"{\i}che},~A.:
{\em SubRiemannian Geometry}, Progress in Mathematics
{\bf 144}, Birkh\"auser Verlag, Basel/Switzerland,
1996, 1--78.

\bibitem[Bs1991]{bs1991}
{\sc Beloshapka}, V.K.:
{\em On holomorphic transformations of a quadric}, 
Mat. Sb. {\bf 182} (1991), no.2, 203--219;
English transl. in Math. USSR Sb. {\bf 72} (1992), no.1, 189--205.

\bibitem[Ca1932]{ca1932}
{\sc Cartan}, \'E.:
{\em Sur la g\'eom\'etrie pseudo-conforme des hypersurfaces
de l'espace de deux variables complexes, I},
Ann. Math. Pura Appl. {\bf 11} (1932), 17--90.

\bibitem[CM1974]{cm1974}
{\sc Chern}, S.S.; {\sc Moser}, J.K.:
{\em Real hypersurfaces in complex manifolds},
Acta Math. {\bf 133} (1974), no.2, 219--271.

\bibitem[CPS2000]{cps2000}
{\sc Coupet}, B.;
{\sc Pinchuk}, S.; 
{\sc Sukhov}, A.:
{\em On partial analyticity of CR mappings},
Math. Z. {\bf 235} (2000), 541--557.

\bibitem[DP2003]{dp2003}
{\sc Diederich}, D.;
{\sc Pinchuk}, S.:
{\em Regularity of continuous CR-maps in arbitrary dimension}, to appear in Michigan Math. J. (2003).

\bibitem[Eb1996]{e1996}
{\sc Ebenfelt}, P.:
{\em On the unique continuation problem for CR mappings 
into nonminimal hypersurfaces}, J. Geom. Anal. {\bf 6} (1996), no.3,
385--405.

\bibitem[GM2001a]{gm2001a}
{\sc Gaussier}, H.; {\sc Merker}, J.:
{\em Estimates on the dimension of the symmetry group of a system of
$k$-order partial differential equations}, Universit\'e de Provence,
Pr\'epublication LATP, {\bf 12}, 2001, 27~pp.

\bibitem[GM2001b]{gm2001b}
{\sc Gaussier}, H.; {\sc Merker}, J.:
{\em A new example of uniformly Levi degenerate hypersurface
in $\mathbb C^3$}, Ark. Mat. (to appear).

\bibitem[GM2001c]{gm2001c}
{\sc Gaussier}, H.; {\sc Merker}, J.:
{\em Symmetries of differential equations and infinitesimal
CR automorphisms of real analytic CR submanifolds of $\C^n$},
manuscript, 34~pp.

\bibitem[HJ1998]{hj1998}
{\sc Huang}, X.; {\sc Ji}, S.:
{\em Global holomorphic extension of a local map and
a Riemann mapping theorem for algebraic domains},
Math. Res. Lett. {\bf 5} (1998), no.1-2,
247--260.

\bibitem[HJY2001]{hjy2001}
{\sc Huang}, X.; {\sc Ji}, S.; {\sc Yau}, S.T.:
{\em An example of a real analytic strongly pseudoconvex
hypersurface which is not holomorphically equivalent to
any algebraic hypersurface},
Ark. Mat. {\bf 39} (2001), no.1, 75--93.

\bibitem[Lie1880]{lie1880}
{\sc Lie}, S.:
{\em Theorie der Transformationsgruppen},
Math. Ann. {\bf 16} (1880), 441--528.

\bibitem[Me1998]{me1998}
{\sc Merker}, J.:
{\em Vector field construction of Segre sets}, 
Preprint 1998, augmented in 2000. Downloadable at
{\tt arXiv.org/abs/math.CV/9901010}.

\bibitem[Me2001]{me2001}
{\sc Merker}, J.:
{\em On the partial algebraicity of holomorphic mappings
between two real algebraic sets},
Bull. Soc. Math. France {\bf 129} (2001), no.3, 547--591.

\bibitem[MW1983]{mw1983}
{\sc Moser}, J.K.; Webster S.M.:
{\em Normal forms for real surfaces in $\C^2$ near complex
tangents and hyperbolic surface transformations}, 
Acta Math. {\bf 150} (1983), no.3-4, 255--296.

\bibitem[Ol1986]{ol1986} 
{\sc Olver},~P.J.:
{\em Applications of Lie groups to differential equations}.
Springer Verlag, Heidelberg, 1986.

\bibitem[Pi1975]{p1975}
{\sc Pinchuk}, S.: 
{\em On the analytic
continuation of holomorphic mappings}
(Russian), Mat. Sb. (N.S.)
{\bf 98(140)} (1975)
no.3(11), 375--392, 416--435, 495--496.

\bibitem[Pi1978]{p1978}
{\sc Pinchuk}, S.:
{\em Holomorphic mappings of real-analytic hypersurfaces} (Russian),
Mat. Sb. (N.S.) {\bf 105(147)} (1978), no. 4,
574--593, 640.

\bibitem[Se1931]{se1931}
{\sc Segre}, B.:
{\em Intorno al problema di Poincar\'e della 
rappresentazione pseudoconforme}, 
Rend. Acc. Lincei, VI, Ser. {\bf 13} (1931), 676--683.

\bibitem[Sha2000]{sha2000} 
{\sc Shafikov}, R.: 
{\em Analytic
continuation of germs of holomorphic mappings between 
real hypesurfaces in $\C^n$}, Michigan Math. J. {\bf 47}
(2000), no.1, 133--149.

\bibitem[Sha2002]{sha2002} 
{\sc Shafikov}, R.: 
{\em Analytic
continuation of holomorphic correspondences and equivalence of domains in $\C^n$}, to appear in Inventiones Math.

\bibitem[SS1996]{ss1996}
{\sc Sharipov}, R.; {\sc Sukhov}, A.:
{\em On CR mappings between algebraic Cauchy-Riemann manifolds and
separate algebraicity for holomorphic functions}, 
Trans. Amer. Math. Soc. {\bf 348} (1996), no.2, 767--780.

\bibitem[St1991]{st1991}
{\sc Stanton}, N.:  
{\em Infinitesimal CR automorphisms of rigid 
hypersurfaces in $\C^2$},
J. Geom. Anal. {\bf 1} (1991), no.3, 231--267.

\bibitem[Sto2000]{sto2000}
{\sc Stormark}, O.:
{\em Lie's structural approach to PDE systems}, 
Encyclop{\ae}dia of mathematics and
its applications, vol. 80, 
Cambridge University Press, Cambridge, 2000, xv+572~pp.

\bibitem[Su2001a]{su2001a} 
{\sc Sukhov},~A.:
{\em Segre varieties and Lie symmetries},
Math. Z. {\bf 238} (2001), no.3, 483--492.

\bibitem[Su2001b]{su2001b}
{\sc Sukhov},~A.:
{\em On transformations of analytic CR structures}, 
Pub. Irma, Lille 2001, Vol. {\bf 56}, no. II.

\bibitem[Ve1999]{v1999}
{\sc Verma}, K:
{\em Boundary regularity of correspondences in $\C^2$}, 
Math. Z. {\bf 231} (1999), no.2, 253--299.

\bibitem[We1977]{w1977}
{\sc Webster}, S.M.:
{\em On the mapping problem for algebraic real hypersurfaces},
Invent. Math. {\bf 43} (1977), no.1, 53--68.

\bibitem[We1978]{w1978}
{\sc Webster}, S.M.:
{\em On the reflection principle in several complex variables}, 
Proc. Amer. Math. Soc. {\bf 71} (1978), no.1, 26--28.

\bibitem[Za1995]{za1995}
{\sc Zaitsev}, D.:
{\em On the automorphism groups of algebraic bounded domains}, 
Math. Ann. {\bf 302} (1995), no.1, 105--129.

\end{thebibliography}
